\documentclass[a4paper,12pt]{article}
\usepackage[T2A]{fontenc}
\usepackage[cp1251]{inputenc}
\usepackage[cp1251]{inputenc}
\usepackage[dvips]{graphicx}
\usepackage{amsfonts,amsmath,amssymb,amscd,array, euscript}

\begin{document}
\title{Algebras of operators in Banach spaces over the quaternion skew field
and the octonion algebra}
\author{
S.V. Ludkovsky}
\date{27.03.2005}
\maketitle
\par \vspace{0.5cm}
\par UDC 517.983 + 517.984

\par \vspace{0.5cm}
\section{Introduction}
The skew field of quaternions $\bf H$ and the algebra of octonions
$\bf O$ are algebras over the real field $\bf R$, but they are not
algebras over the complex field $\bf C$, since each embedding of
$\bf C$ into $\bf H$ or into $\bf O$ is not central. Therefore, a
research of algebras of operators over $\bf H$ or $\bf O$ can not be
reduced to algebras of operators over $\bf C$. On the other hand,
the developed below theory of algebras of operators over $\bf H$ and
$\bf O$ has many specific features in comparison with the general
theory of algebras of operators over $\bf R$ due to the graded
structures of algebras $\bf H$ and $\bf O$. Moreover, the algebra of
octonions can not be realized as the subalgebra of any algebra of
matrices over $\bf R$, since the algebra of octonions is not
associative. At the same time the skew product of the algebra of
matrices over $\bf C$, known as the particular case of the skew
product of algebras of operators, produces only algebras over $\bf
C$ and can not give the algebra $\bf H$ or $\bf O$.
\par The results of this work can be used also for the development
of noncommutative geometry, superanalysis, quantum mechanics over
$\bf H$ and $\bf O$, and the theory of representations of non
locally compact groups, for example, of the type of the group of
diffeomorphisms and the group of loops of quaternion or octonion
manifolds (see \cite{connes,oystaey,emch,lulgcm,lupm}). The wast
majority of preceding works on superanalysis was devoted to
supercommutative superalgebras of the type of Grassman algebra, but
for the noncommutative superalgebras it was almost undeveloped. The
skew field of quaternions and the algebra of octonions serve as very
important examples of superalgebras, which are not supercommutative,
moreover, the algebra of octonions is not associative. In the
present work there are used results of the preceding articles of the
author on this theme, in particular, noncommutative integral over
$\bf H$ and $\bf O$ \cite{luoyst,lufsqv,luoystoc} serving as the
analog of the integral of the Cauchy type known over $\bf C$. As
examples of quaternion or octonion unbounded operators serve
differential operators, including that of in the partial
derivatives. They arise naturally, for example, the
Klein-Gordon-Fock equation may be written in  the form $(\partial
^2/ \partial z^2+\partial ^2/ \partial {\tilde z}^2)f=0$ on the
space of quaternion locally $(z,{\tilde z})$-analytic functions $f$,
where $z$ is the quaternion variable, $\tilde z$ is the conjugated
variable, $z{\tilde z}=|z|^2$. The Dirack operator for the spin
systems over $\bf H^2$ can be written in the form ${ {0\quad
\partial / \partial z} \choose {{- \partial / \partial {\tilde z}}
\quad 0}} $, that is used in the theory of of spin manifolds
\cite{lawmich}, but each spin manifold can be embedded in the
quaternion manifold \cite{lufsqv}.
\par  The skew field of quaternions $\bf H$ has the automorphism $\eta
$ of the order two $\eta: z\mapsto {\tilde z}$, where ${\tilde z}=
w_1-w_ii-w_jj-w_kk$, $z=w_1+w_ii+w_jj+w_kk$; $w_1,...,w_k\in \bf R$.
There is the norm in $\bf H$ such that $|z|=|z{\tilde z}|^{1/2}$,
consequently, ${\tilde z}=|z|^2z^{-1}$.
\par The algebra $\bf K$ of octonions (octave, the Cayley algebra)
is defined as the eight dimensional algebra over $\bf R$ with he
basis, for example,
\par $(1)$ ${\bf b}_3:={\bf b}:= \{ 1, i, j, k, l, il, jl, kl \} $
such that
\par $(2)$ $i^2=j^2=k^2=l^2=-1$, $ij=k$, $ji=-k$,
$jk=i$, $kj=-i$, $ki=j$, $ik=-j$, $li=-il$, $jl=-lj$, $kl=-lk$;
\par $(3)$ $(\alpha +\beta l)(\gamma +\delta l)=(\alpha \gamma
-{\tilde {\delta }}\beta )+(\delta \alpha +\beta {\tilde {\gamma }})l$ \\
it is the law of multiplication in $\bf  K$ for each $\alpha $,
$\beta $, $\gamma $, $\delta \in \bf H$, $\xi :=\alpha +\beta l\in
\bf K$, $\eta :=\gamma +\delta l\in \bf K$, ${\tilde z}:=v-wi-xj-yk$
for quaternions $z=v+wi+xj+yk \in \bf H$ with $v, w, x, y\in \bf R$.
\par The algebra of octonions is neither commutative, nor associative,
since $(ij)l=kl$, $i(jl)=-kl$, but it is distributive and ${\bf R}1$
is its centre. If $\xi :=\alpha +\beta l\in \bf K$, then
\par $(4)$ ${\tilde {\xi }}:={\tilde {\alpha }}-\beta l$ is called the adjoint element
for $\xi $, where $\alpha , \beta \in \bf H$. Then
\par $(5)$ $(\xi \eta )^{\tilde .}={\tilde {\eta }}{\tilde {\xi }}$,
${\tilde {\xi }} + {\tilde {\eta }}= (\xi +\eta )^{\tilde .}$ and
$\xi {\tilde {\xi }}=|\alpha |^2+|\beta |^2$, \\
where $|\alpha |^2=\alpha {\tilde {\alpha }}$, such that
\par $(6)$ $\xi {\tilde {\xi }}=:|\xi |^2$ and $|\xi |$
is the norm in $\bf K$. Therefore,
\par $(7)$ $|\xi \eta |=|\xi | |\eta |$, \\
consequently, $\bf K$ does not contain the divisors of zero (see
also \cite{kurosh}). The multiplication of octonions satisfies the
equations:
\par $(8)$ $(\xi \eta )\eta =\xi (\eta \eta )$,
\par $(9)$ $\xi (\xi \eta )=(\xi \xi )\eta $, \\
which define the property of alternativity of algebras. In
particular, $(\xi \xi )\xi =\xi (\xi \xi )$. Put ${\tilde {\xi
}}=2a-\xi $, where $a=Re (\xi ):=(\xi + {\tilde {\xi }})/2 \in \bf
R$. In view of the fact that ${\bf R}1$ is the centre in $\bf K$ and
${\tilde {\xi }}\xi =\xi {\tilde {\xi }}=|\xi |^2$, then from the
equations $(8,9)$ by induction it follows, that for each $\xi \in
\bf K$ and each $n$-times product, $n\in \bf N$, $\xi (\xi (... \xi
\xi )... )=(...(\xi \xi )\xi ...)\xi $ the result does not depend on
the order of brackets (the order of subsequent multiplications),
consequently, the definition of $\xi ^n:= \xi (\xi (... \xi \xi
)...)$ does not depend on the order of brackets. This also shows
that $\xi ^m\xi ^n=\xi ^n\xi ^m$, $\xi ^m{\tilde {\xi }}^m={\tilde
{\xi }}^m\xi ^n$ for each $n, m\in \bf N$ and $\xi \in \bf K$. Apart
from the quaternions, the octonion algebra can not be realized as a
subalgebra of the algebra ${\bf M}_8({\bf R})$ of all $8\times
8$-matrices over $\bf R$, since $\bf K$ is not associative, but
${\bf M}_8({\bf R})$ is associative. The noncommutative
nonassociative algebra of octonions $\bf K$ is $\bf Z_2$-graded $\bf
R$-algebra ${\bf K}={\bf K}_0+{\bf K}_1$, where elements ${\bf K}_0$
are $\underline {even}$ and elements in ${\bf K}_1$ are $\underline
{odd}$ (see, for example, \cite{kurosh}). There are natural
embeddings ${\bf C}\hookrightarrow \bf K$ and ${\bf
H}\hookrightarrow \bf K$, but neither $\bf K$ over $\bf C$, nor $\bf
K$ over $\bf H$, nor $\bf H$ over $\bf C$ are algebras, since their
centres are the following: $Z({\bf H})=Z({\bf K})=\bf R$.
\par In this article there are given the principal features of quaternion and
octonion cases, since in one article it is impossible to give such
broad theory over $\bf H$ or $\bf O$, as the well developed theory
of operators over $\bf C$ \cite{danschw,kadring}.
\par In the second section the are investigated unbounded, as well as bounded
quasilinear operators in the Hilbert spaces $X$ over algebras $\bf
K$. At the same time the analog of the scalar product in $X$ is
defined with values in the algebra $\bf K$. There are defined and
investigated also graded operators of projections and graded
projection valued measures. The linearity of operators over the
algebra $\bf K$ is already not worthwhile because of the
noncommutativity of $\bf K$, therefore there is introduced the
notion of quasilinear operators. Further theorems about spectral
representations of projection valued graded measures of normal
quasilinear operators, which can be unbounded, are proved. At the
same time graded projection valued measures in the general case may
be noncommutative and nonassociative. Because of the
noncommutativity of the skew field of quaternions and the algebra of
octonions commutative algebras $\cal Z$ over them withdraw their
importance, since they can be only trivial ${\cal Z}{\bf K}={\cal
Z}= \{ 0 \} $, moreover, over the algebra of octonions associativity
loose its importance. Therefore, in the second section there is
introduced the notion of quasicommutative algebra over $\bf K$ and
for them it is proved the analog of the Gelfand-Mazur theorem.
Moreover, there are considered spectral decompositions of the
unitary and self-adjoint operators, it is proved the analog of the
Stone theorem over the skew field of quaternions and the algebra of
octonions for the one parameter families of unitary operators.
\par In the third section the are investigated  $C^*$-algebras over $\bf K$.
For them there are proved analogs of the theorems of
Gelfand-Naimark-Segal, von Neuman, Kaplansky  and so on. There is
studied the question about topological and algebraic irreducibility
of the action of a $C^*$-algebra of quasilinear operators in the
Hilbert space over $\bf K$.

\section{Theory of unbounded operators}
{\bf 1. Definitions and Notes.} Let $X$ be a vector space over the
skew field of quaternions $\bf H$ or the algebra of octonions $\bf
O$, that is, $X$ is the additive group, and the multiplication of
the vectors $v\in X$ on scalars $a, b\in \bf H$ from the left and
the right satisfies the axioms of the associativity and
distributivity, but in the case of the octonion algebra  $\bf O$
associativity is replaced on the alternativity preserving the
ditributivity, where the left distributivity means, that
$(a+b)v=av+bv$ and $a(v+w)=av+aw$ for each $a, b\in \bf K$, $v, w
\in X$, ${\bf K}=\bf H$ or ${\bf K}=\bf O$, the right distributivity
means that $v(a+b)=va+vb$ and $(v+w)a=va+wa$ for each $a, b\in \bf
K$, $v, w \in X$, associativity from the left means that
$a(bv)=(ab)v$ for each $a, b \in \bf H$ and $v\in X$, the
associativity from the right means that $(vb)a=v(ba)$ for each $a, b
\in \bf H$ and $v\in X$, the alternativity from the left means that
$b(bv)=(b^2)v$ for each $b \in \bf O$ and $v\in X$, the
alternativity from the right means that $(vb)b=v(b^2)$ for each $b
\in \bf K$ and $v\in X$. \par If there exists the norm $\| v
\|_X=:\| v \| $ in $X$ relative to which $X$ is complete, where $\|
av \| =|a|_{\bf K}\| v \| $, $\| vb \| =|b|_{\bf K}\| v \| $, $\|
v+w\| \le \| v \| +\| w \| $ for each $v, w\in X$, $a, b\in \bf K$,
then $X$ is called the Banach space over $\bf K$, where ${\bf K}=\bf
H$ or ${\bf K}=\bf O$. Then $X$ has also the structure of the Banach
space $X_{\bf R}$ over $\bf R$, since $\bf K$ is the algebra over
$\bf R$ of the dimension $4$ for ${\bf K}=\bf H$ and $8$ for ${\bf
K}=\bf O$.
\par An operator $T$ on the dense vector subspace ${\cal D}(T)$
in $X$ with values in the Banach space $Y$ over $\bf K$ is called
($\bf K$-)right linear, if $T(va)=(T(v))a$, $T(v+w)=(T(v))+(T(w))$
for each $a\in \bf K$ and every $v, w\in {\cal D}(T)$ and in
addition $T$ is $\bf R$-linear on ${\cal D}(T)_{\bf R}$. For the
right linear operator we also write $Tv$ instead of $T(v)$. If $T$
is $\bf R$-linear and $T(av)=a(T(v))$, $T(v+w)= (T(v))+(T(w))$ for
each $a\in \bf K$, then $T$ is called ($\bf K$-)left linear. For the
left linear operator we also write $vT$ instead of $T(v)$.
\par An operator $T$ is called ($\bf K$-)linear, if it is left and right
linear operator simultaneously. (Nevertheless, $av$ in the general
case is not equal $va$ for $v\in X$ and $a\in \bf K$, also in the
general case $cz$ is not equal $zc$ for $z\in Y$ and $c\in \bf K$,
therefore, the linearity of an operator over $\bf K$ differs from
the linearity of an operator over commutative fields $\bf R$ or $\bf
C$, that is, it is better to avoid such terminology of the linearity
over $\bf K$).
\par An operator $T: {\cal D}(T)\to {\cal R}(T)$ we call ($\bf
K$-)quasilinear, if it is additive $T(v+w)=T(v)+T(w)$ and $\bf
R$-homogeneous $T(av)=aT(v)=T(va)$ for each $a\in \bf R$, $v$ and
$w\in X$, where ${\cal R}(T)\subset Y$ means the region of values of
the operator. For example, derivatives of the quaternion holomorphic
functions are $\bf H$-quasilinear (see \cite{luoyst}).
\par Let  $L_q(X,Y)$ ($L_r(X,Y)$; $L_l(X,Y)$) denote Banach space of all
bounded quasilinear operators $T$ from $X$ into $Y$ (right or left
linear operators correspondingly), $\| T\| :=\sup_{v\ne 0} \| Tv \|
/\| v \| $; $L_q(X) :=L_q(X,X)$, $L_r(X):=L_r(X,X)$ and
$L_l(X):=L_l(X,X)$. The resolvent set $\rho (T)$ of a quasilinear
operator $T$ is defined as $$\rho (T):=\{ z\in {\bf K}: \mbox{ there
exists } R(z;T) \in L_q(X,{\cal D}(T)) \} ,$$ where
$R(z;T):=R_z(T):=(zI-T)^{-1}$, $I$ deontes the unit operator on $X$,
$Iv=v$ for each $v\in X$, analogously for right or leftlinear
operators, where ${\cal D}(T)$ is dense in $X$. The spectrum is
defined as $\sigma (T):={\bf K}\setminus \rho (T) $.
\par {\bf 2. Lemma.} {\it For each $z_1$ and $z_2\in \rho (T)$: \\
$(i)\quad R(z_2;T)-R(z_1;T)=R(z_1;T)((z_1-z_2)R(z_2;T))=
(R(z_1;T)(z_1-z_2))R(z_2;T)$.}
\par {\bf Proof.} The algebra $\bf O$ is alternative, therefore,
$a^{-1}(av)=v$ for each $a\in \bf O$ with $a\ne 0$ and $v\in X$,
since $a^{-1}={\tilde a}/|a|$, $|a|\in \bf R$, $Z({\bf O})=\bf R$.
If an operator $A: X\to X$ is invertible, where $A\in L_q(X,X)$,
then there exists $A^{-1}\in L_q(A(X),X)$. In particular, if
$A(X)=X$ for an invertible operator $A$, then $A(X)$ and $X$ are
isomorphic as the $\bf K$-vector spaces. In view of the fact that
the vector span over $\bf K$ for $A(X)$ is the $\bf K$-vector
subspace in $X$, then $A^{-1}(Av)=v$ for each $v\in X$, in
particular, for $v=R(z,T)u$, where $u\in X$, $z\in \rho (T)$.
Therefore, the following identities are satisfied:
$R(z_1;T)((z_1I-z_2I)R(z_2;T))=(z_1I-T)^{-1}((z_1I-T-(z_2I-T_2))(z_2I-T)^{-1})=(z_2I-T)^{-1}
-(z_1I-T)^{-1}$, since the multiplication in the skew field of
quaternions $\bf H$ is associative, and the algebra $\bf O$ is
alternative.
\par {\bf 3. Lemma.} {\it If $T$ is the closed quasilinear operator, then $\rho (T)$
is open in $\bf K$ and $R(z;T)$ is $\bf K$-holomorphic on $\rho
(T)$, where ${\bf K}=\bf H$ or ${\bf K}=\bf O$.}
\par {\bf Proof.} Let $z_0\in \rho (T)$, then the
operator $R(z_0;T)$ is closed and by the closed graph theorem it is
bounded, since it is defined everywhere. If $z\in \bf K$ and
$|z_0-z|< \| (z_0I-T)^{-1} \|^{-1}$, then
$(zI-T)=(z_0I-T)(I-R(z_0;T)(z_0-z))$, moreover, $$(i)\quad R(z;T)=
\{ \sum_{n=0}^{\infty }[R(z_0;T)(z_0-z)]^n \} R(z_0;T)\in L_s(X),$$
since this series $[R(z_0;T)(z_0-z)]^n R(z_0;T)$ converges relative
to the norm topology in $L_s(X)$. In view of $(i)$ $R(z;T)$ $\bf K$
is locally $z$-analytic, consequently, holomorphic on $\rho (T)$ in
accordance with theorem 2.16 \cite{luoyst,luoystoc}.
\par {\bf 4. Notes and Definitions.} For a normed (or Banach)
space $X$ over $\bf K$ its topologically right adjoint space $X^*_r$
is defined as consisting of all functionals $f: X\to \bf K$ such
that $f$ is $\bf R$-linear and $\bf K$-right linear, where ${\bf
K}=\bf H$ or ${\bf K}=\bf O$. Analogously we put $X^*_q:=L_q(X,{\bf
K})$ and $X^*_l:=L_l(X,{\bf K})$, where $X^*_q$ is the topologically
quasi adjoint space, $X^*_l$ is the topologically left adjoint
space. Then $X^*_s$ is the Banach space over $\bf K$ with the norm
$\| f \| :=\sup_{x\ne 0}|fx|/ \| x \|$. If $X$ and $Y$ are normed
(or Banach) spaces over $\bf K$ and $T: X\to Y$ belongs to
$L_s(X,Y)$, then $T^*$ is defiend by the equation:
$(T^*y^*)(x):=y^*\circ T(x)$ for each $y^*\in Y^*_s$, $x\in X$. Then
$T^*\in L_l(Y^*_r,X^*_r)$ for each $T\in L_r(X,Y)$. If $T\in
L_l(X,Y)$, then $T^*\in L_r(Y^*_l,X^*_l)$, since $y^*\circ T(x)=
xT^*y^*$ in the symmetrical notations, where $x\in X$. If $T\in
L_q(X,Y)$, then $T^*\in L_q(Y^*_q,X^*_q)$.
\par Let $\hat X$ and $\hat Y$ be images relative to the natural embeddings
of $X$ and $Y$ in $(X^*_q)^*_q$ and $(Y^*_q)^*_q$ respectively. For
each $T\in L_s(X,Y)$ define ${\hat T}\in L_s({\hat X},{\hat Y})$ by
the equation ${\hat T}({\hat x})=\hat y$, where $y=T(x)$. Each
function $S$ defined on a domain $(X^*_q)^*_q\supset dom (S)\supset
\hat X$ and such that $S({\hat x})=T({\hat x})$ for each ${\hat
x}\in \hat X$ is called the extension of $T$.
\par  Put by the induction
${\cal D}(T^n):= \{ x: x\in {\cal D}(T^{n-1}), T^{n-1}(x)\in {\cal
D}(T) \} $, ${\cal D}(T^{\infty }):=\bigcap_{n=1}^{\infty }{\cal
D}(T^n)$, where $T^0:=I$, $T^n(x):=T(T^{n-1}(x))$.
\par {\bf 4.1. Theorem.} {\it Let $Y$ be a $\bf K$-vector subspace
in a $\bf K$-vector normed space $X$, where ${\bf K}=\bf H$ or ${\bf
K}=\bf O$. Then to each $\bf K$-quasilinear functional $g\in Y^*_q$
there corresponds $f\in X^*_q$ such that $|f|=|g|$ and $f(y)=g(y)$
for each $y\in Y$.}
\par {\bf Proof.} If $X$ is a $\bf R$-vector
space, then the statement of this theorem follows from the
Hahn-Banach theorem for $p(x) := |g| |x|$. Remind, that the theorem
of Hahn-Banach  states, that if a real valued function $p$ on the
real linear space $X$ satisfies the conditions: $p(x+y)\le
p(x)+p(y)$, $p(ax)=ap(x)$ for each $x, y\in X$, $a\ge 0$, where $g$
it is a real valued linear functional on the $\bf R$-linear subspace
$Y$ in $X$ such that $g(x)\le p(x)$ for each $x\in Y$, then there
exists the real valued linear functional $f$ on the entire $X$ such
that $f(x)=g(x)$ for each $x\in Y$, moreover, $f(x)\le p(x)$ for
each $x\in X$.
\par Consider now the general case of the normed $\bf
K$-vector space $X$. Then the quasilinear functional $g$ has
the decomposition $g(y)=g_0(y)+g_1(y)i_1+...+g_m(y)i_m$ for each $y\in
Y$, where $g_0(y),...,g_m(y)$ are real valued functionals
on $Y$, $ \{ i_0,...,i_m \} $ denotes the set of the standard generators
of the algebra $\bf K$, $i_0=1$, $m=3$ for $\bf H$, $m=7$ for $\bf O$.
Then for the real numbers $a, b\in \bf R$ and arbitrary vectors
$x, y\in Y$ we get: $g_j(ax+by) = ag_j(x)+bg_j(y)$, $|g_j(y)|\le
|g(y)| \le |g| |y|$. There exists the decomposition of the space $X$ of the form
$X=X_0\oplus X_1i_1\oplus ... \oplus X_mi_m$, where $X_0, ...,X_m$ are pairwise
isomorphic real linear normed spaces. Then
from the beginning of the proof of this theorem it follows, that
real valued linear functionals $f_j$ on $X$ such that
$|f_j|\le |g_j|$ for each $j=0,1,...,m$, $f_j(y)=g_j(y)$ for
each $y\in Y$. Then the quasilinear functional $f$ on $X$ it is possible
to prescribe by the formula: $f(x) := f_0(x)+f_1(x)i_1+...+f_m(x)i_m$. Let
$f(x)=r \exp (Mt)$, where $t\in \bf R$, $0<r<\infty $, $M\in {\bf
K}$, $Re (M)=0$ (see the polar decomposition of the elements $z\in \bf K$ in
the Corollary 3.6 in \cite{luoystoc}). Then $|f(x)| = \exp (-tM) f(x)\le
|g| |x|$, since the algebra $\bf K$ is alternative, consequently,
$|f|\le |g|$. On the other hand, since $f$ is the extension of $g$, then
$|f|\ge |g|$, consequently, $|f|=|g|$.
\par {\bf 4.2. Lemma.} {\it Let $Y$ be a $\bf K$-vector subspace
of the $\bf K$-vector normed space $X$, where ${\bf K}=\bf H$ or
${\bf K}=\bf O$, moreover, $\inf_{y\in Y} |y-x| = d>0$. Then there
exists $f\in X^*_q$ such that $f(x)=1$, $|f|=1/d$, $f(y)=0$ for
each $y\in Y$.}
\par {\bf Proof.}  If $x\notin Y$, then ${\bf
K}x=x{\bf K}$ in view of the alternativity of the algebra $\bf K$,
consequently, $({\bf K}x)\cap Y= \{ 0 \} $, since $Y$ is the $\bf
K$-vector space. Then each vector $z$ from the $\bf K$-vector span
$Z:=span_{\bf K} (Y, \{ x \} )$ is uniquely represented as:
$z=y+ax$, where $a\in \bf K$, $y\in Y$. For $z=y+ax$ it is possible
to give $f(z):=a$, then $f$ is the $\bf K$-quasilinear functional on
$Z$. If $a\ne 0$, then $|z|=|y+ax|=|a| |a^{-1}y+x|\ge |a| d$,
consequently, $|f(z)| \le |z| /d$, $|f|\le 1/d$. Let $y_n\in Y$ be a
sequence of vectors such that $\lim_{n\to \infty } |x-y_n|=d$,
consequently, $1=f(x-y_n)\le |f| |x-y_n|$, Therefore, $d^{-1}\le
|f|$, then $|f|=1/d$. In view of Theorem 4.1 there exists the
extension $\bf K$-quasilinear functional $f$ on $X$.
\par {\bf 4.3. Corollary.} {\it Let $x$ be a vector, not belonging to
the closed $\bf K$-vector subpace $Y$ of a $\bf K$-vector
normed space $X$, where ${\bf K}=\bf H$ or ${\bf K}=\bf O$. Then
there exists $f\in X^*_q$ such that $f(x)=1$, $f(y)=0$ for each
$y\in Y$.}
\par {\bf 4.4. Corollary.} {\it For each $x\ne 0$ in a $\bf K$-vector
normed space $X$, where ${\bf K}=\bf H$ or ${\bf K}=\bf O$, there
exists $g\in X^*_q$ such that $|g|=1$, $g(x)=|x|$.}
\par {\bf Proof.} The application of Lemma 4.2 for $Y= \{ 0 \} $
gives the $\bf K$-quasilinear functional $g(y):=|x|f(y)$ for each
$y\in X$, where $f\in X^*_q$ with $f(x)=1$.
\par {\bf 4.5. Corollary.} {\it For each $x\in X$ in a $\bf
K$-vector normed space $X$, where ${\bf K}=\bf H$ or ${\bf
K}=\bf O$, there exists $|x|=\sup_{f\in S^*_q} |f(x)|$, where
$S^*_q$ is the closed sphere of the unit radius with the centre at zero
$X^*_q$. If $f(x)=f(y)$ for each $f\in X^*_q$, then $x=y$.}
\par {\bf 4.6. Lemma.} {\it A mapping $T\mapsto T^*$ is an isometric isomorphism
 of the space $L_q(X,Y)$ into $L_q(Y^*_q,X^*_q)$, where $X$ and $Y$
are normed spaces over ${\bf K}=\bf H$ or ${\bf K}=\bf O$.}
\par {\bf Proof.} If a $\bf K$-quasilinear functional
$f\circ T$ is continuous, then $T^*\circ f^*\in X^*_q$, moreover,
the mapping $T\mapsto T^*$ is $\bf K$-quasilinear. Due to Corollary
4.5 $|T(x)|= \sup_{|f|\le 1} |f(T(x))|$, consequently,
\par $|T^*|=\sup_{|f|\le 1} |T^*\circ f^*|=$  $\sup _{|f|\le 1,
|x|\le 1} |f(T(x))|=$ \par  $\sup_{|x|\le 1} \sup_{|f|\le 1}
|f(T(x))|=\sup_{|x|\le 1} |T(x)|=|T|$,\\  therefore, $T\mapsto T^*$
is the isometric isomorphism from $L_q(X,Y)$ into
$L_q(Y^*_q,X^*_q)$.
\par {\bf 4.7. Lemma.} {\it  Let $X$ and $Y$ be two $\bf K$-vector normed
spaces, where ${\bf K}=\bf H$ or ${\bf K}=\bf O$, then an operator
$T^*$ adjoint to $T\in L_q(X,Y)$ is a continuous $\bf K$-quasilinear
operator from $Y^*_q$ into $X^*_q$.}
\par {\bf 4.8. Lemma.} {\it Let $X$, $Y$, $Z$ be $\bf K$-vector normed
spaces, where ${\bf K}=\bf H$ or ${\bf K}=\bf O$. If $T\in
L_q(X,Y)$, $U\in L_q(Y,Z)$, then $(UT)^*= T^*U^*$, moreover, the
operator $I^*$ for $I\in L_q(X):=L_q(X,X)$ is the unit operator in
$L_q(X^*_q)$.}
\par {\bf Proof.} If $f\in Z^*_q$, $g\in X^*_q$, then
$((UT)^*\circ f)(x)=f(U(T(x)))$ $=(U^*f)(T(x))=$ $((T^*U^*)\circ
f)(x)$, also $I^*g=g\circ I=g$. Thus, the mapping $T\mapsto T^*$ of
the ring $L_q(X)$ into the ring $L_q(X^*)$ is the antiisomorphism.
\par {\bf 4.9. Lemma.} {\it If $T\in L_q(X,Y)$, where $X$ and $Y$ are
$\bf K$-vector normed spaces, ${\bf K}=\bf H$ or ${\bf K}=\bf O$,
then the repeated adjoint operator $T^{**}: (X^*_q)^*_q\to
(Y^*_q)^*_q$ is the extension of a $\bf K$-quasilinear operator $T$
.}
\par {\bf Proof.} Let $x\in X$, $g\in Y^*_q$, then
$(T^{**}{\hat x})g={\hat x}(T^*g)=(T^*g)(x)=g(T(x))=({\hat T}{\hat
x})(g)$.
\par {\bf 4.10. Lemma.} {\it Let $X$ and $Y$ be two $\bf K$-vector
normed spaces over ${\bf K}=\bf H$ or ${\bf K}=\bf O$. Then the
following two conditions are equivalent:
\par $(1)$ there exists the inverse operator $T^{-1}\in L_q(Y,X)$ and
$T^{-1}$ is defined on the entire $Y$; \par $(2)$ an adjoint
operator $T^*$ has an inverse operator $(T^*)^{-1}\in
L_q(X^*_q,Y^*_q)$ and $(T^*)^{-1}$ is defined on the entire
$X^*_q$.}
\par {\bf Proof.} Let $T^{-1}\in L_q(Y,X)$, then due to
Lemma 4.8 there are satisfied the equations $(TT^{-1})^* =
(T^{-1})^*T^* = I$ on $Y^*_q$, and $(T^{-1}T)^*=T^*(T^{-1})^*$ is
the unit operator on $X^*_q$. Thus, there exists $(T^*)^{-1}\in
L_q(X^*_q,Y^*_q)$ and $(T^*)^{-1}=(T^{-1})^*$.
\par Vice versa, if there exists $(T^*)^{-1}\in L_q(X^*_q,Y^*_q)$, then
due to the proof just given above there exists $(T^{**})^{-1}\in
L_q((Y^*_q)^*_q,(X^*_q)^*_q)$. Thus, the mapping $T^{**}$ is the
homeomorphism. For $T\in L_q(X,Y)$ we have
\par $((b_1f_1+b_2f_2)\circ T)(x)=((b_1f_1)\circ T)(x)+((b_2f_2)\circ
T)(x)$ $b_1(f_1(T(x))+b_2(f_2(T(x))$ \par
$=(T^*(b_1f_1))(x)+(T^*(b_2f_2))(x)$
$=b_1(T^*(f_1)(x))+b_2(T^*(f_2(x))$ \\ for each $x\in X$ and each
$b_1, b_2\in \bf K$, $f_1, f_2\in Y^*$. Since $T^*$ has the inverse
operator $(T^*)^{-1}$ defined on the entire $X^*_q$, then $T(X)$ is
the $\bf K$-vector space. Due to Lemma 4.9 $T^{**}$ serves as the
extension of the operator $T$, consequently, $T$ is bijective and
$T(X)$ is closed. If $y\in Y$ and $y\notin T(X)$, then due to
Corollary 4.3 there exists $g\in Y^*_q$ such that $g\ne 0$,
$gT=T^*g=0$. But this contradicts to the bijectivity of the operator
$T^*$, that is, $T(X)=Y$.
\par {\bf 5. Lemma.} {\it Let $T\in L_s(X)$, then
$\sigma (T^*)=(\sigma (T))^{\tilde .}$ and $(R(\lambda ,T))^*=
R(\lambda ^*,T^*)$, where $\lambda ^*I:=(\lambda I)^*$, $s\in \{ q,
r, l \} $.}
\par {\bf Proof.} If $S\in L_s(X,Y)$ and there exists
$S^{-1}\in L_s(Y,X)$, then $S^* \in L_u(Y^*_s,X^*_s)$ has the
inverse operator $(S^*)^{-1} \in L_u(X^*_s,Y^*_s)$ and $(S^{-1})^*=
(S^*)^{-1}$, where $(s,u)\in \{ (q,q); (r,l); (l,r) \} $. Then
\par $(\lambda I-T)^*y^*=y^*\circ (\lambda I-T)$ $=y^*\circ (\lambda
I)-y^*\circ T$, consequently, \par $(\lambda ^*I-T^*)[(\lambda
I-T)^*]^{-1}=I$ and $(R(\lambda ,T))^* =R(\lambda ^*,T^*)$.
\par {\bf 6. Definition and Note.} Denote by ${\cal H}(T)$
the family of all holomorphic functions $f$ in the variable $z\in
\bf K$ on neighborhoods $V_f$ for $\sigma (T)$, where $T\in L_s(X)$,
$s\in \{ q, r, l \} $, ${\bf K}=\bf H$ or ${\bf K}=\bf O$, and for a
quasilinear operator $T$ let ${\cal H}_{\infty }(T)$ be the family
of all holomorphic functions on neighborhoods $U_f$, $U_f\subset \bf
K$,  for $\sigma (T)$ and $\infty $ in the one point
compactification $\bf {\hat K}$ of the algebra $\bf K$ as the
locally compact topological space. Choose a marked point $z_0\in
\sigma (T)$. For each $M = x_0i_0 + ... + x_mi_m\in \bf K$ with
$|M|=1$ and $Re (M)=0$, where $x_0, ...,x_m\in \bf R$, there exists
a closed rectifiable path $\eta $ consisting of the finite union of
paths $\eta (s)=z_0+r_p\exp (2\pi sM)$ with $s\in [a_p,b_p]\subset
[0,1]\subset \bf R$ and segments of the straight lines $\{ z\in {\bf
K}: z=z_0+(r_pt+r_{p+1}(1-t)) \exp (2\pi b_pM), t\in [0,1] \} $
joining them, moreover, $\eta \subset U\setminus \sigma (T)$, where
$a_p<b_p$ and $0<r_p<\infty $ for each $p=1,...,m$, $m\in \bf N$,
$b_p=a_{p+1}$ for each $p=1,...,m-1$, $a_1=0$, $b_m=1$. Then there
exists a rectifiable closed path $\psi $ homotopic to $\eta $ and a
neighborhood $U$ satisfying conditions of Theorem $3.9$
\cite{luoyst} and such that $\psi \subset U\setminus \sigma (T)$.
For $T\in L_q(X)$ it is possible to define  $$(i)\quad f(T):=(2\pi
)^{-1} (\int_{\psi }f(\zeta )R(\zeta ;T)d\zeta )M^{-1},$$ where the
convergence is defined relative to the weak operator topology. This
integral depends on $f, T$ and it is independent from $U$, $\psi $,
$\eta $, $\gamma $, $M$. For an unbounded quasilinear operator $T$
let $A:=-R(a;T)$ and $\Psi : {\bf {\hat K}}\to {\bf {\hat K}}$,
$\Psi (z):=(z-a)^{-1}$, $\Psi (\infty )=0$, $\Psi (a)=\infty $,
where $a\in \rho (T)$. For $f\in {\cal H}_{\infty }(T)$ we define
$f(T):=\phi (A)$, where $\phi \in {\cal H}_{\infty }(A)$ is given by
the equation $\phi (z):=f(\Psi ^{-1}(z))$.
\par {\bf 7. Note.} Consider  the Banach space $X$ over the skew field of quaternions
$\bf H$ as the Banach space $X_{\bf C}$ over $\bf C$, then
\par $(i)$ $X_{\bf C}= X_1\oplus X_2j$,\\
where $X_1$ and $X_2$ are Banach spaces over $\bf C$ such that $X_1$
is isomorphic with $X_2$. The complex conjugation in $\bf C$ induces
the complex conjugation of vectors in $X_m$, where $m=1$ or $m=2$.
Each vector $x\in X$ we can write in the matrix form
\par $(ii)$ $x={{x_1\quad x_2}\choose {-{\bar x}_2\quad {\bar x}_1}}$, \\
where $x_1\in X_1$ and $x_2\in X_2$. Then each quasilinear operator
$T$ we can write in the form
\par $(iii)$ $T={{T_1\quad T_2}\choose {-{\bar T}_2\quad {\bar T}_1}}$, \\
where $T_1: X_1\supset {\cal D}(T_1)\to Y_1$, $T_2: X_1\supset {\cal
D}(T_2)\to Y_2$, $T(x)=Tx$ for $s\in \{ q, r \} $, $T(x)=xT$ for
$s=l$.
\par $(iv)$ ${\bar T}_mx:=\overline {T_m{\bar x}}$, where $m=1$ or $m=2$.
\par In particular, for the commutator  $[\zeta I,T]$
for $\zeta ={{b\quad 0}\choose {0\quad \bar b}}\in \bf H$, where
$b\in \bf C$, the following formula is accomplished:
\par $(v)$ $[\zeta I,T]=2(-1)^{1/2} Im(b) {{0\quad T_2}\choose
{-{\bar T}_2 \quad 0}}$, where $Im (b)$ denotes the imaginary part of $b$, \\
$2(-1)^{1/2} Im(b)=(b-{\bar b})$.
\par {\bf 8. Theorem.} {\it If $f\in {\cal H}_{\infty }(T)$ over the algebra
${\bf K}=\bf H$ or ${\bf K}=\bf O$, then $f(T)$ does not depend on
$a\in \rho (T)$ and    $$(i)\quad f(T)=f(\infty )I+(2\pi )^{-1}
(\int_{\psi } f(\lambda ) R(\lambda ;T) d\lambda ) M^{-1}.$$}
\par {\bf Proof.} If $a\in \rho (T)$, then $0\ne b=(\lambda
-a)^{-1}$ for $\lambda \ne a$ and \par $(T-aI)(T-\lambda I)^{-1}=
(bI-A)^{-1}b$,  therefore, \par $I+b^{-1}(T-\lambda
I)^{-1}=b(bI-A)^{-1}b-bI$ \\
and $b\in \rho (A)$. If $0\ne b\in \rho
(A)$, then $A(bI-A)^{-1}= (T-\lambda I)^{-1}b^{-1}$ and $\lambda \in
\rho (T)$. The point $b=0$ is in $\sigma (A)$, since $A^{-1}=T-aI$
is unbounded. Let $a\notin V$, then $U=\psi ^{-1}(V)\supset \sigma
(A)$ and $U$ is open in $\bf \hat K$, but $\phi (z):=f(\psi
^{-1}(z)) \in {\cal H}_{\infty }(U)$. Due to corollary 3.26 in
\cite{luoyst} and 3.25 in \cite{luoystoc} we get the formula $(i)$.
\par {\bf 9. Theorem.} {\it Let $a, b, c, e\in \bf K$,
$f, g \in {\cal H}_{\infty }(T)$, where ${\bf K}=\bf H$ or ${\bf
K}=\bf O$. Then
\par $(i)$ $a(fc)+b(ge)\in {\cal H}_{\infty }(T)$ and $a(fc)(T)+b(ge)(T)=
(a(fc)+b(ge))(T)$;
\par $(ii)$ $fg\in {\cal H}_{\infty }(T)$ and $f(T)g(T)=(fg)(T)$;
\par $(iii)$ If $f$ is decomposable into a converging series
$$f(z)=\sum_k \{ (b_k,z^k) \} _{q(2m(k))}$$ in a neighborhood $\sigma
(T)$, then $$f(T)=\sum_k \{ (b_k,T^k) \} _{q(2m(k))}$$ on ${\cal
D}(T^{\infty })$, where $b_k=(b_{k,1},...,b_{k,m(k)})$, \\
$\{ (b_k,z^k) \} _{q(2m(k))} := \{
b_{k,1}z^{k_1}...b_{k,m(k)}z^{k_{m(k)}} \} _{q(2m(k))}$, $b_{k,j}\in
\bf K$;
\par $(iv)$ $f\in {\cal H}_{\infty }(T^*)$ and $f^*(T^*)=(f(T))^*$, where
$f^*(z):=(f(z^*))^*$, $q(v)$ is a vector indicating on an order of
the product in $ \{ * \} $ in the nonassociative case of ${\bf
K}=\bf O$, as in \S 2.1 \cite{luoystoc}, which can be omitted in the
associative case of ${\bf K}=\bf H$.}
\par {\bf Proof.} $\bf (i).$ Take $V_f\cap V_g=:V$
and for it construct $U$, $\eta $ and $\psi $ as in \S 2.6. Then the
first statement follows from $2.6.(i)$.
\par $\bf (ii).$ Due to Theorem $3.28$ \cite{luoyst}
or Theorem $3.27$ \cite{luoystoc} function ${\tilde g}({\tilde
z})=:\phi (z)$ belongs to ${\cal H}_{\infty }(T)$ and ${\tilde \phi
}({\tilde z})=g(z)$, where ${\tilde z}=x_0i_0-x_1i_1-...-x_mi_m$,
$z=x_0i_0+x_1i_1+...+x_mi_m$, $x_0, x_1,...,x_m \in \bf R$, $z\in
V_g\subset \bf K$. Using $\phi (A)$, we consider the case of a
bounded operator $T$. Due to $2.6.(i)$: \par $({\tilde y}^*(\phi
({\tilde T}) {\tilde h}))^{~}=(2\pi )^{-1}y^* (M\int_{\psi
}(d{\tilde \zeta } R({\tilde \zeta },T)g({\tilde \zeta }))h)$, \\
where ${\tilde y}^*({\tilde T}{\tilde h}):=(y^*(Th))^{\tilde .}$ and
${\tilde y}^*{\tilde h}:=(y^*h)^{\tilde .}$ for each $y^*\in X^*_q$
and $h\in X$. Therefore, $$({\tilde y}^*(\phi ({\tilde T}){\tilde
h}))^{\tilde .}=(2\pi )^{-1}y^* (M\int_{\tilde \psi }(d\zeta R(\zeta
,T)g(\zeta ))h),\mbox{ consequently},$$  $$(\phi ({\tilde
T}))^{\tilde .}=(2\pi )^{-1}M (\int_{\tilde \psi }d\zeta R(\zeta
,T)g(\zeta ))=g(T),$$ since left and right integrals coincide in the
space of holomorphic functions over $\bf K$. The function $fg$ is
holomorphic on $V$ (see \S \S $2.1$ and $2.12$
\cite{luoyst,luoystoc}). There exist $\psi _f$ and $\psi _g$ as in
\S $2.6$ and contained in $U\setminus \sigma (T)$, where ${\bar
U}\subset V$. Due to the Fubini theorem there exists
$$(v)\quad f(T)g(T)=(2\pi )^{-2}\int_{\psi _f}\int_{\psi _g}
(f(\zeta _1)R(\zeta _1;T)(d\zeta _1))((d\zeta _2)R(\zeta
_2,T)g(\zeta _2)),$$ where $\zeta _1\in \psi _f$ and $\zeta _2\in
\psi _g$. There are satisfied the identities $R(\zeta ;T)d\zeta
=d_{\zeta } Ln (\zeta I-T)$ and $(d\zeta ) R(\zeta ;T)=d_{\zeta} Ln
(\zeta I-T)$ for a chosen branch of $Ln$ (see \S \S $3.7, 3.8$
\cite{luoyst,luoystoc}), consequently,\par  $(R(\zeta _1;T)(d\zeta
_1))((d\zeta _2)R(\zeta _2;T))=$ $d_{\zeta _1}(d_{\zeta _2} (Ln
(\zeta _1I-T) Ln (\zeta _2I-T)))$ \par $=((d\zeta _1)R(\zeta
_1;T))(R(\zeta _2;T)d\zeta _2)$. Due to Lemma $2.2$:
$$(vi)\quad R(a;T)R(b;T)=[R(a;T)-R(b;T)](b-a)^{-1}$$
$$ +R(a;T)([R(b;T),(b-a)I](b-a)^{-1}),$$
$$(vii)\quad [R(b;T),(b-a)I]=R(b;T)([T,(b-a)I]R(b;T)),$$
since $[(bI-T),(b-a)I]=-[T,(b-a)I]$, where $a, b\in \rho (T)$. Let,
in particular, $\psi _f$ and $\psi _g$ be contained in the plane
${\bf R}\oplus i\bf R$ in $\bf H$, where $i, j, k$ are the
generators of $\bf H$ such that $i^2=j^2=k^2=-1$, $ij=k$, $jk=i$,
$ki=j$. In the case of ${\bf K}=\bf O$ we can embed the
corresponding copy of $\bf H$ into $\bf O$ an the integral over $\bf
O$ along a path then reduces to the integral over $\bf H$ along the
path. In addition, it is possible to use the invariance of the
integral under the corresponding conditions relative to the choice
of homotopic paths (see Theorems 2.11 and 3.9 \cite{luoystoc}). It
is necessary to note, that the operator $T$ is uniquely defined by
all its restrictions on subspaces $X_{\bf H}$ over $\bf H$ in $X$
under all possible embeddings of the skew field of quaternions $\bf
H$ as the subalgebra into the algebra of octonions $\bf O$, where
$X=X_0\oplus X_1l$, $X_0$ and $X_1$ are Banach spaces over $\bf H$
isomorphic $X_{\bf H}$, $l\in {\bf O}\setminus \bf H$ such that
${\bf H}\oplus {\bf H}l$ is the algebra isomorphic with $\bf O$,
$|l|=1$, $l$ is the generator of the doubling procedure of $\bf H$
to $\bf O$ up to an isomorphism of algebras. Indeed, for each $x\in
X$ there exists $X_{\bf H}$ with $x\in X_{\bf H}$. Due to $2.7.(v)$
and $2.8.(vii)$:
$$(viii)\quad \int_{\psi _f}\int_{\psi _g}f(\zeta _1)R(\zeta _1;T)
[R(\zeta _2;T),(\zeta _2-\zeta _1)I](\zeta _2-\zeta _1)^{-1}d\zeta
_1 d\zeta _2 g(\zeta _2)=0,$$ since a branch of the function $Ln$
can be chosen along the entire axis $j$ in $\bf H$, due to the
argument principle $3.30$ \cite{luoyst} it corresponds to the
residue $c(\zeta _1-z)^{-1}(\zeta _2-z)^{-1} b(\zeta _2-\zeta
_1)(\zeta _2-z)^{-1}(\zeta _2-\zeta _1)^{-1}$. In the case ${\bf
K}=\bf O$ in this equation we mean the restriction of $T$ on all
possible $X_{\bf H}$. Then from $(v,vi,viii)$ it follows:
$$(ix)\quad f(T)g(T)=(2\pi )^{-2}\int_{\psi _f}\int_{\psi _g}
f(\zeta _1)[R(\zeta _1;T)-R(\zeta _2,T)](\zeta _2-\zeta _1)^{-1}
d\zeta _1d\zeta _2 g(\zeta _2).$$ Choose $\psi _g $ such that $|\psi
_g(s)-z_0| > |\psi _f(s)-z_0|$ for each $s\in [0,1]$. From the
additivity of the integral along a path and the Fubini theorem:
$$(x) \quad f(T)g(T)=(2\pi )^{-2}\int_{\psi _f}
f(\zeta _1)R(\zeta _1;T)d\zeta _1(\int_{\psi _g}(\zeta _2-\zeta _1)^{-1}
d\zeta _2 g(\zeta _2))$$
$$-(2\pi )^{-2} \int_{\psi _g}
(\int_{\psi _f}f(\zeta _1)d\zeta _1 R(\zeta _2,T) (\zeta _2-\zeta
_1)^{-1}) d\zeta _2 g(\zeta _2).$$ Due to Theorems $3.9, 3.24$
\cite{luoyst} or $3.9, 3.28$ \cite{luoystoc} the second integral on
the right of $(x)$ is equal to zero, since $\int_{\psi _f}f(\zeta
_1)d\zeta _1 R(\zeta _2;T) (\zeta _2-\zeta _1)^{-1}=0$, but the
first integral gives
$$f(T)g(T)=(2\pi )^{-1} (\int_{\psi _f}
f(\zeta _1) R(\zeta _1;T) d\zeta _1 g(\zeta _1))M^{-1}$$
$$=(2\pi )^{-1} (\int_{\psi _f}
f(\zeta )g(\zeta ) d Ln (\zeta I-T)) M^{-1},$$ where $\zeta \in \psi
_f$.
\par $\bf (iii)$ follows from the application of $\bf (i,ii)$
by induction and the convergence of the series in the strong
operator topology.
\par $\bf (iv).$ Due to Lemma $2.5$ $\sigma (T^*)=(\sigma (T))^{\tilde .}$,
Then $f\in {\cal H}_{\infty }(T^*)$. In view of
$(f(T))^*y^*=y^*\circ f(T)$ for each $y^*\in X^*_q$, then due to
Lemma 2.5 $R(\zeta ^*;T^*)=(R(\zeta ;T))^*$, consequently,
$$(f(T))^*y^*=(2\pi )^{-1}(M^{-1})^*(\int_{\psi } (d\zeta ^*)R(\zeta
^*;T^*) ((f(\zeta ))^*y^*)),$$ where $(f(\zeta ))^*y^*:=y^*\circ
f(\zeta )$. If $f(\zeta )$ is represented by the series converging
int the ball: $$f(\zeta )=\sum_n \{ (a_n,\zeta ^n) \}
_{q(2m(n))},\mbox{ then } $$   $$f(\zeta ^*)=\sum_n \{ a_{n,1}\zeta
^{*n_1}...a_{n,m(n)}\zeta ^{*n_{m(n)}} \} _{q(2m(n))},\mbox{
consequently,} $$   $$[f(\zeta ^*)]^*=\sum_n \{ \zeta
^{n_{m(n)}}a_{n,m(n)}^*...\zeta ^{n_1} a_{n,1}^* \} _{q'(2m(n))},$$
where $q'(v)$ corresponds to the order of the associated
multiplication in the adjoint product in comparison with $q(v)$,
also  $$(f(T))^*=(2\pi )^{-1}(M^{-1})^*\int_{\psi }(d\zeta
^*)R(\zeta ^*;T^*) f^*(\zeta ^*)=f^*(T^*).$$
\par {\bf 10. Theorem.} {\it Let an
operator $T$ be bounded, $T\in L_q(X)$, where $X$ is a  Banach space
over ${\bf K}=\bf H$ or ${\bf K}=\bf O$, $f\in {\cal H}(T)$, then
$f(\sigma (T))=\sigma (f(T)).$}
\par {\bf Proof.} Suppose, that $b\in sp(T)$, then
$L_q(X)(T-bI)$ or $(T-bI)L_q(X)$ is the proper left or right ideal
in $L_q(X)$. Let for the definiteness $L_q(X)(T-bI)$ be the proper
left ideal in $L_q(X)$. In accordance with 2.6(i): $$f(T)=(2\pi
)^{-1}(\int_{\psi }f(z)R(z;T)dz)M^*.$$
 If ${\bf K}=\bf O$, then for
a given $b\in \bf K$ choose a copy of $\bf H$ embedded into $\bf O$,
such that $b\in \bf H$. Using the invariance of the integral along a
path for homotopic paths in $\rho (T)$ satisfying conditions of the
Theorem 3.9 in \cite{luoystoc}, it is possible to put $\psi \subset
\bf H$.
\par Due to formula 2.9(vii) for each bounded $T|_{X_{\bf H}}$
there is satisfied the equality $$\int_{\psi } f(z)
[R(z;T),(z-b)I]dz=0,$$ hence also on the entire Banach space $X$,
where $b\in \rho (T)$. There are satisfied the identities:
\par $[(zI-T)(z-b)]^{-1}(T-bI) = [(z-b)^{-1}(zI-T)^{-1}] [(T-zI)+(z-b)I]$
$= - (z-b)^{-1} + (z-b)^{-1}((zI-T)^{-1}(z-b))$, \\
since the algebra $\bf K$ is alternative, $(vK)(x) = v(K(x))$ for
each $K\in L_q(X)$, $v\in \bf K$, $x\in X$, also $(vB^{-1})\circ
B(x)=vx$ for each invertible operator $B$, for each vector $x\in X$
and each $v\in \bf K$. At the same time
\par $(z-b)^{-1}((zI-T)^{-1}(z-b))=
(z-b)^{-1}((z-b)(zI-T)^{-1})+(z-b)^{-1}([(zI-T)^{-1},(z-b)I])$,\\
but due to the alternativity of the algebra $\bf K$ there is the
equality: \par $ (z-b)^{-1}((z-b)(z I-T)^{-1})=(z I-T)^{-1} $,
consequently, \par $(z-b)^{-1}((zI-T)^{-1}(z-b))=(z I-T)^{-1}
+(z-b)^{-1}([(zI-T)^{-1},(z-b)I])$, therefore,
\par $[(zI-T)^{-1}-(z-b)^{-1}]=[(z I-T)^{-1}- (z-b)^{-1}I]
+(z-b)^{-1}([(zI-T)^{-1},(z-b)I])$. Then $$[f(T)-f(I)] = (2\pi
)^{-1} (\int_{\psi }f(z)[(zI-T)^{-1}-(z-b)^{-1}]dz)M^*$$
$$=((\int_{\psi }f(z) [(zI-T)(z-b)]^{-1}dz)(T-bI))M^*\in
L_q(X)(T-bI),$$ since the integral along a path is left linear, but
the integral of the term with the commutator is equal to zero, as in
the proof 2.9, consequently, $f(b)\in sp(f(T))$, that is,
$f(sp(T))\subset sp (f(T))$.
\par If $b\notin f(sp(T))$, then $(f(z)-b)^{-1}=:g(z)$ is holomorphic on
an open neighborhood for $sp(T)$. Due to Theorem 2.9 $g(T)$ is two
sided inverse operator for $f(T)-bI$ in $L_q(X)$, since
$g(z)(f(z)-b)=1$ on an open neighborhood of the spectrum $sp (T)$.
Thus, $b\notin sp (f(T))$, that is, $sp(f(T))\subset f(sp(T))$.
Therefore, $sp(f(T))=f(sp(T))$.
\par {\bf 11. Theorem.} {\it Let $f\in {\cal H}_{\infty }(T)$,
$T\in L_q(X)$, $X$ be a Banach space over ${\bf K}=\bf H$ or ${\bf
K}=\bf O$, $f(U)$ is open for some open subset $U\subset dom
(f)\subset \bf \hat H$, $g\in {\cal H}_{\infty }(T)$ and
$f(U)\supset \sigma (T)$, $dom (g)\supset f(U)$, then $F:=g\circ
f\in {\cal H}(T)$ and $F(T)=g(f(T)$.}
\par {\bf Proof.} By the supposition of this theorem the function $g$ is holomorphic
on an open subset $U_1$ containing $sp(f(T))$ and $f$ is holomorphic
on an open subset $U_2$ containing $sp(T)$. Due to Theorem 2.10
$f(sp(T))=sp(f(T))\subset U_1$, such that $sp(T)\subset
f^{-1}(U_1)$. Due to the continuity of $f$, $f^{-1}(U_1)\cap U_2=:U$
is the open subset in $\bf K$, on which there is given the
holomorphic function $g\circ f$, $sp (T)\subset U$. Thus, $ g\circ
f\in {\cal H}(T)$. Choose rectifiable closed paths (that is, loops)
$\psi $ and $\psi _g$ encompassing $sp(T)$ and $sp(f(T))$
respectively, where $\psi \subset U$, $\psi _g\subset U_2$. Then for
each $\zeta \in \psi _g$ the function $h_{\zeta }:=[\zeta -
f(z)]^{-1}$ is holomorphic on the open subset $V$ in $\bf K$, where
$\psi \subset V\subset U$ and $f(V)\subset U_1$. Then due to Theorem
2.9 $h_{\zeta }(T)=[zI-T]^{-1}$ for each $\zeta \in \psi _g$.
Therefore, $$g(f(T)) = (2\pi )^{-1}(\int_{\psi _g}f(\zeta )[\zeta
I-f(T)]^{-1}d\zeta )M^*_g = (2\pi )^{-1}(\int_{\psi _g}f(\zeta
)h_{\zeta }(T)d\zeta )M^*_g$$ $$=(2\pi )^{-1}(\int_{\psi _g} (((2\pi
)^{-1}\int _{\psi } h_{\zeta }(z)(zI-T)^{-1}dz)M^*) d\zeta )M^*_g,$$
where $M, M_g\in \bf K$, $|M|=|M_g|=1$, $Re (M)=Re (M_g)=0$, $M$
corresponds to the loop $\psi $, also $M_g$ corresponds to the loop
$\psi _g$. Consider bounded $T$ on all possible subspaces $X_{\bf
H}$ embedded into $X$ for all possible embeddings of $\bf H$ into
$\bf O$ as the subalgebra. In particular, for ${\bf K}=\bf H$ take
$X_{\bf H}=X$. Using the homotopy of the loops while accomplishing
conditions of Theorem 3.9 \cite{luoyst,luoystoc}, it is possible to
suppose $\psi $ and $\psi _g$ be contained in the same copy of the
skew field of quaternions $\bf H$. Choose $M=M_g$, then
$$g(f(T))=(2\pi )^{-1}(\int_{\psi }((2\pi )^{-1}\int_{\psi _g}
g(\zeta )[\zeta -f(z)]^{-1}d\zeta )M^*)(zI-T)^{-1}dz)M^*$$ $$=(2\pi
)^{-1}(\int_{\psi }(g\circ f)(z)(zI-T)^{-1}dz)M^*=(g\circ f)(T),$$
since the integral along paths is left and right linear. Due to
arbitrariness of $X_{\bf H}$ from this the statement of this theorem
follows in the case of $\bf O$ also, since the union of all copies
of the skew field of quaternions embedded into the algebra of
octonions covers it as the set: $\bigcup_{\theta } \theta ({\bf
H})=\bf O$, where $\theta : {\bf H}\hookrightarrow \bf O$ are
embeddings of $\bf H$ into $\bf O$ as subalgebras.
\par {\bf 12. Definition and Note.} Let $\cal A$ be a Banach space
and an algebra over $\bf K$, where ${\bf K}=\bf H$ or ${\bf K}=\bf
O$, with the unity $e$ having properties: $|e|=1$ and $|xy|\le |x|
|y|$ for each $x$ and $y\in \cal A$, then $\cal A$ is called the
Banach algebra or the $C$-algebra over $\bf K$. A Banach algebra
$\cal A$ we call quasicommutative, if there exists a commutative
algebra ${\cal A}_0$ over $\bf R$, such that ${\cal A}$ is
isomorphic with the algebra ${\cal A}_0\oplus {\cal A}_1i_1\oplus
... \oplus {\cal A}_mi_m$, where the algebras ${\cal
A}_0$,...,${\cal A}_m$ are pairwise isomorphic over $\bf R$, $ \{
i_0,i_1,...,i_m \} $ are standard generators of the algebra $\bf K$,
$m=3$ for $\bf H$, $m=7$ for $\bf O$.
\par Consider $X$ over $\bf R$: $\quad X=X_0i_0\oplus
X_1i_1\oplus ...\oplus X_mi_m$, where $X_0,...,X_m$ are pairwise
isomorphic Banach spaces over $\bf R$. Then ${\cal A}={\cal
A}_0\oplus {\cal A}_1i_1\oplus ...\oplus {\cal A}_mi_m$, where
${\cal A}_0$, ..., ${\cal A}_m$ are algebras over $\bf R$.
Multiplying $\cal A$ on $S\in \{ i_0,...,i_m \} $ we get
automorphisms of $\cal A$, consequently, ${\cal A}_0$, ..., ${\cal
A}_m$ are pairwise isomorphic.
\par {\bf 13. Definition and Note.}
A Banach algebra $\cal A$ over ${\bf K}=\bf H$ or ${\bf K}=\bf O$ is
supplied with the involution, when there exists the operation $*:
{\cal A}\ni T\mapsto T^*\in \cal A$, such that $(T^*)^*=T$,
$(T+V)^*=T^*+V^*$, $(TV)^*=V^*T^*$, $(\alpha T)^*=T^*\tilde \alpha $
for each $\alpha \in \bf K$, where ${\bf K}=\bf H$ or ${\bf K}=\bf
O$.
\par An element $x\in \cal A$  is called regular, if there exists
$x^{-1}\in \cal A$. In the contrary case it is called singular. Then
the spectrum $\sigma (x)$ for $x$ is defined as the set of all $z\in
\bf K$, for which $ze-x$ is singular, it spectral radius is the
following $|\sigma (x)|:=\sup_{z\in \sigma (x)} |z|$. A resolvent
set is defined as \par $\rho (x):= \{ z\in {\bf K}:$ $ze-x$
$\mbox{is regular} \} $ \\
and a resolvent is defined as
$R(z;x):=(ze-x)^{-1}$ for each $z\in \rho (x)$.
\par If it is given a Banach algebra $\cal A$ over $\bf K$ and its
closed subalgebra $\cal B$ over $\bf K$, then the factoralgebra
${\cal A}/\cal B$ by the definition consists of all elements, which
are classes of coset elements $[x]:=x+{\cal B}$, where $x\in \cal
A$, $S+V := \{ z: z=s+v, s\in S, v\in V \} $ for $S, V \subset \cal
A$. A subalgebra $\cal B$ is called left (right) ideal, if ${\cal
A}{\cal B}\subset \cal B$ (${\cal B}{\cal A}\subset \cal B$
respectively), if a subalgebra ${\cal B}$ is simultaneously a left
and right ideal, then $\cal B$ is called a two sided ideal or ideal.
\par {\bf 14. Lemma.} {\it A spectrum $\sigma (x)$ of an
element $x\in \cal A$ is a nonvoid compact subset in $\bf K$. Its
resolvent $x(z):=R(z;x)$ is holomorphic on $\rho (x)$, $x(z)$
converges to zero while $|z|\to \infty $ and
\par $x(z)-x(y)=x(z)(y-z)x(y)$ for each $y, z\in \rho (x)$.}
\par {\bf Proof} follows from the equalities $((ze-x) x(z))x(y)=x(y)$,
$x(z)(x(y) (ye-x))=x(z)$, $(ze-x)((x(z)-x(y))(ye-x))=(ye-x)-(ze-x)$
$=(y-z)e$, consequently, \par  $x(z)-x(y)=R(z;x)((y-z)R(y;x))$
$=x(z)((y-z)x(y))$. \\
Therefore $x(z)$ is continuous by $z$ on $\rho
(x)$ and there exists $(\partial [x(z+y) x^{-1}(y)]/\partial
z).h=-x(y)h$ for each $h\in \bf K$. For each marked point $y$ the
term $x^{-1}(y)$ is constant on $\cal A$, moreover, $(\partial
[x(z+y)x^{-1}(y)]/\partial {\tilde z})=0$, consequently, $x(z) \in
{\cal H} (\rho (x))$. The second statement follows from the
consideration of the complexification ${\bf C}\otimes \cal A$.
\par {\bf 15. Theorem.} {\it Let $\cal B$ be a closed proper ideal
over $\bf K$ in a quasicommutative Banach algebra $\cal A$ over the
algebra ${\bf K}=\bf H$ or ${\bf K}=\bf O$. The factoralgebra ${\cal
A}/\cal B$ is isometrically isomorphic with $\bf K$ if and only if
the ideal $\cal B$ is maximal.}
\par {\bf Proof.} The skew field of quaternions is produced from the field
of complex numbers by the application of the procedure of doubling
with the help of the generator $j=i_2$, also the algebra of
octonions is produced by the application of the doubling procedure
of the skew field of quaternions with the help of the generator
$l=i_4$. Then the algebra $\cal A$ is isomorphic to the direct sum
${\cal A}_0\oplus i_1{\cal A}_1\oplus ... \oplus i_m{\cal A}_m$,
where algebras ${\cal A}_0$,...,${\cal A}_m$ over the real field
$\bf R$ are pairwise isomorphic, $\{ i_0,i_1,...,i_m \} $ are the
standard generators of the algebra $\bf K$. Then the algebras ${\cal
A}_{0,1}:={\cal A}_0\oplus i_1{\cal A}_1$ and ${\cal B}_{0,1}:={\cal
B}_0\oplus i_1{\cal B}_1$ are commutative over the field of complex
numbers. Moreover, ${\cal B}_{0,1}$ is the ideal in ${\cal
A}_{0,1}$, since ${\cal A}{\bf K}= {\bf K}{\cal A}=\cal A$ and
${\cal B}{\bf K}= {\bf K}{\cal B}=\cal B$. Due to the theorem of the
Gelfand-Mazur (see Theorem 1 in \S 11 of  Chapter III \cite{nai})
the factoralgebra ${\cal A}_{0,1}/{\cal B}_{0,1}$ is isomorphic with
$\bf C$ if and only if a proper ideal ${\cal B}_{0,1}$ is maximal in
${\cal A}_{0,1}$. On the other hand, ${\cal A}$ and ${\cal B}$ are
produced from ${\cal A}_{0,1}$ and ${\cal B}_{0,1}$ respectively
with the help of the doubling procedure, one time in the case of
${\bf K}=\bf H$ and two times in the case of ${\bf K}=\bf O$,
therefore the factoralgebra ${\cal A}/{\cal B}$ is isomorphic with
$\bf K$ if and only if $\cal B$ is the maximal proper ideal in $\cal
A$.
\par {\bf 16. Definitions.} A $C^*$-algebra $\cal A$ over ${\bf K}=\bf H$
or ${\bf K}=\bf O$ is a Banach algebra over $\bf K$ with the
involution $*$, such that $|x^*x|=|x|^2$ for each $x\in \cal A$.
\par A scalar product in a vector space $X$ over $\bf H$
(that is, linear relative to the right and left multiplications
separately on scalars from $\bf H$) it is a biadditivite $\bf
R$-bilinear mapping $<*;*>: X^2\to \bf H$, such that
\par $(1)$ $<x;x>=\alpha _0$, where $\alpha _0\in \bf R$;
\par $(2)$ $<x;x>=0$ if and only if $x=0$;
\par $(3)$ $<x;y>=<y;x>^{\tilde .}$ for each $x, y \in X$;
\par $(4)$ $<x+z;y>=<x;y>+<z;y>$;
\par $(5)$ $<xa;yb>={\tilde a}<x;y>b$ for each
$x, y, z\in X$, $a, b\in \bf H$.
\par In the case of a vector space $X$ over $\bf O$ we consider an
$\bf O$-valued function on $X^2$ such that \par $(1')\quad <\zeta
,\zeta >=a$ with $a\ge 0$ and $<\zeta ,\zeta >=0$ if and only if
$\zeta =0$, \par $(2')\quad <\zeta ,z+\xi >=<\zeta ,z>+<\zeta ,\xi
>$, \par $(3')\quad <\zeta +\xi ,z>=<\zeta ,z>+<\xi ,z>$,
\par $(4')\quad <\alpha \zeta ,z>=\alpha <\zeta ,z>=<\zeta ,\alpha
z>$ for each $\alpha \in \bf R$ and $<\zeta a ,\zeta >={\tilde
a}<\zeta ,\zeta >$ for each $a\in \bf O$, \par $(5')\quad <\zeta
,z>^{\tilde .}=<z,\zeta >$ for each $\zeta , \xi $ and $z\in X$.
\par While the representation of $X$ in the form $X=X_0i_0\oplus X_1i_i\oplus ...
\oplus X_mi_m$, where $m=3$ for $\bf H$ and $m=7$ for $\bf O$, $X_0,
...,X_m$ are pairwise isomorphic $\bf R$-linear spaces, $ \{
i_0,...,i_m \} $ is the family of standard generators of the algebra
${\bf K}=\bf H$ or ${\bf K}=\bf O$, $i_0=1$, we shall suppose
naturally, that
\par $(6)$ $<x_p,y_q>\in \bf R$ and
\par $(7)$ $<x_pi_p,y_qi_q>=<x_p,y_q>i_p^*i_q$
for each $x_p\in X_p, y_q\in X_q$, $p, q\in \{ 0,1,2,...,m \} $.
\par If $X$ is complete relative to the norm topology
\par $(8)$ $|x|:=<x;x>^{1/2}$,  then $X$ is called $\bf K$ Hilbert space.
\par In particular, for $X=\bf K^n$ we can take the canonical scalar
product: \\
$(9)\quad <\zeta ;z>:= (\zeta ,z)=\sum_{l=1}^n\mbox{ }^l{\tilde
\zeta }\mbox{ }^lz$, where $z=(\mbox{ }^1z,...,\mbox{ }^nz)$,
$\mbox{ }^lz\in \bf K$.
\par {\bf 17. Lemma.} {\it The Banach algebra $L_q(X)$ on a Hilbert space
$X$ over ${\bf K}=\bf H$ or ${\bf K}=\bf O$ with the involution:
\par $(1)$ $<Tx;y>=:<x;T^*y>$ for each $x, y \in X$ \\
is the $C^*$-algebra.}
\par {\bf Proof.} From Formula $(1)$, Definitions 2.13 and 2.16
it follows the involutivity of the algebra $L_q(X)$, since
$(bT)(x):=b(T(x))$ and $(Tb)(x):=(T(x))b$ for each $T\in L_q(X)$,
$x\in X$, $b\in \bf K$, together with the Conditions $2.16(1-7)$
this gives $(bT)^*=T^*b^*$ due to the decomposition $T=T_0\oplus
T_1i_1\oplus ... \oplus T_mi_m$, $T_j: X\to X_j$ for each
$j=0,1,...,m$. At the same time $<T^*(T(x)),y>=<T(x),T(y)>$ for each
$x, y \in X$, in particular, $<T^*(T(x)),x>=\| T(x) \|^2$.
Therefore, \par $\| T \| ^2\le \sup _{\| x \| =1} |<T^*(T(x)),x>|\le
\sup_{\| x \| =1, \| y \| =1} |<T^*(T(x)),y>| $ \par $= \| T^*T \|
\le \| T ^* \| \| T \| $ \\ and analogous inequalities follow from
the susbtitution of $T$ on $T^*$, that gives \par $ \| T \| ^2\le \|
T^* \| \| T \| $ and $ \| T^* \| ^2 \le \| T \| \| T^* \| $, \\
consequently, $ \| T \| = \| T^* \| $ and $ \| T^*T \| = \| T \|
^2$.
\par {\bf 18. Lemma.} {\it If $\cal A$ is a quasicommutative
$C^*$-algebra over ${\bf K}=\bf H$ or ${\bf K}=\bf O$, then
$|x^2|=|x|^2$, $|x|=|x^*|$ and $I^*=I$, where $I$ is the unit in
$\cal A$.}
\par {\bf Proof.} Each vector $x\in \cal A$ can be represented in
the form: $x=x_0i_0+...+x_mi_m$, $x_j\in {\cal A}_j$ for each
$j=0,1,...,m$. Then $x^*=x_0^*i_0-x_1^*i_1-...-x_m^*i_m$, since
$(x_ji_j)^*=(-1)^{\kappa (i_j)}x_ji_j$, where $i_j\in \{ i_0,
i_1,...,i_m \}$ for each $j=0,1,...,m$, $\kappa (i_0)=0$, $\kappa
(i_j)=1$ for $j=1,2,...,m$. Therefore, $\quad [x,x^*]=0$ and
$|x_j^2|=|x_j|^2$. Then $|x|^2=|x_0|^2+|x_1|^2+...+|x_m|^2$ and
$|x^2|^2=|(x^2)^*x^2|=|(x^*)^2x^2|=|(xx^*)(xx^*)|=|x|^4$,
consequently, $|x^2|=|x|^2$. In view of $I=I_0$ there is
accomplished $I^*=I_0^*=I_0=I$.
\par {\bf 19. Definition.} A homomorphism $h: {\cal A}\to
\cal B$  of a $C^*$-algebra $\cal A$ into $\cal B$ over ${\bf K}=\bf
H$ or ${\bf K}=\bf O$, preserving the involution: $h(x^*)=(h(x))^*$
is called a $*$-homomorphism. If $h$ is a bijective $*$-homomorphism
$\cal A$ onto $\cal B$, then $h$ is called the $*$-isomorphism, also
$\cal A$ and $\cal B$ are called $*$-isomorphic. By $\sigma ({\cal
A})$ denote the structural space for $\cal A$ and it is called also
the spectrum for $\cal A$. The structural space is defined
analogously to the complex case with the help of Theorem 2.15 above.
\par {\bf 20. Proposition.} {\it For a Hilbert space $X$ over $\bf H$
spaces $L_l(X,{\bf H})$ and $L_r(X,{\bf H})$ are isomorphic with
$X$, for a Banach space $X$ a space $L_q(X)$ is isomorphic with
$L_l(X^2)$, $L_q({\bf H})=\bf H^4$, there exists a bijection between
the family of quasilinear operators $T$ on ${\cal D}(T)\subset X$
and the family of left linear operators $V$ on ${\cal D}(V)\subset
X^2$.}
\par {\bf Proof.} For $X$ there exists the decomposition $X=X_e\oplus
X_jj$, where $X_0$ and $X_j$ are isomorphic Hilbert spaces over $\bf
C$, that is, $z=z_1+z_2j$ for each $z\in X$, where $z_1\in X_e$,
$z_2\in X_j$. At the same time for $T\in L_q(X,Y)$ it is
accomplished the decomposition $T=T_e\oplus T_jj$, where $T_e: X\to
Y_e$, $T_j: X\to Y_j$, $Y$ is a Hilbert space over $\bf H$. In view
of $T(z)=T(z_1)+T(z_jj)$, there are equalities
$T_e(z)=T_{1,1}(z_1)+T_{1,2}(z_2)$,
$T_j(z)=T_{2,1}(z_1)+T_{2,2}(z_2)$, where each operator $T_{b,c}$ is
defined on a complex Hilbert space $X_c$ and takes the value in the
complex Hilbert space $X_b$ for each $b, c \in \{ 1, 2 \} $. Then
each operator $T_{b,c}$ due to its $\bf R$-linearity can be
represented in the form of two operators
$T_{b,c}=T_{b,c,1}+T_{b,c,2}$, where $T_{b,c,1}(x)$ is linear by
$x\in X_c$, also $T_{b,c,2}(x)$ is conjugate linear by $x\in X_c$,
that is, $T_{b,c,1}(px+qy)=pT_{b,c,1}(x)+qT_{b,c,1}(y)$ and
$T_{b,c,2}(px+qy)= {\bar p}T_{b,c,2}(x)+{\bar q}T_{b,c,2}(y)$ for
each $x, y\in X_c$ and $p, q\in \bf C$, where $\bar p$ denotes the
complex conjugated number $p\in \bf C$.
 Let $L_q(X)\ni \alpha =\alpha _ee+ \alpha
_ii+\alpha _jj+\alpha _kk$. In view of $S\alpha S\in L_q(X)$ for
each $S\in \bf H$, there exist quaternion constants $S_{m,l,n}$,
such that
\par $(1)\quad \alpha _m(x)m=
\sum_nS_{m,1,n}\alpha (x)S_{m,2,n}$ for each $m\in \{ e,i,j,k \} $, \\
where $S_{m,1,n}=\gamma _{m,n}S_{m,2,n}$ with $\gamma
_{m,n}=(-1)^{\phi (m,n)}/4\in \bf R$, $\phi (m,n)\in \{ 1, 2 \} $,
$S_{m,l,n}\in {\bf R}n$ for each $n\in \{ e,i,j,k \} $ (see \S \S
3.7, 3.28 \cite{luoyst}). Applying for $x$ the decomposition from \S
2.12, due to $(1)$ we get $4\times 4$-block form of operators over
$\bf R$ and the isomorphism of $L_q(X)$ with $L_l(X^2)$.
\par {\bf 21. Note.} For a Hilbert space $X$ over $\bf
O$ it can be used the decomposition $X=X_1\oplus X_2l$, where $X_1$
and $X_2$ are two isomorphic Hilbert spaces over $\bf H$, $l$ is the
generator of the doubling procedure of the skew field $\bf H$ up to
the algebra $\bf O$, then each operator $T$ on $X$ with values in
the Hilbert space $Y$ over $\bf O$ can be written in the form
$T(z)=T_1(z)\oplus T_2(z)l$ for each $z\in X$, where $T_1: X\to
Y_1$, $T_2: X\to Y_2$. At the same time
$T_p(z)=T_{p,1}(z_1)+T_{p,2}(z_2)$, where $T_{p,q}: X_q\to Y_p$ for
each $p, q \in \{ 1, 2 \} $, $z=z_1+z_2l$, $z_1\in X_1$, $z_2\in
X_2$.
\par For a left linear operator $T: D(T)\to X$ on a $\bf K$-linear
subspace $D(T)\subset X$ for a Banach space $X$ over ${\bf K}=\bf H$
or ${\bf K}=\bf O$ a nonzero vector $x\in D(T)$, satisfying the
condition $T(x)=\lambda x$ for some $\lambda \in \bf K$ is called
the eigenvector, also $\lambda $ is called the eigenvalue. Then
$T(bx)=bT(x)=b(\lambda x)$ for each $b\in \bf K$ and left linear
operator $T$, moreover, $\lambda I-T$ is not injective and $\lambda
\in \sigma (T)$. The set of all eigenvalues is called the point
spectrum and it is denoted by $\sigma _p(T)$. If $\lambda $ is not
an eigenvalue and if $Range (\lambda I-T):=(\lambda I-T)(D(T))$ is
not dense in $X$, then we say that $\lambda $ lies in the residue
spectrum, which is denoted by $\sigma _r(T)$. Then the notions of
the point $\sigma _p(T)$ and the residue $\sigma _r(T)$ spectra over
$\bf H$ due to Proposition 2.20 can be spread on quasilinear
operators.
\par Due to the Lebesgue theorem about decomposition of the measure:
if $\mu $ is a Borel measure on $\bf K$, then $\mu $ is uniquelly
representable in the form $\mu =\mu _{ac}+\mu _{sing}+\mu _p$, where
$\mu _{ac}$ is absolutely continuous relative to the Lebesgue
measure, also $\mu _{sing }$ is a continuous (that is, without
atoms) measure and singular relative to the Lebesgue measure on $\bf
K$ as the topological space ${\bf R}^{m+1}$, $\mu _p$ is purely
atomic (point) measure and singular relative to the Lebesgue
measure.
\par An operator $T$ in a Hilbert space $X$ over ${\bf K}=\bf
H$ or ${\bf K}=\bf O$ is called Hermitian, if $<Tx,y>=<x,Ty>$ for
each $x$ and $y$ from the domain of the definition ${\cal D}(T)$ of
the operator $T$, where ${\cal D}(T)\subset X$. An operator $T$ in
$X$ is called symmetrical, if it is Hermitian and ${\cal D}(T)$ is
dense in $X$. An operator $T$ in $X$ is called selfadjoint, when
${\cal D}(T)$ is dense in $X$ and $T=T^*$ on ${\cal D}(T)$.
\par {\bf 22. Theorem.} {\it A quasicommutative $C^*$-algebra
$\cal A$ over ${\bf K}=\bf H$ or ${\bf K}=\bf O$ is isometrically
$*$-isomorphic to the algebra $C(\Lambda ,{\bf K})$ of all
continuous $\bf K$-valued functions on its spectrum $\Lambda $.}
\par The {\bf Proof} follows from the fact, that the mapping
$x\mapsto x(.)$ from $\cal A$ into $C(\Lambda ,{\bf K})$ is the
$*$-homomorphism, where $x({\cal M})$ is defined by the equality
$x+{\cal M}=x({\cal M})+{\cal M}$ for each maximal ideal $\cal M$.
Let $x(\lambda )=\alpha _0i_0+ \alpha _1i_1+...+\alpha _mi_m$, then
$x^*(\lambda )= \beta _0i_0+\beta _1i_1+...+\beta _mi_m$, where
$\alpha _0,..., \beta _m\in \bf R$, $i_0=1$, $\{ i_0,i_1,...,i_m \}
$ are standard generators of the algebra $\bf K$. There exists the
decomposition for $X:=C(\Lambda ,{\bf K})$ from \S 2.12 with
$X_0=C(\Lambda ,{\bf R})$ and $x_p=\sum_nS_{p,1,n}zS_{p,2,n}{\tilde
i_p}$ for each $z\in \bf K$, where $z=x_0i_0+x_1i_1+...+x_mi_m$,
$x_j\in \bf R$ for each $j=0,1,...,m$ (see Proposition 2.20).
Therefore, it is applicable the Stone-Weiestrass theorem for $\bf
K$-valued functions. If $\lambda _1\ne \lambda _2$ are two maximal
ideals in $\Lambda $, then $y(\lambda _1)\ne y(\lambda _2)$ for
$y\in \lambda _1\setminus \lambda _2$. Consequently, the algebra of
functions $x(.)$ coincides with $C(\Lambda ,{\bf K})$.
\par {\bf 22.1. Definitions.} In the algebra $L_q(X)$ on a Hilbert
space $X$ over ${\bf K}=\bf H$ or ${\bf K}=\bf O$ we consider a base
of neighborhoods of an operator $A$ consisting of sets of the form
\par $V(A; x_1,...,x_n;b) := \{ B\in L_q(X): \| (B-A)x_j \| <b,
j=1,2,...,m \} $, \\
where  $x_1,...,x_n\in X$, $b>0$. A topology,
formed by such base is called the strong operator topology.
\par We introduce also the the base of all neighborhoods \par $W(A; x_1,...,x_n;
y_1,...,y_n; b):= \{ B\in L_q(X): | <x_j,(A-B)y_j> | <b, j=1,2,...,n
\} $, \\
a generated by it topology we call the weak operator
topology, where $x_1,...,x_n, y_1,...,y_n \in X$, $b>0$. This means
that it is induced by functionals $w_{x,y}(A) := <x,A(y)>$.
\par We recall, that a $C^*$-algebra over $\bf C$ is called the von Neumann
algebra, if it is contained in the algebra of $\bf C$-linear bounded
operators on a Hilbert space over $\bf C$, such that this algebra is
closed relative to the weak operator topology and contains the unit.
\par {\bf 22.2. Theorem.} {\it Strong and weak closures of a $\bf R$-convex
subset $Y$ in $L_q(X)$ for a Hilbert space $X$ over ${\bf K}=\bf H$
or ${\bf K}=\bf O$ coincide.} \par {\bf Proof.} In view of
$(A_1a_1+A_2a_2)(y) = (A_1(y))a_1 + (A_2(y))a_2$ for each $y\in X$,
where $A_1, A_2\in L_q(X)$, $a_1, a_2 \in \bf K$, then over $\bf H$
each functional $w_{x,y}(A)$ is right $\bf H$-linear, that is not in
general accomplished over $\bf O$. If an operator $A_s$ is
selfadjoint, $x_p\in X_0$, then
\par $<x_pi_p,A_s(y)i_s>=<x_p,A_s(y)>i_p^*i_s$, \\
where $\{ i_0,i_1,...,i_m \} $ are the standard generators of the
algebra $\bf K$, $i_0=1$, $X=x_0\oplus X_1i_1\oplus ... \oplus
X_mi_m$, where each $X_s$ is the Hilbert space over $\bf R$
isomorphic with $X_0$. Then $$<x,A(y)> = \sum_{p,s,q}
<x_p,A_s(y_qi_q)>i_p^*i_s,$$ where it is used the decomposition $A =
\sum_{s=0}^mA_si_s$ with the selfadjoint operators $A_s$ on $X$ with
values in $X_s$, where $X=X_0\oplus X_1i_1\oplus ... \oplus X_mi_m$,
$ \{ i_0,i_1,...,i_m \} $ are standard generators of the algebra
$\bf K$, $X_0,...,X_m$ are pairwise isomorphic Hilbert spaces over
$\bf K$, also each operator $A_s$ is selfadjoint, this means, that
it commutes with each generator $i_p$.
\par If an operator belongs to the strong operator closure
of a subset $Y$, then it belongs to a weak operator closure. Let now
$A$ belongs to the weak operator closure of $Y$. Take arbitrary
vectors $x_1,...,x_n,$ $y_1,...,y_n \in X$. Consider the space
$X^{\oplus n}:= X\oplus ... \oplus X$, which is equal to the
$n$-times direct sum of the space $X$. Denote by $B^{\oplus
n}(y_1,...,y_n):=(B(y_1),...,B(y_n))\in X^{\otimes n}$, $B^{\oplus
n} := B\oplus ... \oplus B$, where $B \in L_q(X)$. Then the set $
Y^{\oplus n}:= \{ B^{\oplus n}: B\in Y \} $ is convex over $\bf R$
in $L_q(X^{\oplus })$. In view of the fact that $A^{\oplus }$
belongs to the weak operator closure, the operator $A^{\oplus
n}(y_1,...,y_n)$ belongs to the weak closure of $Y^{\oplus n}
(y_1,...,y_n)$. Due to the  Hahn-Banach theorem over $\bf R$, also
considering the underlying spaces over $\bf R$, we get that
$A^{\oplus n}(y_1,...,y_n)$ belongs to the operator closure for
$Y^{\oplus n} (y_1,...,y_n)$ in $X^{\oplus n}$. Therefore, for each
$b>0$ there exists $F\in Y$, such that $ \| (F-A)x_j \| <b$ for each
$j=1,2,...,n$. Thus, the weak and strong operator closures for $Y$
coincide.
\par {\bf 22.3. Theorem.} {\it The unit ball in $L_q(X)$ for a Hilbert space
$X$ over ${\bf K}=\bf H$ or ${\bf K}=\bf O$ is compact in the weak
operator topolog.}
\par {\bf Proof.} Due to \S 2.22.2 the image $B_{x,y}$
of the unit ball $B(L_q(X),0,1)$ under the mapping $w_{x,y}(A)$ is
the ball in $\bf K$ of the radius $ \| x \| \| y \| $. The family $
\{ w_{x,y}: x, y \in X \} $ distinguishes points in $L_q(X)$,
therefore, the mapping $$\theta (A):= \prod_{x, y\in {\bf K}, \| x
\| =1, \| y \| =1} w_{x,y}(A)$$ is the homeomorphism on a subset in
$S:=\prod_{x, y\in {\bf K}, \| x \| =1, \| y \| =1} B_{x,y}$. The
topological space $S$ is compact due to the Tychonoff theorem, since
each ball $B_{x,y}$ is compact. In view of $$<x,A(y)> = \sum_{p,s,q}
<x_p,A_s(y_qi_q)>i_p^*i_s,$$ also since each continuous $\bf
R$-bilinear form $b_{p,s,q}$ on $X_p\times X_qi_q$ with values in
${\bf R}i_p^*i_s$ due to Riesz representation theorem has the form
$b_{p,s,q}(x_p,y_qi_q)=<x_p,A_s(y_qi_q)>i_p^*i_s$ for some bounded
$\bf R$-linear operator $A_s$ on $X_qi_q$. If $b$ is the limit point
for $\theta (B(L_q(X),0,1))$, then $b(x,y) = \sum_{p,s,q} b_{p,s,q}
(x_p,y_qi_q)$, where each form $b_{p,s,q}$ is $\bf R$-bilinear on
$X_p\times X_qi_q$ and takes values in ${\bf R}i_p^*i_s$, since for
each $t>0$ there exist $F_s\in L_q(X)$, $x_1, x_2\in X_p, y_1,
y_2\in X_q$, such that \par $|a b_{p,s,q}(x_n,y_k)- a<x_n,
T_s(y_ki_q)>|<t$, $\quad |b_{p,s,q}(x_n,y_k )- <x_n,
T_s(y_ki_q)>|<t$,
\par $| b_{p,s,q}(ax_1+x_2,y_k)- <ax_1+x_2, T_s(y_ki_q)>|<t$, \par $|
b_{p,s,q}(x_n,ay_1+y_2)- <x_n, T_s((ay_1+y_2)i_q)>|<t$ \\
for each $j=1, 2; k=1, 2$, where $a\in \bf R$. Then
$|b_{p,s,q}(ax_1+x_2,y_k)- ab_{p,s,q}(x_1,
y_ki_q)-b_{p,s,q}(x_2,y_k)|<3t$ and $|b_{p,s,q}(x_j,ay_1+y_2)-
ab_{p,s,q}(x_j, y_1i_q)-b_{p,s,q}(x_j,y_2)|<3t$. In view of
arbitrariness of a small $t>0$ each $b_{p,s,q}$ is the $\bf
R$-bilinear mapping on $X_p\times X_qi_q$.  Thus, there exists $A\in
L_q(X)$, for which $b(x,y)=<x,A(y)>$. Thus, $\theta (B(L_q(X),0,1))$
is closed in $S$, consequently, $\theta (B(L_q(X),0,1))$ is compact.
\par {\bf 22.4. Lemma.} {\it  If $ \{ S_a: a \in \Upsilon \} $ is a
monotonely increasing net of selfadjoint operators on a Hilbert
space $X$ over ${\bf K}=\bf H$ or ${\bf K}=\bf O$ and $S_a\le kI$
for each $a$, where $S_a\in L_q(X)$, $0<k<\infty $, $\Upsilon $ is a
nonvoid directed set, then $ \{ S_a: a\in \Upsilon \} $ converges
relative to the strong operator topology to the selfadjoint operator
$S$, moreover, $S$ is the least upper bound for $ \{ S_a: a\in
\Upsilon \} $.}
\par {\bf Proof.} Since the convergence of the net
$ \{ S_a: a\in \Upsilon \} $ is equivalent to the convergence of the
net $ \{ S_a: a\in \Upsilon , a\ge a_0 \} $ for some $a_0\in
\Upsilon $, then without restriction of the generality it can be
supposed, that the net $ \{ S_a : a\ge a_0 \} $ is bounded from
below, therefore $ - \| S_{a_0} \| I \le S_a\le kI$ and $ \{ S_a: a
\} $ is the bounded set of operators. Due to Theorem 2.22.3 there
exists a subnet $ \{ S_c: c \} $ converging in the weak operator
topology to the operator $S\in L_q(X)$. In view of the fact that
$S_a$ is monotonely increasing, that is, $<S_ax,x>\ge <S_cx,x>$ for
each $a\ge c$ and $x\in X$, then $S\ge S_a$ for each $a$. If $a\ge
c$, then $0\le S-S_a\le S-S_c$, and also \par $0\le <(S-S_a)x,x> =
\| (S-S_a)^{1/2} x\| ^2 \le <(S-S_c)x,x>$ \\
converges to zero by
$c\in \Upsilon $ for each $x\in X$. Thus, $ \{ (S-S_a)^{1/2}: a \} $
converges relative to the strong operator topology to zero, also due
to continuity of the multiplication relative to the strong operator
topology on bounded subsets of operators, $ \{ (S-S_a): a \} $ also
converges to zero. If $F\ge S_a$ for each $a$, then $<Fx,x>\ge
<Sx,x>$ for each $x\in X$, since $<Sx,x>=\lim_a <S_ax,x>$, also
$<Fx,x>\ge <S_ax,x>$ for each $a$. Thus, $F\ge S$, that is, $S$ is
is the least upper bound.
\par {\bf 22.5. Lemma.} {\it  If $A$ is a bounded operator on
a Hilbert space $X$ over ${\bf K}=\bf H$ or ${\bf K}=\bf O$, $0\le
A\le I$, then $ \{ A^{1/n}: n\in {\bf N} \} $ is the monotonely
increasing sequence of operators, moreover, $\lim_{n\to \infty
}A^{1/n}$ is the projection operator on the closure of the region of
values of the operator $A$.}
\par {\bf Proof.} Consider a $C^*$-algebra generated by
operators $A$ and $I$ (see Theorem 2.22), then $\{ A^{1/n}: n \} $
is a monotonely increasing sequence bounded from above by the unit
operator $I$. Due to Lemma 2.22.4 $\{ A^{1/n}: n \} $ has the limit
$P$ relative to the strong operator topology, where $A^{1/n}\le P$
for each $n$. At the same time, $P^2$ is the limit of the sequence $
\{ A^{2/n}: n \} $, which is subsequence of the sequence $\{
A^{1/n}: n \} $, then $P^2=P$, that is, $P$ is the projection
operator. To the algebra $C(\Lambda ,{\bf K})$ it can be applied the
Stone-Weierstrass Theorem over $\bf K$ (see also \S 2.7 in
\cite{luoystoc}), therefore $A^{1/n}$ is the limit in the strong
operator topology of a family of polynomials of $A$ over $\bf K$
with a nozero constant term. Therefore, $A^{1/n}x=0$, if $Ax=0$,
consequently, $Px=0$. If $Px=0$, then \par $0=<Px,x>\ge
<A^{1/n}x,x>= \| A^{1/(2n)} x \| ^2\ge 0$, \\
that is, $A^{1/(2n)}x=0$, also then $Ax=0$. Thus,
$A^{-1}(0)=P^{-1}(0)$, consequently, $P$ is the projection operator
onto $cl (Range (A))$.
\par {\bf 22.6. Note.} Let $X$ be a Hilbert space
over ${\bf K}=\bf H$ or ${\bf K}=\bf O$, then there exists an
underlying Hilbert space $X_{\bf R}$ over $\bf R$. The scalar
product on $X$ with values in $\bf K$ from \S 2.16 induces the
scalar product $(x,y):=Re <x,y>$ in $X_{\bf R}$ with values in $\bf
R$. Then for a $\bf K$-vector subspace $Y$ in $X$ it can be created
the orthogonal complement $Y^{\perp }$ relative to the scalar
product $(x,y)$. At the same time the equation ${\hat E}_e(y+z)=y$
for $y\in Y$, $z\in Y^{\perp }$ defines the $\bf R$-linear operator
${\hat E}_e$ acting on a Hilbert space $X$, that is, ${\hat E}_e$ is
the projection on $Y$ parallel to $Y^{\perp }$. Moreover, ${\hat
E}_e^2={\hat E}_e$,
\par $(i)$ $Y= \{ y: y={\hat E}_e(x), x\in X \} = \{ y\in X: {\hat E}_e(y)=y \} $, also
\par $(ii)$ $Y^{\perp } = \{ z\in X: {\hat E}_e(z)=0 \} $.
\par Mention, that $(I-{\hat E}_e)$ is the orthogonal projection from $X$ on
$Y^{\perp }$, since $(Y^{\perp })^{\perp }=Y$ and $(I-{\hat
E}_e)(z+y)=z$ for each $z\in Y^{\perp }$ and $y\in Y$. Consider the
decomposition of a Hilbert space $X=X_0\oplus X_1i_1\oplus ...
\oplus X_mi_m$ over $\bf K$, where $ \{ i_0,i_1,...,i_m \} $ is the
family of standard generators of the algebra $\bf K$, $i_0=1$,
$X_0,...,X_m$ are pairwise isomorphic Hilbert spaces over $\bf R$.
\par The family of operators ${\hat E}(a,*)\in L_q(X)$, $a\in \bf K$ we call the
$\bf K$-graded projection operator on $X$, if it satisfies
conditions $(1-7)$ below:
\par $(1)$ ${\hat E}(x)={\hat E}_e(x)$ for each $x\in X_0$, where
\par $(2)$ ${\hat E}(x):={\hat E}(1,x)$ for each $x\in X$,
\par $(3)$ ${\hat E}_e^2(x)={\hat E}_e(x)$ for each $x\in X$,
\par $(4)$ ${\hat E}(s,x)={\hat E}(1,s x)$ for each $s\in \bf K$
and $x\in X_0$,
\par $(5)$ ${\hat E}(s x)=:{\hat E}_s(x)$ for each $s\in \bf K$ and $x\in X_0$,
\par $(6)$ ${\hat E}_s({\hat E}_q(x))={\hat E}_{sq}(x)$ for each
$s, q \in \{ i_0,...,i_m \} $ and $x\in X_0$,
\par $(7)$ ${\hat E}_s^*=(-1)^{\kappa (s)}{\hat E}_s$ on $X_0$
for each $s\in \{ i_0,...,i_m \} $, where $\kappa (s)=2$ for
$s=i_0$, $\kappa (s)=1$ for $s\in \{ i_1,i_2,...,i_m \} $, ${\hat
E}_s^*$ denotes the adjoint operator.
\par {\bf 22.7. Proposition.} {\it There exists a bijective correspondence
(up to an isometrical \\
automorphism of the algebra $\bf K$) between
closed $\bf K$-vector subspaces $Y$ of a Hilbert space $X$ over
${\bf K}=\bf H$ or ${\bf K}=\bf O$ and $\bf K$-graded projection
operators satisfying conditions $2.22.6(i,ii)$. At the same time
projection operators ${\hat E}_e$ are nonnegative, also $ \| E \|
=1$, if $E\ne 0$.}
\par {\bf Proof.} Let be given a closed $\bf K$-vector
subspace $Y$ in $X$, then there exists a bijective correspondence
between subspaces $Y_0$ over $\bf R$ in $X_0$ and projection
operators ${\hat E}_e$ in $X_0$, which is established with the help
of Relations $2.22.6(i,ii)$. On the other hand, between $Y_0$ and
$Y$ there exists a bijective correspondence. At the same time $\bf
K$-graded projective operator on $X$ is uniquely characterized by
Relations $2.22.6(1-7)$ by the way of setting ${\hat E}_e$ on $X_0$.
The restriction ${\hat E}_e|_{X_0}$ is the selfadjoint operator.
Moreover, ${\hat E}_e(X_0)=Y_0$ if and only if ${\hat E}(X)=Y$,
since \par ${\hat E}(x)=\sum_{p=0}^m{\hat E}(i_px_p)$, \\
where $x_p\in X_p$, but $X_0,...,X_m$ are pairwise isomorphic, also
${\hat E}(i_px_p)={\hat E}_{i_p}(x_p)$, ${\hat E}(a_1+a_2,x)= {\hat
E}(1,(a_1+a_2)x)={\hat E}(1,a_1x)+{\hat E}(1,a_2x)={\hat
E}(a_1,x)+{\hat E}(a_2,x)$ for each $a_1, a_2\in \bf K$, $x\in X_0$.
\par  From the relation ${\hat E}_s({\hat E}_s(x))={\hat E}_{s^2}(x)$
it follows, that ${\hat E}_s^2|_{X_0}=-{\hat E}_e|_{X_0}=-{\hat
E}_e^2|_{X_0}$ for each $s=i_1,i_2,...,i_m$. Therefore, ${\hat
E}^2(i_pX_p)\subset X_0$, but $i_pX_p$ is the orthogonal to $X_0$
subspace relative to the scalar product $(x,y)$, therefore, ${\hat
E}_e(X_0)$ is orthogonal to ${\hat E}(i_pX_p)$ relative to the
scalar product $(x,y)$, also due to 2.22.6(6) ${\hat E}(i_pX_p)$ is
orthogonal to ${\hat E}(i_qX_q)$ for each $p\ne q\in \{ i_0,...,i_m
\} $. The spaces $i_pX_p$ and $i_qX_q$ are pairwise isomorphic over
$\bf R$, consequently, ${\hat E}(i_pX_p)$ and ${\hat E}(i_qX_q)$ are
pairwise isomorphic over $\bf R$ and orthogonal relative to scalar
products $(x,y)$. Then together with ${\hat E}_s^2|_{X_0}=-{\hat
E}_e|_{X_0}=-{\hat E}_e^2|_{X_0}$ this gives that ${\hat E}(X)$ is
the $\bf K$-vector space. Moreover, \par $<{\hat
E}_{i_p}(x_p),i_qx_q>=$ $<{\hat E}(i_px_p),i_qx_q>=$ $<x_p,{\hat
E}^*_{i_p}(i_qx_q)>=$ \par $(-1)^{\kappa (i_p)} <x_p,{\hat
E}_{i_p}(i_qx_q)>=$ $(-1)^{\kappa (i_p)}<x_p,{\hat E}_e(
i_pi_qx_q)>=$ $<x_p,{\hat E}_{ip^*i_q}(x_q)>=$ \par $<x_p,{\hat
E}_e(i_p^*i_qx_q)>=$ $<{\hat E}_e(x_p),i_p^*i_qx_q>$ \\
for each $p\ne q$. Then $\| {\hat E} \| \le 1$, on the other hand,
$\| {\hat E}|_{X_0} \| =1$, for ${\hat E}\ne 0$. \par  For each $x,
y\in X$ due to Condition 2.22.6 the identities  \par $<{\hat
E}(x),y>=\sum_{p,q=0}^m <{\hat E}(i_px_p),i_qy_q>=$
$\sum_{p,q=0}^m<{\hat E}(i_px_p),y_q>i_q=$ \par
$\sum_{p,q=0}^m<i_px_p,{\hat E}^*(y_q)>i_q=$
$\sum_{p,q=0}^m<i_px_p,{\hat E}^*_e(y_q)>i_q=$ \par
$\sum_{p,q=0}^m<x_p,{\hat E}(y_q)>i_p^*i_q=$ $\sum_{p,q=0}^m<{\hat
E}(x_p),y_q>i_p^*i_q$ \\
are accomplished, since $<x_p,y_p>\in \bf
R$. In view of ${\hat E}(X_p)=Y_0$ for each $p$ it follows, that
$<{\hat E}(x),y>=0$ if and only if $y_q\in Y_0^{\perp }$ for each
$q$. Thus, $Y={\hat E}(X)={\bf K}Y_0$, $X=Y\oplus Y^{\perp }$, where
$Y^{\perp }={\bf K}Y_0^{\perp}$, $X_0=Y_0\oplus Y_0^{\perp }$, that
is, $Y\cap Y^{\perp }= \{ 0 \} $. Thus, $<{\hat E}(x),y>=0$ for each
$x\in X$ if and only if $y\in Y^{\perp }$. Evidently the graded
projection operator is defined uniquely up to a $\bf R$-linear
isometric automorphism of the algebra $\bf K$.
\par {\bf 22.8. Definitions.} For a given family $ \{ Y_a: a \in \Upsilon \} $
of closed $\bf K$-vector subspaces of a Hilbert space $X$ over ${\bf
K}=\bf H$ or ${\bf K}=\bf O$ there exists the greatest closed $\bf
K$-vector subspace $\wedge _aY_a$, which is contained in each of
$Y_a$, also there exists the least closed $\bf K$-vector subspace
$\vee _aY_a$, which contains each subspace $Y_a$, $\wedge
_aY_a=\bigcap _aY_a$, $\vee _aY_a=cl (\bigcup_aY_a)$. Due to
Proposition 2.22.7 the relation of partial ordering of $\bf
K$-vector subspaces $Y_a$ in $X$ induces the ordering of the family
of graded projection operators $\mbox{ }_a{\hat E}\le \mbox{
}_c{\hat E}$ if and only if $Y_a\subset Y_c$, moreover, there exists
the infimum $\wedge _a\mbox{ }_a{\hat E}_a$, also the supremum $\vee
_a \mbox{ }_a{\hat E}_a$ in the family of graded projection
operators. In particular, ${\hat E}\wedge {\hat F}$ and ${\hat
E}\vee {\hat F}$ are called the intersection and the union of two
graded projection operators $\hat E$ and $\hat F$.
\par {\bf 22.9. Proposition.} {\it  If ${\hat E}$ and $\hat F$ are
graded projection operators from a Hilbert space $X$ over ${\bf
K}=\bf H$ or ${\bf K}=\bf O$ onto a closed $\bf K$-vector subspace
$Y$ and $Z$ in $X$ respectively, then the following conditions are
equivalent:
\par $(1)$ $Y\subset Z$, \par $(2)$ ${\hat F}\hat E$ coincides with
$\hat E$ and with ${\hat E}\hat F$ (up to the isometric automorphism
of the algebra $\bf K$), \par $(3)$ $\| {\hat E}(x) \| \le \| {\hat
F}(x) \| $ for each $x\in X$, \par $(4)$ ${\hat E}\le \hat F$.}
\par {\bf Proof.} If $Y\subset Z$ are two closed
$\bf K$-vector subspaces in $X$, then ${\hat E}(x)\in Y\subset Z$
for each $x\in X$, therefore ${\hat F}({\hat E}(X))={\hat E}(X)$,
${\hat F}_e({\hat E}_e(x_0))={\hat E}_e(x_0)$ for each $x_0\in X_0$.
On the other hand, ${\hat E}{\hat F}=({\hat F}{\hat E})^*$, but
${\hat E}_e|_{X_0}$ is selfadjoint, consequently, from $(1)$ it
follows $(2)$. If ${\hat E}\hat F$ coincides with ${\hat E}$ up to
an isometric automorphism of the algebra $\bf K$, then $\| {\hat
E}(x) \| = \| {\hat E}({\hat F}(x)) \| \le \| {\hat F}(x) \| $ for
each $x\in X$, since $ \| {\hat E} \| \le 1$, therefore, from $(2)$
it follows $(3)$. If $(3)$ is satisfied, then from $<{\hat
E}(x),x>=<{\hat E}^2(x),{\hat E}(x)> = \| {\hat E}(x) \| ^2$ (see \S
2.22.8) and from the analogous identity $<{\hat F}(x),x>= \| {\hat
F}(x) \| ^2$ it follows $(4)$. If ${\hat E}\le \hat F$, then from
the identity $$<{\hat E}(x),y>= \sum_{p,q=0}^m<{\hat
E}(x_p),y_q>i_p^*i_q$$ for each $x, y \in X$, where $<{\hat
E}(x_p),y_q>\in \bf R$ it follows, that for each $y\in Y$ are
satisfied the equalities \par $\| y \| ^2 =<{\hat E}(y),y>\le {\hat
F}(y),y>= \| {\hat F}(y) \| ^2\le \| y \| ^2$,\\
 whence $\| {\hat
F}(y) \| = \| y \| $ and $y\in Z$, that is, from $(4)$ it follows
$(1)$.
\par {\bf 22.10. Corollary.} {\it Suppose that $\hat E$ and $\hat F$ are
graded operator projections from a Hilbert space $X$ over ${\bf
K}=\bf H$ or ${\bf K}=\bf O$ onto closed $\bf K$-vector subspaces
$Y$ and $Z$ respectively, also ${\hat E}\le \hat F$, then ${\hat
F}-\hat E$ are graded operator projections from $X$ onto $Z\wedge
Y^{\perp }$.}
\par {\bf Proof.} The mapping ${\hat E}\mapsto I-\hat E$
inverses the ordering of graded operator projections, therefore,
\par $\vee _a(I-\mbox{ }_a{\hat E})=I - \wedge _a \mbox{ }_a{\hat E}$,
$\quad \wedge (I-\mbox{ }_a{\hat E})=I - \vee _a \mbox{ }_a{\hat E}$
\\ for each family $ \{ \mbox{ }_a{\hat E}: a\in \Upsilon \} $ of
graded projection operators. This gives relations $\vee _aY_a^{\perp
}=(\wedge _aY_a)^{\perp }$, $\wedge _aY_a^{\perp }=(\vee
_aY_a)^{\perp }$ for each family of $\bf K$-vector closed subspaces
in $X$. For a $\bf K$-vector subspaces $A$ and $B$ in $X$ the
orthogonalities relative to the scalar products $(x,y)$ and $<x,y>$
coincide, since ${\bf K}A=A\bf K$, ${\bf K}B=B\bf K$. In view of
$Y={\hat E}(X)$, $Z={\hat F}(X)$, also due to $$<{\hat E}(x),y>=
\sum_{p,q=0}^m<{\hat E}(x_p),y_q>i_p^*i_q$$ for each $x, y \in X$,
where $<{\hat E}(x_p),y_q>\in \bf R$, then $<{\hat E}{\hat
F}(x),y>=<{\hat F}(x),{\hat E}(y)>$ for each $x, y\in X$, therefore,
$Y$ is orthogonal to $Z$ if and only if ${\hat E}{\hat F}=0$. If it
is accomplished the latter equality, then ${\hat F}{\hat E}=({\hat
E}{\hat F})^*=0$, hence ${\hat E}_e|_{X_0}$ and ${\hat F}_e|_{X_0}$
commute, moreover, their product is equal to zero, consequently,
${\hat E}\vee {\hat F}={\hat E}+{\hat F}$, that is, $Y\vee Z=Y+Z$.
\par {\bf 22.11. Corollary.} {\it  If $\hat E$ and $\hat F$ are graded operators
of projections from a Hilbert space $X$ over ${\bf K}=\bf H$ or
${\bf K}=\bf O$ onto closed subspaces $Y$ and $Z$ in $X$
respectively, also ${\hat E}\le \hat F$, then ${\hat F}-\hat E$ are
graded projection operators from $X$ on $Z\vee Y^{\perp }$.}
\par {\bf Proof.} Due to Proposition 2.22.9 ${\hat E}{\hat
F}={\hat F}{\hat E}=\hat E$, then graded operator projections ${\hat
F}(I-{\hat E})$, $(I-{\hat E}){\hat F}$ and ${\hat F}-\hat E$
coincide up to an isometric automorphism of the algebra $\bf K$,
then ${\hat F} - {\hat E}$ is the graded operator of projection
${\hat F}\vee (I-{\hat E})$ (up to an isometric automorphism of the
algebra $\bf K$) from $X$ on $Z\vee Y^{\perp }$.
\par {\bf 22.12. Proposition.} {\it  If $
\{ \mbox{ }_a{\hat E}: a\in \Upsilon \} $ is an increasing net of
graded operators  of projections from a Hilbert space $X$ over ${\bf
K}=\bf H$ or ${\bf K}=\bf O$, ${\hat E}=\vee  \mbox{ }_a{\hat E}$,
then ${\hat E}(x)=\lim_a \mbox{ }_a{\hat E}(x)$ for each $x\in X$.}
\par {\bf Proof.} In view of the fact that $ \{ \mbox{ }_a{\hat E}(X): a \} $
is an increasing net of closed subspaces in $X$ it follows, that
$\bigcup \mbox{ }_a{\hat E}(X)$ is the $\bf K$-vector subspace in
$X$ with the closure relative to the norm equal to ${\hat E}(X)$.
Suppose that $x\in X$ and $b>0$. In view of ${\hat E}(x)\in {\hat
E}(X)$ there exists $y\in \mbox{ }_a{\hat E}(X)$ for some $a$ such
that $ \| {\hat E}(x)-y \| <b$. For $c\ge a$ the inequalities
$\mbox{ }_a{\hat E}\le \mbox{ }_c{\hat E}\le {\hat E}$, $y\in \mbox{
}_a{\hat E}(X)\subset \mbox{ }_c{\hat E}(X)\subset {\hat E}(X)$ are
accomplished. Thus, \par $ \| {\hat E}(x) - \mbox{ }_c{\hat E}(x) \|
= \| {\hat E}({\hat E}(x) - y)-\mbox{ }_c{\hat E}({\hat E}(x)-y) \|
\le \| {\hat E} -\mbox{ }_c{\hat E} \| \| {\hat E}(x)-y \| <b$.
\par {\bf 22.13. Corollary.} {\it If $\{ \mbox{ }_a{\hat E} : a\in \Upsilon \} $
is a descending net of graded operators of projections from a
Hilbert space $X$ over ${\bf K}=\bf H$ or ${\bf K}=\bf O$, ${\hat
E}=\wedge \mbox{ }_a{\hat E}$, then ${\hat E}(x)=\lim_{a\in \Upsilon
} \mbox{ }_a{\hat E}(x)$ for each $x\in X$.}
\par The {\bf Proof} follows by the way of the application of Proposition 2.22.12
to the increasing net $ \{ I- \mbox{ }_a{\hat E}: a\in \Upsilon  \}
$.
\par {\bf 22.14. Definition.} If a family of graded operators of projections
$ \{ \mbox{ }_a{\hat E}: a\in \Upsilon \} $ satisfies the relations
$\mbox{ }_a{\hat E} \mbox{ }_c{\hat E}=0$ for each $a\ne c\in
\Upsilon $ (that is equivalent to the orthogonality of $\bf
K$-vector subspaces $\mbox{ }_a{\hat E}(X)$ and $\mbox{ }_c{\hat
E}(X)$), then such family is called orthogonal.
\par {\bf 22.15. Proposition.} {\it If it is given an orthogonal family
of graded operators of projections $ \{ \mbox{ }_a{\hat E}: a\in
\Upsilon \} $ acting on a Hilbert space $X$ over ${\bf K}=\bf H$ or
${\bf K}=\bf O$, ${\hat E}:=\vee \mbox{ }_a{\hat E}$, then ${\hat
E}(x)=\sum \mbox{ }_a{\hat E}(x)$ for each $x\in X$.}
\par {\bf Proof.} For a finite family $\Upsilon $ this statement follows
from Corollary 2.22.10 by the way of mathematical induction by the
number of elements in $\Upsilon $. For an infinite $\Upsilon $ let
$\cal G$ denotes the family of all finite subsets in $\Upsilon $.
For each $G\in \cal G$ consider $F_G := \{ \vee \mbox{ }_a{\hat E}:
a\in G \} $, then ${\hat F}_G = \{ \sum \mbox{ }_a{\hat E}: a\in G
\} $ and ${\hat F}_G$ is the increasing family of graded projection
operators, where $\cal G$ is directed by the inclusion of subsets
from $\cal G$, then ${\hat E} = \vee _{G\in \cal G}{\hat F}_G$. Due
to Proposition 2.22.12 the net ${\hat F}_Gx$ converges to ${\hat
E}(x)$ by the norm in $X$.
\par {\bf 22.16. Definition.} For an operator $T\in L_q(X)$ for
a Hilbert space $X$ over ${\bf K} = \bf H$ or ${\bf K}=\bf O$ we
define two subspaces $ker_{\bf K} (T) = cl (span_{\bf K} \{ x\in X:
T(x)=0 \} )$ which is the kernel of the operator and $[T(X)]_{\bf K}
:= cl (span _{\bf K} \{ T(x): x \in X \} )$ which is the space of
values of the operator $T$, where because of in general $\bf
K$-nonlinearity of the operator $T$ there are taken the $\bf
K$-vector spans $span_{\bf K}(T^{-1}(0))$ and $span_{\bf K} (T(X))$
for $T^{-1}(0)$ and $T(X)$, where $cl(A)$ denotes the closure of a
subset $A$ in $X$. With these subspaces we associate two graded
operator projections ${\hat N}(T)$ and ${\hat R}(T)$ from $X$ on
$ker_{\bf K} (T)$ and $[T(X)]_{\bf K}$ respectively.
\par {\bf 22.17. Proposition.} {\it If subspaces $T^{-1}(0)$
and $T(X)$ are \\
$\bf K$-vector spaces, then $ker_{\bf K}=T^{-1}(0)$, $[T(X)]_{\bf
K}= cl (T(X))$, ${\hat R}(T)=I-{\hat N}(T^*)$, \\
${\hat N}(T)=I-{\hat
R}(T)$, ${\hat R}(T^*T)={\hat R}(T^*)$, ${\hat N}(T^*T)={\hat
N}(T)$.}
\par {\bf Proof.} If a subspace $S$ is a
$\bf K$-vector space, then $span_{\bf K}S=S$, since $span_{\bf
K}(S):= \{ z\in X:$
$z=a_1(x_1b_1)+...+a_m(x_mb_m)+(c_1y_1)t_1+...+(c_ly_l)t_l$
$a_1,...,a_m, b_1,...,b_m, c_1,...,c_l, t_1,...,t_l\in {\bf K},$
$x_1,...,x_m, y_1,...,y_l\in X \} $. Therefore, $[S]_{\bf K}=cl
(S)=cl (span_{\bf K}(S))$, since the closure of the $\bf K$-vector
subspace $S$ in the Banach space over $\bf K$ also is the $\bf
K$-vector space. If the suppositions of this proposition are
accomplished, then $T^{-1}(0)=ker_{\bf K}(T)$ due to the continuity
of the operator $T$ and $[T(X)]_{\bf K}= cl (T(X))$. Then $ \{ x\in
X: T(x)=0 \} $ $= \{ x\in X: <T(x),y> =0$ $\mbox{for each}$ $y\in X
\} $  $= \{ x\in X:$ $<x,T^*y>=0$ $\mbox{for each }$ $y\in X \} $ $=
(T^*(X))^{\perp } =[T^*(X)]^{\perp }$, that is, $T^*(X)$ is the $\bf
K$-vector subspace. From this it follows, that ${\hat N}(T)=I-{\hat
R}(T^*)$. The substitution of $T$ on $T^*$ gives the equality ${\hat
N}(T^*)=I-{\hat R}(T)$ and $(T^*)^{-1}(0)=(T(X))^{\perp }$ also is
the $\bf K$-vector subspace. In view of $\| T(x) \| ^2= <T(x),T(x)>
=<T^*(T(x)),x>$, the equality $T(x)=0$ is satisfied if and only if
$T^*(T(x))=0$, that is, ${\hat N}(T)={\hat N}(T^*T)$. Then from the
proved above it follows, that ${\hat R}(T^*T)=I-{\hat
N}(T^*T)=I-{\hat N}(T)={\hat R}(T^*)$.
\par In particular, for the selfadjoint element $T\in L_q(X)$ it is accomplished
the equality ${\hat R}(T)=I - {\hat N}(T)$.
\par {\bf 22.18. Proposition.} {\it If $\hat E$ and $\hat F$ are graded
operator projections for a Hilbert space $X$ over ${\bf K}=\bf H$ or
${\bf K}=\bf O$, then ${\hat R}({\hat E}+{\hat F})={\hat E}\vee \hat
F$, ${\hat R}({\hat E}{\hat F})={\hat E}-{\hat E}\wedge (I-{\hat
F})$ up to the isometric automorphism of the algebra $\bf K$.}
\par {\bf  Proof.} The subspaces ${\hat E}(X)$ and ${\hat
F}(X)$, ${\hat E}^{-1}(0)$ and ${\hat F}^{-1}(0)$ are the $\bf
K$-vector subspaces in $X$ (see Proposition 22.7). In view of $\|
{\hat E}(x) \| ^2 + \| {\hat F}(x) \| ^2=<{\hat E}(x),x>+<{\hat
F}(x),x>=<({\hat E}+{\hat F})(x),x>$ for each vector $x\in X$, then
$({\hat E}+{\hat F})(x)=0$ if and only if ${\hat E}(x)=0$ and ${\hat
F}(x)=0$. Thus, ${\hat N}({\hat E}+{\hat F})={\hat N}({\hat
E})\wedge {\hat N}({\hat F})$ $=(I-{\hat E})\wedge (I-{\hat F})$,
also due to Proposition 22.17 and \S 22.10 ${\hat R}({\hat E}+{\hat
F}) =I-(I-{\hat E})\wedge (I-{\hat F})={\hat E}\vee \hat F$ up to an
isometric automorphism of the algebra $\bf K$, since $$<{\hat
E}(x),y>= \sum_{p,q=0}^m<{\hat E}(x_p),y_q>i_p^*i_q=
\sum_{p,q=0}^m<x_p,{\hat E}(y_q)>i_p^*i_q=<x,{\hat E}(y)>$$ for each
$x, y\in X$. \par The graded operators of projections $I-\hat E$ and
${\hat E}\wedge (I-{\hat F})$ are mutually orthogonal. If $x\in X$
and ${\hat F}{\hat E}=0$, then ${\hat E}(x)=(I-{\hat F})({\hat
E}(x))\in ({\hat E}\wedge (I-{\hat F}))(X)$, therefore, $x=(I-{\hat
E})(x)+{\hat E}(x)\in (I-{\hat E}+{\hat E}\wedge (I-{\hat F}))(X)$.
Each vector $x\in (I-{\hat E}+{\hat E}\wedge (I-{\hat F}))(X)$ can
be expressed in the form $y+z$, where $y=(I-{\hat E})(y)$, $z={\hat
E}(z)=(I-{\hat F})(z)$, then ${\hat F}({\hat E} (x))={\hat F}({\hat
E}(y))+{\hat F}({\hat E}(z))={\hat F}({\hat E}(y))+{\hat F}({\hat
E}(z))={\hat F}({\hat E}((I-{\hat E})(y)))+ {\hat F}((I-{\hat
F})(z))=0$. Therefore, ${\hat N}({\hat F}{\hat E})=I-{\hat E}+{\hat
E}\wedge (I-{\hat F})$, that together with Proposition 22.17 gives
${\hat R}({\hat E}{\hat F})=I-{\hat N}({\hat F}{\hat E})={\hat
E}-{\hat E}\wedge (I-{\hat F})$.
\par {\bf 22.19. Proposition.} {\it If $ \{ y_a: a\in \Upsilon \} $ is an
orthonormal basis over $\bf K$ in a closed $\bf K$-vector subspace
$Y$ of a Hilbert space $X$ over ${\bf K}=\bf H$ or ${\bf K}=\bf O$,
also $\hat E$ is a graded operator of projection from $X$ onto $Y$,
then $${\hat E}(x)=\sum_{a\in \Upsilon }<x,y_a>y_a$$ for each $x\in
X$.}
\par {\bf Proof.} A Hilbert space $X$ has the decomposition
$X=X_0\oplus X_1i_1\oplus ... \oplus X_mi_m$, where $X_0,...,X_m$
are pairwise isomorphic Hilbert spaces over $\bf R$, $\{
i_0,i_1,...,i_m \} $ is the family of standard generators of the
algebra $\bf K$, $i_0=1$, $m=3$ for $\bf H$, $m=7$ for $\bf O$. Due
to the known structural theorem about Hilbert spaces over $\bf R$
the space $X_0$ is isomorphic with the Hilbert space $l_2(\omega
,{\bf R})$ over $\bf R$ with the orthonormed basis $e_j$ of the
cardinality $card (\omega )$, where the set $\omega $ can be taken
as some ordinal due to the Kuratowski-Zorn lemma, $j\in \omega $,
$e_j$ is the net of numbers $\{ x_{i,j}\in \{ 0, 1 \}: i\in \omega
\} $, $x_{j,j}=1$, $x_{i,j}=0$ for each $i\ne j\in \omega $. At the
same time $card (\omega )=\chi (X)$ is equal to the topological
character $X$. Then $X$ is isomorphic with $l_2(\omega ,{\bf K})$,
also in view of the fact that any orthonormed basis in $Y$ can be
enlarged up to the orthonormal basis of the entire space $X$, hence
up to the isometric isomorphism $Y$ is isomorphic with the subspace
$l_2(\nu ,{\bf K})$ for some ordinal $\nu $ such that $\nu \le
\omega $, $card (\nu )=card (\Upsilon )$. On the other hand, ${\bf
K}e_j=e_j{\bf K}$ for each $j$, therefore, ${\hat E}(x)\in Y$ has
the decomposition $${\hat E}(x)=\sum_{a\in \Upsilon }<{\hat
E}(x),y_a>y_a= \sum_{a\in \Upsilon }<x,{\hat E}(y_a)>y_a= \sum_{a\in
\Upsilon }<x,y_a>y_a.$$
\par {\bf 22.20. Proposition.} {\it If $\cal H$ is a selfadjoint $C^*$-algebra
of operators from $L_q(X)$ acting in a Hilbert space $X$ over ${\bf
K}=\bf H$ or ${\bf K}=\bf O$, also $\cal H$ is closed relative to
the weak operator topology, then the union and the intersection of
each family of graded operators of projections for $\cal H$ belong
to $\cal H$. There exists a graded operator of projection ${\hat
P}\in \cal H$ such that it is greater, then others graded operators
of projections in $\cal H$ and ${\hat P}A=A{\hat P}=A$ for each
$A\in \cal H$ such that $A^{-1}(0)$ and $A(X)$ are $\bf K$-vector
subspaces in $X$.}
\par {\bf Proof.} Due to Proposition 22.17 ${\hat R}(T^*T)=
{\hat R}(T^*)$, if $T(X)$ and $T^{-1}(0)$ are $\bf K$-vector
subspaces in $X$ and a graded operator of projection for such
selfadjoint operator from $\cal H$ belongs to $\cal H$,
consequently, for each such operator $T\in \cal H$, for which $T(X)$
and $T^{-1}(0)$ are $\bf K$-vector subspaces in $X$. In particular,
for an arbitrary graded operator of projection $\hat E$ from $\cal
H$ it can be taken as such operator $T=\hat E$. Due to Lemma 2.22.5
${\hat R}(A)$ is the limit relative to the strong operator topology
(an hence relative to the weak operator topology also) of the
sequence $A^{1/n}$, when $0\le A\le I$. If $A\in \cal H$, then
$A^{1/n}\in \cal H$, since the $C^*$-algebra is $\cal H$ closed
relative to the operator norm. Thus, ${\hat R}(A)\in \cal H$,
consequently, ${\hat R}(T)\in \cal H$ for each $T\in \cal H$, if
$T(X)$ and $T^{-1}(0)$ are $\bf K$-vector subspaces in $X$.
\par If $\hat E$ and $\hat F$ are graded operators of projections, then
${\hat R}({\hat E}+{\hat F})= {\hat E}\vee {\hat F}$ due to
Proposition 2.22.18. Thus, ${\hat E}\vee {\hat F}\in \cal H$, if
${\hat E}, {\hat F}\in \cal H$. If $\{ \mbox{ }_a{\hat E}: a \in
\Upsilon \} $ is a family of graded operators of projections from
$\cal H$, then their finite unions belong to $\cal H$ and form an
increasing net of bounded from above by the unit operator $I$. Due
to Lemma 2.22.4 this net has the supremum, which is its limit
relative to the strong operator topology. Thus, ${\hat E}\in \cal
H$, ${\hat E}=\vee \mbox{ }_a{\hat E}$. Due to \S 2.22.10 applied to
${\hat P}(X)$, where $\hat P$ is the union of all graded operators
of projections from $\cal H$, $\wedge \mbox{ }_a{\hat E}={\hat P}-
\vee _a ({\hat P}- \mbox{ }_a{\hat E})$. Due to the preceding proof
${\hat P}\in \cal H$ and $\wedge \mbox{ }_a{\hat E}\in \cal H$. The
Hilbert space $X$ has the decomposition $X=X_0\oplus X_1i_1\oplus
... \oplus X_mi_m$, where $X_0,...,X_m$ are pairwise isomorphic
Hilbert spaces over $\bf R$, $\{ i_0,i_1,...,i_m \} $ is the family
of standard generators of the algebra $\bf K$, $i_0=1$, $m=3$ for
$\bf H$, $m=7$ for $\bf O$. In view of the fact that $\hat P$
contains ${\hat R}(A)$ for each $A\in \cal H$ such that $A^{-1}(0)$
and $A(X)$ are $\bf K$-vector subspaces in $X$, then ${\hat P}A=A$
and ${\hat P}A^*=A^*$, that is, ${\hat P}_sA_q=(-1)^{\kappa (s)
\kappa (q)} A_q{\hat P}_s$ for each $s\ne q\in \{ i_0,i_1,...,i_m \}
$ from the standard set of generators of the algebra $\bf K$, where
$\kappa (1)=0$, $\kappa (i_p)=1$ for each $p\ge 1$, ${\hat
P}_sA_s=A_s{\hat P}_s$ for each $s\in \{ i_0,...,i_m \} $,
$A_{i_p}(x_p):=A(x_pi_p)$ and ${\hat P}_{i_p}(x_p):={\hat
P}(x_pi_p)$ for each $x_p\in X_p$. Thus, ${\hat P}A=A{\hat P}=A$.
\par {\bf 23. Definition.} Let $X$ and $Y$ be Hilbert spaces over ${\bf K}=\bf H$
or ${\bf K}=\bf O$ and $\cal B$ be a $\sigma $-algebra of Borel
subsets of the Hausdorff topological space $\Lambda $. We introduce
constants $S_{v,p,n}\in \bf K$ such that each operator $T\in
L_q(X,Y)$ has the decomposition:
\par $(i)$   $T(x) = \sum_{v=0}^m T_v(x)i_v$, moreover,
\par $(ii)$  $T_v(x)i_v = \sum_{n=0}^m S_{v,1,n}T(x)S_{v,2,n}$ for each $x\in X$,
$v=0,...,m $, \\
where $S_{v,1,n}=\gamma _{v,n}S_{v,2,n}$ with $\gamma _{v,n}\in \bf
R$, $S_{v,l,n}\in {\bf R}i_n$ for each $n=0,1,...,m $, where each
operator $T_v$ is defined on $X$ and takes values in $Y_v$,
$X=X_0i_0\oplus X_1i_1\oplus ... \oplus X_mi_m$, $X_0,...,X_m$ are
Hilbert pairwise isomorphic space over $\bf R$. These constants
arise from the identity:
\par $(iii)$ $z_0=(z+(2^r-2)^{-1} \{ -z + \sum
_{n=1}^{2^r-1}i_n(zi_n^*) \} )/2$,
\par $(iv)$ $z_v = (-zi_v+i_v(2^r-2)^{-1} \{ -z + \sum _{n=1}^{2^r-1}i_n(zi_n^*)
\} )/2$ \\
for each $v=1,2,...,2^r-1$ and $z\in \bf K$, where
$z=z_0i_0+z_1i_1+...+z_mi_m$, $z_0,...,z_m\in \bf R$, $r=2$ for $\bf
H$, $r=3$ for $\bf O$, $m=2^r-1$.
\par  Consider a mapping $\hat E$ defined on ${\cal
B}\times X^2$ and defining a unique $X$-projection valued spectral
measure $\hat E$ such that
\par $(i)$ $<{\hat E}(\delta )x;y>={\hat \mu }(\delta ;x,y)$
is regular (noncommutative) $\bf K$-valued measure for each $x, y
\in X$, where $\delta \in \cal B$. By our definition this means that
\par $(1)\quad {\hat \mu }(\delta ;x,y)= {\hat \mu }(x,y).\chi
_{\delta }$ and
\par $(2)\quad {\hat \mu }(x,y).f:=
\sum_{v,l,n}\int_{\Lambda }S_{v,1,n}f(\lambda )S_{v,2,n}{\tilde
i_v}i_l \mu _{i_v,i_l}(d\lambda ;x,y)$, \\
where $\chi _{\delta }$ is the characteristic function for $\delta
\in \cal B$, $\mu _{i_v,i_l}$ is the regular realvalued measure, $q,
n, l =0,1,...,m $, $p=1$ or $p=2$, $f$ is an arbitrary $\bf
K$-valued function on $\Lambda $, which is $\mu
_{i_v,i_l}$-integrable for each $v, l$;
\par $(3)\quad ({\hat E}_S(\delta ))^*=(-1)^{\kappa (S)}
{\hat E}_S(\delta ),$ where $({\hat E}_S(\delta ).e)x:= ({\hat
E}_S(\delta ))x=({\hat E}(\delta ).S)x$, $S=ci_q$, $q=0,1,...,m $,
$c=const\in \bf R$, $x\in X_0$;
\par $(4)\quad {\hat E}_{bS}=b{\hat E}_S$ for each $b\in \bf R$
and each pure vector $S=ci_q$;
\par $(5)\quad {\hat E}_{S_1S_2}(\delta \cap \gamma )=
{\hat E}_{S_1}(\delta ){\hat E}_{S_2}(\gamma )$ for each pure $\bf
K$ vectors $S_1$ and $S_2$ and every $\delta , \gamma \in \cal B$.
\par  Though from $(3,4)$ it follows, that ${\hat E}(\delta ).\lambda = \lambda
_0{\hat E}_{i_0}(\delta )+\lambda _1{\hat E}_{i_1}(\delta )
+...+\lambda _m{\hat E}_{i_m}(\delta ) =: {\hat E}_{\lambda }(\delta
)$, but in the general case it may happen that $({\hat E}(\delta
).\lambda )x\ne ({\hat E}(\delta ))\lambda x$, where $\lambda
_0,...,\lambda _m\in \bf R$.
\par {\bf 24. Theorem.} {\it Each quasicommutative $C^*$-algebra
$\cal A$ contained in $L_q(X)$ for a Hilbert space $X$ over ${\bf
K}=\bf H$ or ${\bf K}=\bf O$ is isometrically $*$-equivalent to the
algebra $C(\Lambda ,{\bf K})$, where $\Lambda $ is its spectrum.
Moreover, each isometric $*$-isomorphism $f\mapsto T(f)$ between
$C(\Lambda ,{\bf K})$ and $\cal A$ defines a unique $X$-projection
valued spectral measure $\hat E$ on ${\cal B}(\Lambda )$, such that
\par $(i)$ $<{\hat E}(\delta )x;y>={\hat \mu }(\delta ;x,y)$
is a regular $\bf K$-valued measure for each $x, y \in X$, where
$\delta \in \cal B$;
\par $(ii)$ ${\hat E}_{S_1}(\delta ).T(S_2f)=(-1)^{\kappa (S_1)
\kappa (S_2)}T(S_2{\hat E}_{S_1}(\delta ).f)$ for each $f\in
C(\Lambda ,{\bf R})$, $\delta \in \cal B$ and pure vectors $S_1\ne
S_2\in \bf K$, $\kappa (1):=0$, $\kappa (i_v)=1$ when $v=1,2,...,n$;
\par $(iii)\quad T(f)=\int_{\Lambda }{\hat E}(d\lambda ).f(\lambda )$
for each $f\in C(\Lambda ,{\bf K})$, moreover, $\hat E$ is $\sigma
$-additive in the strong operator topology.}
\par {\bf Proof.} Mention at first that $\Lambda $ is compact.
There exists the decomposition $C(\Lambda ,{\bf K})$ as in \S \S
2.12, 18. Each $\psi \in C^*_q(\Lambda ,{\bf K})$ has the
decomposition $\psi (f)=\psi _0(f)i_0+\psi _1(f)i_1+...+\psi
_m(f)i_m$, where $f\in C(\Lambda ,{\bf K})$. Moreover, $\psi _l(f) =
\psi _l(f_0i_0)+\psi _l(f_1i_1)+...+ \psi _l(f_mi_m)$, where $f_v\in
C(\Lambda ,{\bf R})$, $v, l= 0,1,...,m $. Then
$$(1)\quad \psi (f)=\sum_{v,n,l}\psi _l(S_{v,1,n}fS_{v,2,n})i_l,$$
where $v, n, l =0,1,...,m $.  Due to the Riesz representation
theorem $IV.6.3$ \cite{danschw}: \par $\psi _l(gi_v)=\int_{\Lambda
}g(\lambda ) \mu _{i_v,i_l}(d\lambda )$ \\
for each $g\in C(\Lambda ,{\bf R})$, where $\mu _{i_v,i_l}$ is the
$\sigma $-additive realvalued measure. Accomplishing the integration
componentwise of the matrix valued function gives $$(2)\quad \psi
(f)=\sum_{v,n,l} \int_{\Lambda } S_{v,1,n}f(\lambda
)S_{v,2,n}{\tilde i}_vi_l\mu _{i_v,i_l}(d\lambda ).$$ For $\psi
(f):=<T(f)x;y>$ for each $f\in C(\Lambda ,{\bf K})$ and marked
points $x, y\in X$ from $(2)$ it follows, that $$(3)\quad
<T(f)x;y>=\sum_{v,n,l} \int_{\Lambda } S_{v,1,n}f(\lambda
)S_{v,2,n}{\tilde i}_vi_l\mu _{i_v,i_l}(d\lambda ;x,y),$$ since
$|<T(f)x;y>| \le |f| |x| |y|$, consequently, $\mu _{i_v,i_l}(\delta
;xa,yb)=a\mu _{i_v,i_l}(\delta ;x,y)b$ for each $a, b\in \bf R$,
moreover, $$(4)\quad \sup_{\delta \in \cal B} (\sum_l |\sum_v
z_vi_v\mu _{i_v,i_l}(\delta ;x,y)|^2)^{1/2}\le |z| |x| |y|$$ for
each $z=z_0i_0+z_1i_1+...+z_mi_m\in \bf K$, since $|i_l|=1$.
\par From $(3)$ it follows that $\mu _{i_v,i_l}(\delta ;x,y)$
is $\bf R$-bihomogeneous and biadditive by $x, y$. If $f(\lambda
)\in {\bf R}i_v$ for $\mu $-almost every $\lambda \in \Lambda $ for
some $v=0,1,...,m$, then $T(f)=T((-1)^{\kappa (i_v)}{\tilde
f})=(-1)^{\kappa (i_v)}T(f)^*$, consequently,
$<T(f)x;y>=(-1)^{\kappa (i_v)}<T(f)y;x>^{\tilde .}$. Therefore, $\mu
_{i_v,i_l}(\delta ;x,y)=(-1)^{\kappa (i_v) \kappa (i_l)}\mu
_{i_v,i_l} (\delta ;y,x)$ for each $v\ne l=0,1,...,m $, $x, y\in X$.
\par {\bf 25. Definition.} An operator $T$ in a Hilbert space $X$
over ${\bf K}=\bf H$ or ${\bf K}=\bf O$ is called normal, if
$TT^*=T^*T$; $T$ is unitary, if $TT^*=I$ and $T^*T=I$; $T$ is
symmetrical, if $<Tx;y>=<x;Ty>$ for each $x, y \in {\cal D}(T)$, $T$
is selfadjoint, if $T^*=T$. Further, for $T^*$ it is supposed that
${\cal D}(T)$ is dense in $X$.
\par {\bf 26. Lemma.} {\it An operator $T\in L_q(X)$ for a Hilbert
space $X$ over ${\bf K}=\bf H$ or ${\bf K}=\bf O$ is normal if and
only if a minimal ($\bf K$-)subalgebra $\cal A$ in $L_q(X)$
containing $T$ and $T^*$ is quasicommutative.}
\par {\bf Proof.} Let $T$ be normal, then in
$X=X_0i_0\oplus X_1i_1\oplus ... \oplus X_mi_m$ it can be presented
in the form $T=T_0i_0+T_1i_1+...+T_mi_m$, where $Range (T_v)\subset
X_v$ and $T_v\in \cal A$ for each $v=0,1,...,m $. Therefore,
$<s(Ts)(x);s(Ts)(y)>=<s(T^*s)(x);s(T^*s)(y)>$ for each $s\in \bf K$
с $|s|=1$ and $Re (s)=0$, consequently,
$(s(Ts))(s(Ts))^*=(s(Ts))^*(s(Ts))$. The space $X$ is isomorphic
with $l_2(\upsilon ,{\bf K})$, in which $<x;y>=\sum_{b\in \upsilon
}\mbox{ }^b{\tilde x}\mbox{ }^by$, where $\upsilon $ is a set, $x=
\{ \mbox{ }^lx: \mbox{ }^lx\in {\bf K}, l\in \upsilon \} \in
l_2(\upsilon ,{\bf K})$. Then in $L_q(l_2(\upsilon ,{\bf K}))$ it is
accomplished $T^*=\tilde T$, moreover, ${\bar A}^*=A$ and ${\bar
B}^*=B$. Therefore, $TT^*=T^*T$ gives $AA^*=A^*A$, also the
automorphism $j: X\to X$ and the equality
$(s(Ts))(s(Ts))^*=(s(Ts))^*(s(Ts))$ with $s=\theta $,
$s=(i_v+j_w)/2$ and $s=(i_v+i_w)\theta /2$ for each $v<w \in \{
1,..., m \} $, $\theta =\exp (\pi i_1/4)$ lead to the pairwise
commuting $\{ T_0,T_1,...,T_m \} $.
\par Vise versa, if $\cal A$ quasicommutative, then
$\{ T_v: v=0,1,...,m \} $ are pairwise commuting, consequently,
$TT^*=T^*T$.
\par {\bf 27. Lemma.} {\it Let $T$ be a symmetrical operator and
$a\in {\bf K}\setminus {\bf R}e$, where ${\bf K}=\bf H$ or ${\bf
K}=\bf O$, then there exists $R(a;T)$ and $|x|\le 2 |R(a;T)x| /
|a-{\tilde a}|$ for each $x\in {\cal D}(T)$. Let $T$ be a closed
operator, then the sets $\rho (T)$, $\sigma _p(T)$, $\sigma _c(T)$
and $\sigma _r(T)$ are not intersecting and their union is the
entire $\bf K$. For a selfadjoint quasilinear operator $T$ there is
the inclusion $\sigma (T)\subset {\bf R}e$, moreover,
$R(a;T)^*=R(a^*;T)$.}
\par The {\bf Proof} is analogous to the complex case. Indeed, in the general
case for a quasilinear operator (not necessarily symmetrical) for
the decomposition of the components $\mu _{i_v,i_l}(d\lambda ;x,y)$
of the projection valued measure ${\hat \mu }(x,y).f$ defined in \S
23 as the sum of the point $(\mu _{i_v,i_l})_p$, absolutely
continuous $(\mu _{i_v,i_l})_{ac}$ and continuous singular $(\mu
_{i_v,i_l})_{sing}$ measures in accordance with the Lebesgue theorem
cited in \S 21 gives $L^2({\bf K},\mu _{i_v,i_l},{\bf K})={\cal
H}_{p,v,l}\oplus {\cal H}_{ac,v,l}\oplus {\cal H}_{sing,v,l}$, where
${\cal H}_{p,v,l}:=L^2({\bf K},(\mu _{i_v,i_l})_p,{\bf K})$,  ${\cal
H}_{ac,v,l}:=L^2({\bf K},(\mu _{i_v,i_l})_{ac},{\bf K})$ and ${\cal
H}_{sing,v,l}:=L^2({\bf K},(\mu _{i_v,i_l})_{sing},{\bf K})$. At the
same time for a symmetrical operator there is the inclusion $\sigma
(T)\subset \bf R$ and due to the relations given in \S 23 for
components of the projection valued measure the supports of all
these measures for different $v,l$ are consistent and are contained
in $\bf R$.
\par {\bf 28. Theorem.} {\it For a selfadjoint quasilinear operator $T$
there exists a uniquely defined regular countably additive
selfadjoint spectral measure $\hat E$ on ${\cal B}({\bf K})$,
${\hat E}|_{\rho (T)}=0$ such that \\
$(a)$ ${\cal D}(T):= \{ x: x\in X; \int_{\sigma (T)}<({\hat E}(dz).z^2)x;x>
<\infty \} $ and \\
$(b)$ $Tx=\lim_{n\to \infty }\int^n_{-n}({\hat E}(dz).z)x$,
$x\in {\cal D}(T)$.}
\par {\bf Proof.} Due to Proposition
2.20, Note 2.21 and Definition 2.16(5) the space ${\cal D}(T)$ is
$\bf K$-vector. Use the Lemma 2.27 and take a marked element $q\in
\{ i_0,i_1,...,i_m \} $, then $h(z):=(q-z)^{-1}$ is the
homeomorphism of the sphere $S^m := \{ z\in {\bf K}: |z|=1 \} $ and
for $A:=(q-z)(R(z;T)(q-z))+(q-z)I$ for each $z\in \rho (T)\setminus
\{ q \} $ is accomplished the identity $(hI-R(q;T))A=I$. If $z=q$,
then $h=\infty $, consequently, $h\notin \sigma (R(q;T))$. Let $0\ne
h\in \rho (R(q;T))$, then there exists $B:= R(q;T)A$, where
$A:=(hI-R(z;T))^{-1}$, consequently, $B$ is bijective, ${\cal R}
(B)={\cal D}(T)$ and $(zI-T)B=(z-q)I$, that is, $z\in \rho (T)$. For
$h=0\in \rho (R(q;T))$ the operator $R(q;T)^{-1}=(hI-T)$ is bounded
and is defined everywhere operator and this case is considered in
Theorem 2.24. For each $\delta \in {\cal B}({\bf K})$ put ${\hat
E}(\delta ):={\hat E}^1(h(\delta ))$, where ${\hat E}^1$ is the
decoposition of the unity for a normal operator $R(q;T)$, then the
end of the proof is analogous to Theorem XII.2.3 \cite{danschw}.
\par {\bf 29. Note and Definition.} A unique spectral measure,
related with a selfadjoint quasilinear operator $T$ is called a
decomposition of the unity for $T$. For a $\bf K$-valued Borel
function $f$ defined $\hat E$-almost everywhere on $\bf K$ $f(T)$ is
defined by the relations: \par ${\cal D}(f(T)):=\{ x:$ $\mbox{there
exists}$ $\lim_nf_n(T)x \} $, \\
where $f_n(z):=f(z)$ for $|f(z)|\le
n$; $f_n(z):=0$ while $|f(z)|>n$; $f(T)x:=\lim_nf_n(T)x$, $x\in
{\cal D}(f(T))$, $n\in \bf N$.
\par {\bf 30. Theorem.} {\it Let $\hat E$ be a decomposition of the unity
for a selfadjoint quasilinear operator $T$ and $f$ is from \S 2.29.
Then $f(T)$ is a closed quasilinear operator with an everywhere dense domain of its
definition, moreover: \\
$(a)$ ${\cal D}(f(T))=\{ x: \int^{\infty }_{-\infty }|f(z)|^2
<{\hat E}(dz)x;x> <\infty \} $; \\
$(b)$ $<f(T)x;y>=\int^{\infty }_{-\infty }<{\hat E}(dz).f(z)x;y>$,
$x\in {\cal D}(f(T))$; \\
$(c)$ $|f(T)x|^2=\int^{\infty }_{-\infty }|f(z)|^2<{\hat E}(dz)x;x>$,
$x\in {\cal D}(f(T))$; \\
$(d)$ $f(T)^*={\tilde f}(T)$; \\
$(e)$ $R(q;T)=\int^{\infty
}_{-\infty } {\hat E}(dz).(q-z)$, $q\in \rho (T)$.}
\par {\bf Proof.} Take $f_n$ from \S 2.29
and $\delta _n:= \{ z: |f(z)|\le n \} $. Then $|f(T)x|^2=
\lim_n|f_n(T)x|^2=\int^{\infty}_{-\infty }|f(z)|^2 <{\hat
E}(dz)x;x>$ for each $x\in {\cal D}(f(T))$, from this it follows
$(c)$, also the closedness of $f(T)$ and $(a)$ are verified
analogously to the complex case. To the noncommutative measure $\hat
\mu $ given on the algebra $\Upsilon $ of subsets of the set $\cal
S$ there corresponds a quasilinear operator with values in $\bf K$
and due to Proposition 2.20 this measure is completely characterized
by the $\bf R$-valued measures $\mu _{i_v,i_l}$, such that $\mu
_{i_v,i_l}(f_v)= {\hat \mu }(f_v){\tilde i}_l$ for each $\hat \mu
$-integrable $\bf K$-valued function $f$ with components $f_v$,
where $v,l =0,1,...,m$. Then it can be defined the variation
$V({\hat \mu },U):=\sup_{W_l\subset U} \sum_l|{\hat \mu } (\chi
_{W_l})|$ by all finite disjoint systems $\{ W_l \} $ of subsets
$W_l\in \Upsilon $ in $U$ with $\bigcup_lW_l=U$. If $\hat \mu $ is
bounded, then it is the quasilinear operator of the bounded
variation with $V({\hat \mu },{\cal S})\le 2^{2^{m+2}} \sup_{U\in
\Upsilon } |{\hat \mu }(\chi _U)|$, moreover, $V({\hat \mu },*)$ is
additive on $\Upsilon $. The function $f$ we call $\hat \mu
$-measurable, if each $f_v$ is $\mu _{i_v,i_l}$-measurable for each
$v$ and $l=0,1,...,m .$ The space of all ${\hat \mu }$-measurable
$\bf K$-valued functions $f$ with $V({\hat \mu
},|f|^p)^{1/p}=:|f|_p<\infty $ we denote by $L^p({\hat \mu })$ while
$0<p<\infty $, also $L^{\infty }({\hat \mu })$ is the space of all
$f$ for which there exists $|f|_{\infty }:=ess_{V ({\hat \mu
},*)}-\sup |f|<\infty .$ In details we write $L^p({\cal S},{\Upsilon
},{\hat \mu },{\bf K})$ instead of $L^p({\hat \mu })$. A subset $W$
in $\cal S$ we call $\cal \mu $-zero-set, if $V^*({\hat \mu },W)=0$,
where $V^*$ is the extension of the complete variation $V$ by the
formula $V^*({\hat \mu },A):=\inf_{{\Upsilon }\ni F\supset A}V({\hat
\mu },F)$ for $A\subset \cal S$. A noncommutative measure $\hat
\lambda $ on $\cal S$ we call an absolutely continuous relative to
$\hat \mu $, if $V^*({\hat \lambda },A)=0$ for each subset $A\subset
\cal S$ with $V^*({\hat \mu },A)=0$. A measure $\hat \mu $ we call
positive, if each $\mu _{m,n}$ is nonnegative and $\sum_{m,n}\mu
_{m,n}$ is positive. The using of components $\mu _{m,n}$ and the
classical Radon-Nikodym theorem (see Theorems III.10.2,10.7
\cite{danschw}) leads to the following noncommutative variants.
\par $(i).$ If $({\cal S},{\Upsilon },{\hat \mu })$ is a space with
a $\sigma $-finite positive noncommutative $\bf K$-valued measure
$\hat \mu $, also $\hat \lambda $ is an absolutely continuous
relative to $\hat \mu $ bounded noncommutative measure defined on
$\Upsilon $, then there exists a unique $f\in L^p({\cal S},{\Upsilon
},{\hat \mu },{\bf K})$, such that ${\hat \lambda }(U)={\hat \mu
}(f\chi _U)$ for each $U\in \Upsilon $, moreover, $V({\hat \mu
},{\cal S})=|f|_1$.
\par $(ii).$ If $({\cal S},{\Upsilon },{\hat \mu })$ is a
space with bounded noncommutative $\bf K$-valued measure $\hat \mu
$, also $\hat \lambda $ is absolutely continuous relative to $\hat
\mu $ noncommutative measure defined on $\Upsilon $, then there
exists a unique $f\in L^1({\hat \mu })$, such that ${\hat \lambda }
(U)={\hat \mu }(f\chi _U)$ for each $U\in \Upsilon $. Due to $(ii)$
there exists a Borel measurable function $\phi $, such that ${\hat
\nu }(\delta ):={\hat \mu }_{x,y}(\phi \chi_{\delta })= <{\hat
E}(\phi \chi _{\delta })x;y>$ for each $\delta \in {\cal B}({\bf
R})$. Due to $(i)$ $|\phi (z)|=1$ $\hat \nu $-almost everywhere.
Consider $f^1(z):=|f(z)|\phi (z)$, then due to $(a)$ ${\cal
D}(f^1(T))={\cal D}(f(T))$ and $<f^1(T)x;y>= \int^{\infty }_{-\infty
}|f(z)|{\hat \nu }(dz)$. Therefore, $$<f(T)x;y>=\lim_n\int_{\delta
_n}<{\hat E}(dz).f(z)x;y>=\int^{\infty }_{-\infty }<{\hat
E}(dz).f(z)x;y>$$ and from this it follows $(b)$.
\par $(d)$. From ${\hat E}_S^*=(-1)^{\kappa (S)}{\hat E}_S$
for each  $S=cs$, $0\ne c\in \bf R$, $s\in \{ i_0,...,i_m \} $, it
follows, that ${\hat E}.{\tilde f}={\hat E}^*.f$. Take $x,y\in {\cal
D}({\tilde f}(T))={\cal D}(f(T))$, then $$<{\tilde
f}(T)x;y>=\int^{\infty }_{-\infty } <{\hat E}(dz).{\tilde
f}(z)x;y>=<x;f(T)y>,$$ consequently, ${\tilde f}(T)\subset f(T)^*$.
If $y\in {\cal D}(f(T)^*)$, then for each $x\in X$ and $t\in \bf N$:
${\tilde f}_t(T)y := {\hat E}(\delta _t) .f(T)^*y$ converges to
$f(T)^*y$ for $t\to \infty $, consequently, $y\in {\cal D}({\tilde
f}(T))$. Due to Theorem 2.24 the statement $(e)$ follows from the
fact that $\mbox{ }_t{\hat E}(\delta ):={\hat E}(\delta _t\cap
\delta )$ is the decomposition of the unity for the bounded
restriction $T|_{X_t}$, where $X_t:={\hat E}(\delta _t)X$.
\par {\bf 31. Theorem.} {\it A bounded normal operator
$T$ on a Hilbert space over ${\bf K}=\bf H$ or ${\bf K}=\bf O$ is
unitary, Hermitian or positive if and only if $\sigma (T)$ is
contained in $S^m:= \{ z\in {\bf K}: |z|=1 \} $, $\bf R$ or
$[0,\infty )$ respectively.}
\par {\bf Proof.} Due to Theorem 2.24 the equality $T^*T=TT^*=I$
is equivalent to $z{\tilde z}=1$ for each $z\in \sigma (T)$. If
$\sigma (T)\subset [0,\infty )$, then $<Tx;x>=\int_{\sigma (T)}
<{\hat E}(dz).zx;x>\ge 0$ for each $x\in X$. The final part of the
proof is analogous to the complex case, using the technique given
above.
\par {\bf 32. Definition.} A family $ \{ T(t): 0\le t\in {\bf R} \} $
of bounded quasilinear operators in a vector space $X$ over ${\bf
K}=\bf H$ or ${\bf K}=\bf O$ is called a strongly continuous
semigroup, if \par $(i)$ $T(t+q)=T(t)T(q)$ for each $t, q\ge 0$;
\par $(ii)$ $T(0)=I$; \par $(iii)$ $T(t)x$ is a continuous function by
$t\in [0,\infty )$ for each $x\in X$.
\par {\bf 33. Theorem.} {\it For each strongly continuous
semigroup $ \{ U(t): 0\le t\in {\bf R} \} $ of a unitary quasilinear
operators in a Hilbert space $X$ over ${\bf K}=\bf H$ or ${\bf
K}=\bf O$ there exists a unique selfadjoint quasilinear operator $B$
in $X$, such that $U(t)=\exp (tMB)$, where $M\in \bf K$, $|M|=1$,
$Re (M)=0$.}
\par {\bf Proof.} If $\{ T(t): 0\le t \} $ is a semigroup
continuous in the unform topology (see also the complex case in
Theorem VIII.1.2 \cite{danschw}), then there exists a bounded
operator $A$ on $X$ such that $T(t)=\exp (tA)$ for each $t\ge 0$. If
$Re (z):=(z+{\tilde z})/2 >|A|$, then $|\exp (-t(zI-A))|\le \exp
(t(|A|- Re (z))\to 0$ for $t\to \infty $. For such $z\in \bf K$ due
to the Lebesgue theorem: $(zI-A)\int_0^{\infty }\exp
(-t(zI-A))dt=I$, also by Lemma 2.3 there exists
$R(z;A)=\int_0^{\infty }\exp (-t(zI-A))dt$. For each $\epsilon
>0$ let $A_{\epsilon }x:= (T(\epsilon )x-x)/\epsilon $, where $x\in
X$, for which there exists $\lim_{0<\epsilon \to 0}A_{\epsilon }x$,
the set of all such $x$ we denote by ${\cal D}(A)$. Evidently,
${\cal D}(A)$ is the $\bf K$-vector subspace in $X$. Take some
infinitesimal quasilinear operator $Ax:=\lim_{0<\epsilon \to
0}A_{\epsilon }x$. Considering  $\bf K$ as the Banach space over
$\bf R$, we get the analogs of Lemmas 3,4,7, Corollaries 5, 9 and
Theorem 10 from \S VIII.1 \cite{danschw}, moreover, ${\cal D}(A)$ is
dense in $X$, also $A$ is the closed quasilinear operator on ${\cal
D}(A)$. Let $w_0:= \lim_{t\to \infty } ln (|T(t)|)/t$ and $z\in \bf
K$ with $Re (z)>w_0$. For each $w_0<\delta <Re(z)$ due to Corollary
VIII.1.5 \cite{danschw} there exists as constant $M>0$ such that
$|T(t)|\le M\exp (\delta t)$ for each $t\ge 0$. Then there exists
$R(z)x:=\int_0^{\infty }\exp (-t(zI-A))xdt$ for each $x\in X$ and
$Re (z)>w_0$, consequently, $R(z)x\in {\cal D}(A)$. Let $T_z$ be a
quasilinear operator corresponding to $z^{-1}A$ instead of $T$ for
$A$, where $0\ne z\in \bf K$, moreover, ${\cal D}(A)={\cal
D}(z^{-1}A)$. Then $$z^{-1}A\int_0^{\infty }\exp (-t(I-z^{-1}A))xdt=
\int_0^{\infty }\exp (-t(I-z^{-1}A)z^{-1}Axdt,$$ consequently,
$R(z)(zI-A)x=x$ for each $x\in {\cal D}(A)$ and $R(z)=R(z;A)$. Thus,
$$R(z;A)x=\int_0^{\infty }\exp (-t(zI-A))xdt$$ for each $z\in \rho
(A)$ and $x\in X$.
\par With the help of $2.7.(iii)$ for the quasilinear operator $A$ there exists
the quasilinear operator $B$ such that $A=MB$, where $M\in \bf K$,
$|M|=1$, $Re (M)=0$. In view of $U(t)U(t)^*=U(t)^*U(t)=I$, then $A$
commutes with $A^*$ and $\exp (t(A+A^*))=I$. From
$R(z;B)^*=R({\tilde z},B)$ it follows, that we can choose $B=B^*$.
If ${\hat E}$ is the decomposition of the unity for $B$ and
$V(t):=\exp (MtB)$, then by Theorem 2.30 $$<V(t)x;y>=\int^{\infty
}_{-\infty } <{\hat E}(dz).\exp (Mtz)x;y>,$$ then due to the Fubini
theorem $$\int_0^{\infty }<V(t).\exp (-bt)x;y>dt=\int_0^{\infty }
\int^{\infty }_{-\infty }<{\hat E}(dz).\exp (-(b-Mz)t)x;y>dt$$
$$=\int^{\infty }_{-\infty }<{\hat E}(dz).(b-Mz)^{-1}x;y>=R(b;MB)x;y>$$
for each $b\in \bf K$ with $Re (b)>0$. Therefore,
$$\int_0^{\infty }<V(t).\exp (-bt)x;y>dt= \int_0^{\infty }<U(t).\exp
(-bt)x;y>dt$$ for $Re (b)>0$. Due to Lemma $VIII.1.15$ $<V(t).\exp
(-\epsilon t)x;y>= <U(t).\exp (-\epsilon t)x;y>$ for each $t\ge 0$
and $Re(b)>0$, consequently, $U(t)=V(t)$.
\par {\bf 34. Notations.} Let $X$ be a $\bf K$-vector
locally convex space. Consider left and right two sided $\bf
K$-vector spans of the family of vectors $\{ v^a: a\in {\bf A} \} $,
where $$span_{\bf K}^l \{ v^a: a\in {\bf A} \} := \{ z\in X:
z=\sum_{q_a\in \bf K; a\in \bf A} q_av^a \} ;$$ $$span_{\bf K}^r \{
v^a: a\in {\bf A} \} := \{ z\in X: z=\sum_{q_a\in \bf K; a\in \bf A}
v^aq_a \} ;$$ $$span_{\bf K} \{ v^a: a\in {\bf A} \} := \{ z\in X:
z=\sum_{q_a, r_a \in \bf K; a\in \bf A} \{ q_av^ar_a \}_{q(3)} \}
.$$
\par {\bf 35. Lemma.} {\it In the notation of \S 2.34
$span_{\bf K}^l \{ v^a: a\in {\bf A} \} = span_{\bf K}^r \{ v^a:
a\in {\bf A} \} = span_{\bf K} \{ v^a: a\in {\bf A} \} .$ }
\par {\bf Proof.}  Due to the continuity of the addition and multiplication on scalars of
vectors in $X$ and using the convergence of nets of vectors it is
sufficient to prove the statement of this lemma for a finite set
$\bf A$. Then the space $Y:=span_{\bf K} \{ v^a: a\in {\bf A} \} $
is finite dimensional over $\bf K$ and evidently left and right $\bf
K$-vector spans are contained in it. Then in $Y$ it can be chosen a
basis over $\bf K$ and each vector can be written in the form $v^a=
\{ v^a_1,...,v^a_n \} $, where $n\in \bf N$, $v^a_s\in \bf K$. On
the other hand, the algebra $\bf K$ is alternative, also
$X=X_0i_0\oplus X_1i_1\oplus ...\oplus X_mi_m$, where $X_0,...,X_m$
are pairwise isomorphic $\bf R$-linear locally convex spaces.
Therefore, $span_{\bf K}^l \{ v^a: a\in {\bf A} \} \cap span_{\bf
K}^r \{ v^a: a\in {\bf A} \} \supset span_{\bf K} \{ v^a: a\in {\bf
A} \} $, that together with the inclusion $span_{\bf K}^r \{ v^a:
a\in {\bf A} \} \cup span_{\bf K}^l \{ v^a: a\in {\bf A} \} \subset
span_{\bf K} \{ v^a: a\in {\bf A} \} $ proved above leads to the
statement of this lemma.
\par {\bf 36. Lemma.} {\it Let $X$ be a Hilbert space over ${\bf K}=\bf H$
or ${\bf K}=\bf O$, also $\bf X_{\bf R}$ be the same space
considered over the field $\bf R$. A vector $x\in X$ is orthogonal
to a $\bf K$-vector subspace $Y$ in $X$ relative to the $\bf
K$-valued scalar product in $X$ if and only if $x$ is orthogonal to
$Y_{\bf R}$ relative to the scalar product in $X_{\bf R}$. The space
$X$ is isomorphic to the standard Hilbert space $l_2(\alpha ,{\bf
K})$ over $\bf K$ of converging relative to the norm of sequences or
nets $v=\{ v^a: a\in \alpha \} $ with the scalar product $<v;w>:=
\sum_a{\tilde v}^aw_a$, moreover, $card (\alpha )\aleph _0=w(X)$,
where $card (\alpha )$ is the cardinality of the set $\alpha $,
$\aleph _0=card ({\bf N})$.}
\par {\bf Proof.} Due to Lemma 2.35 and by the transfinite induction
in $Y$ there exists a $\bf K$-linearly independent system of vectors
$\{ v^a: a\in {\bf A } \} $, such that $span_{\bf K}^r \{ v^a: a\in
{\bf A} \} $ is everywhere dense in $Y$. In another words in $Y$
there exists a Hamel basis over $\bf K$. A vector $x$ by the
definition is orthogonal to $Y$ if and only if $<v;x>=0$ for each
$v\in Y$, that is equivalent to $<v^a;x>=0$ for each $a\in \bf A$.
The space $X$ is isomorphic to the direct sum $X_0\oplus
X_1i_1\oplus ...\oplus X_mi_m$, where $X_0$, ..., $X_m$ are the
pairwise isomorphic Hilbert spaces over $\bf R$, also $X_{\bf
R}=X_0\oplus X_1\oplus ... \oplus X_m$. The scalar product $<x;y>$
in $X$ then can be written in the form
\par $(i)$ $<x;y>=\sum_{v,n=0}^m <x_v;y_n>{\tilde i}_vi_n$, \\
where $<x_v;y_n>\in \bf R$ due to $2.16$. Then the scalar product
$<x;y>$ in $X$ induces the scalar product
\par $(ii)$ $<x;y>_{\bf R}:=\sum_{v=0}^m <x_v;y_v>$ \\
in $X_{\bf R}$. Therefore, from the orthogonality of $x$ to the
subspace $Y$ relative to $<x;y>$ it follows the orthogonality of $x$
to the subspace $Y_{\bf R}$ relative to $<x;y>_{\bf R}=(x,y)$. Due
to Lemma 2.35 from $y\in Y$ it follows, that $y_vi_v\in Y$ for each
$v=0,1,...,m $. Then from $<x;y_v>_{\bf R}=0$ for each $y\in Y$ and
$v$ due to 2.16.(5) it follows, that $<x;y>=0$ for each $y\in Y$.
Then by the theorem about transfinite induction \cite{eng} in $X$
there exists the orthogonal basis over $\bf K$, in which every
vector can be presented in the form of the converging series of left
(or right) $\bf K$-vector combinations of basic vectors. For each
$x\in X$ due to the fact that the space $X$ is normed it follows,
that the base of neighborhoods of $x$ is countable, also for the
topological density it is accomplished the equality $d(X)=card
(\alpha )\aleph _0$, since $\bf K$ is separable, therefore,
$w(X)=d(X)$. From this it follows the last statement of this lemma.
\par {\bf 37. Lemma.} {\it For each quasilinear operator $T$ in
a Hilbert space $X$ over ${\bf K}=\bf H$ or ${\bf K}=\bf O$ an
adjoint operator $T^*$ in $X$ relative to the $\bf K$-scalar product
coincides with an adjoint operator $T^*_{\bf R}$ in $X_{\bf R}$
relative to the $\bf R$-valued scalar product in $X_{\bf R}$.}
\par {\bf Proof.} Let ${\cal D}(T)$ be the domain of the definition of
the operator $T$, which is dense in $X$. Due to Formula
$2.36.(i,ii)$ and the existence of the automorphisms $z\mapsto zi_v$
in $\bf K$ for each $v=0,1,...,m $ it follows, that the continuities
of $<Tx;y>$ and $<Tx;y>_{\bf R}$ by $x\in {\cal D}(T)$ are
equivalent, therefore, due to Lemma 2.35 the family of all $y\in X$,
for which $<Tx;y>$ is continuous by $x\in {\cal D}(T)$ forms the
$\bf K$-linear subspace in $X$ and it is also selfadjoint relative
to $<Tx;y>_{\bf R}$, that is the domain of the definition  ${\cal
D}(T^*)$ of the operator $T^*$. Then the adjoint operator $T^*$ is
defined by the equality $<Tx;y>=:<x;T^*y>$, also $T^*_{\bf R}$ is
given by the way of $<Tx;y>_{\bf R}=<x;T^*_{\bf R}y>_{\bf R}$, where
$x\in {\cal D}(T)$, also $y\in {\cal D}(T^*)$. Due to Formula
$2.36.(i,ii)$ $<x_v;(T^*y)_v>=<x_v;(T^*y)_v>_{\bf R}$ for each $x\in
{\cal D}(T)$, $y\in {\cal D}(T^*)$ and $v=0,1,...,m $. In view of
Proposition 2.20 and Lemma 2.35 ${\cal D}(T)$ and ${\cal D}(T^*)$
are $\bf K$-vector spaces, then an automorphisms of the algebra $\bf
K$ given above lead to $T^*=T^*_{\bf R}$.
\par {\bf 38. Definition.} A bounded quasilinear operator $P$ in a
Hilbert space $X$ over $\bf K$ is called a partial $\bf R$- (or $\bf
K$-) isometry, if there exists a closed $\bf R$- (or $\bf K$-)
vector subspace $Y$ such that $\| Px \| = \| x \| $ for $x\in Y$ and
$P(Y^{\perp }_{\bf R}) = \{ 0 \} $ (or $P(Y^{\perp })= \{ 0 \} $)
respectively, where $Y^{\perp }:=\{ z\in X: \quad <z;y>=0 \quad
\forall y\in Y \} $, $Y^{\perp }_{\bf R}:=\{ z\in X_{\bf R}: \quad
<z;y>_{\bf R}=0 \quad \forall y\in Y \} $.
\par {\bf 39. Theorem.} {\it If $T$ is a closed quasilinear
operator in a Hilbert space $X$ over $\bf K$, then $T=PA$, where $P$
is a partial $\bf R$-isometry on $X_{\bf R}$ with the initial domain
$cl (Range (T^*))$, also $A$ is a selfadjoint quasilinear operator
such that $cl (Range (A))=cl (Range (T^*))$. If $T$ is $\bf
K$-vectorial (that is, left or right $\bf K$-linear), then $P$ is a
partial $\bf K$-isometry.}
\par {\bf Proof.} Due to the spectral theorem
2.28 a selfadjoint quasilinear operator $T$ is positive if and only
if its spectrum is contained in $\sigma (T)\subset [0, \infty )$
(see also Lemma XII.7.2 \cite{danschw}). In the algebra $\bf K$ each
polynomial has a root (see Theorem 3.17 \cite{luoystoc}). Therefore,
if $T$ is a positive selfadjoint quasilinear operator, then there
exists a unique positive quasilinear operator $A$, such that $A^2=T$
(see also Lemma XII.7.2 \cite{danschw}). Then there exists a
positive square root $A$ of the operator $T^*T$. At the same time
$A$ is right or left $\bf K$-linear, if $T$ is right or left $\bf
K$-linear. Put $SAx=Tx$ for each $x\in {\cal D}(T^*T)$, also $V$ let
it be an isometric extension of $S$ on $cl (Range (A))$. The space
$cl (Range (A))$ is $\bf R$-linear. If $A$ in addition left (or
right) $\bf K$-linear, then $cl (Range (A))$ is the $\bf K$-vector
subspace due to Lemma 2.35. In view of Lemma 2.36 there exists the
perpendicular projection $\hat E$ from $X$ on $cl (Range (A))$,
moreover, $\hat E$ is right or left $\bf K$-linear, if $cl (Range
(A))$ is the $\bf K$-vector subspace. Then put $P=V\hat E$. From
$<Ax;Ax>=<Tx;Tx>$ for each $x\in {\cal D}(T^*T)$ it follows that
$PAx=Tx$ for each $x\in {\cal D}(T^*T)$. The remainder of the proof
can is accomplished analogously to the proof of Theorem $XII.7.7$
\cite{danschw} with the help of Lemmas 2.35-37.
\par {\bf 40. Note and Definition.} Apart from the case of the field
of complex numbers $\bf C$ nontrivial polynomials of quaternion or
octonion variables can have roots , which are not points, but also
closed submanifolds in $\bf K$ with the codimension over $\bf R$
from $0$ up to $m$ (see \cite{luoyst,luoystoc}).
\par A closed subset $\lambda \subset \sigma (T)$ is called
an isolated subset of a spectrum, if there exists a neighborhood $U$
of a subset $\lambda $ such that $\sigma (T)\cap U=\lambda $. An
isolated subset $\lambda $ of a spectrum $\sigma (T)$ is called a
pole of a spectrum (of an order $p$), if $R(z;T)$ has a zero on
$\lambda $ (of an order $p$, that is, each $z\in \lambda $ is zero
of an order $0<p(z)\le p$ for $R(z;T)$ and $\max_{z\in \lambda
}p(z)=p$). A subset $\lambda $ is clopen (closed and open
simultaneously)  in $\sigma (T)$ is called the spectral set. Let
$\eta _1$ be a closed rectifiable path in $U$, encompassing $\lambda
$ and not intersecting with $\lambda $, characterized by a vector
$M_1\in \bf K$, $|M_1|=1$, $M_1+{\tilde M}_1=0$ (see Theorem 3.22
\cite{luoyst} and Theorem 3.21 \cite{luoystoc}), denote by
$$(i)\quad \phi _n(z,T):=(2\pi )^{-1}
\{ \int_{\eta _1} R(\zeta ;T)((\zeta -a)^{-1}(z-a))^n(\zeta -a)^{-1}
d\zeta )M_1^{-1} \} ,$$ where $\eta _1\subset B({\bf
K},a,R)\setminus B({\bf K},a,r)$, $B({\bf K},a,R)\subset U$,
$\lambda \subset B({\bf H},a,r)$, $0<r<R<\infty $. We say that an
index $\lambda $ is equal to $p$ if and only if there exists a
vector $x\in X$ such that
$$(ii)\quad (zI-T)^{s_1}((v_1{\hat E}(\delta (z);T)v_2)...((zI-T)^{s_m}(v_{2m-1}
{\hat E}(\delta (z);T)(x)))...)=0$$ for each $z\in \lambda $ and
each $0\le s_n\in \bf Z$ with $s_1+...+s_m=p$ and each
$v_1,..,v_{2m-1}$, where $v_1=v_1(\delta ,T)\in \bf
K$,...,$v_{2m-1}=v_{2m-1}(\delta ,T) \in \bf K$, $m\in \bf N$,
$\delta :=\delta (z)\ni z$, $\delta (z)\in {\cal B}(\lambda )$, also
Expression in $(ii)$ is different from zero for some $z\in \lambda $
and $s_1,...,s_m$ with $s_1+...+s_m=p-1$.
\par {\bf 41. Theorem.} {\it A subset $\lambda $ is a pole of an
order $p$ of a quasilinear operator $T\in L_q(X)$ for $U=B({\bf
K},\alpha ,R')$, $0<R<R'<\infty $ in Definition 2.40, where
$0<r<\infty $, if and only if $\lambda $ has an index $p$.}
\par {\bf Proof.} Choose with the help of the homotopicity relative to
$U\setminus \lambda $ closed paths (loops)  $\eta _1$ and $\eta _2$
homotopic to $\gamma _1$ and $\gamma _2$, moreover, $\inf _{\theta }
|\eta _1(\theta )| > \sup_{\theta }|\eta _2(\theta )|$, where
$\gamma _1$ and $\gamma _2$ are chosen as in Theorem $3.22$
\cite{luoyst} or as in Theorem 3.21 \cite{luoystoc} (see also
Theorem $3.9$ there), $\theta \in [0,1]$. Due to these theorems  the
Loran decomposition over $\bf K$ of $R(z;T)$ in the neighborhood
$B({\bf K},a,R)\setminus B({\bf K},a,r)$ has the form
$$R(z;T)=\sum_{n=0}^{\infty }(\phi _n(z,T)+ \psi _n(z,T)),$$ where
$\phi _n$ is given by the Formula $2.40.(i)$, also
$$(i)\quad \psi _n(z,T):=(2\pi )^{-1}
\{ \int_{\eta _2} R(\zeta ;T)(z-a)^{-1}((\zeta -a)(z-a)^{-1})^n
d\zeta )M_2^{-1} \} .$$ If $\lambda $ is a pole of order $p$, then
$\phi _p=0$ and $\phi _{p-1}\ne 0$, therefore, there exists $x\in X$
such that
$$(ii)\quad \phi _{p(z)}(z,T)x=0 \mbox{ for each }z\in \lambda ,
\mbox{ also }$$
$$(iii)\quad \phi _{p-1}(z,T)x\ne 0\mbox{ for some }z\in \lambda .$$
Analogous decompositions are true with the corresponding $\phi _n$
for the product \par $f(T)(R(z;T)g(T))$, where $f$ and $g$ are
holomorphic functions on a neighborhood of $\sigma (T)$ in $\bf K$
different from zero everywhere on $\lambda $. Each function $\phi
_n$ for $R(z;T)$ can be approximated with any accuracy in the strong
operator topology in the form of left $\bf K$-linear combinations of
functions from $2.40.(ii)$ due to Lemma 2.35 and Definition of the
$\bf K$ integral along rectifiable paths, since for $|\xi |> \sup
|\chi |$ the series $R(\xi ;T_{\chi })$ converges uniformly in the
norm operator topology for the spectral set $\chi $ of the spectrum
$\sigma (T)$, where $T_{\chi }=T|_{X_{\chi }}$, $X_{\chi }:={\hat
E}_e(\chi ;T)X$. The variation of $f$ and $g$ gives that the index
of $\lambda $ is not less than $p$. Vise versa let be satisfied
Conditions $(ii)$ for some $n$. The resolvent $R(z;T)x$ is regular
on ${\bf K}\setminus B({\bf K},a,r)$ and
$$x=(2\pi )^{-1}\{ \int_{\eta }R(\zeta ;T)x d\zeta \}
M^{-1}(2\pi )^{-1}\{ \int_{\eta _2}R(\zeta ;T)x d\zeta \}
M_2^{-1}\omega (T)x,$$ where $\omega (T)$ is the function equal to
$1$ on a neighborhood of $\lambda $ and equal to zero on ${\bf
K}\setminus U$, $\eta $ is the corresponding closed rectifiable path
encompassing $\sigma (T)$ and characterized by $M\in \bf K$,
$|M|=1$, $M+{\tilde M}=0$. Then due to $2.40.(ii)$ $\phi
_{p(z)}(z,T)x=0$ for each $z\in \lambda $.
\par {\bf 42. Note.} An isolated point $\lambda $
of a spectrum $\sigma (T)$ for a normal quasilinear operator $T\in
L_q(X)$ in a Hilbert space $X$ over $\bf K$ can be not having
eigenvectors because of the noncommutativity of a projection valued
measure $\hat E$ apart from the case of linear operators on a
Hilbert space over $\bf C$.
\par {\bf 43. Examples. (1).} Consider a Hilbert space
$X := L^2(B({\bf K},0,1),\nu ,{\bf K})$ of all classes of $\nu $
measurable functions on $B({\bf K},0,1)$ with the finite norm and
the scalar product $$(f,g) := \int_{B({\bf K},0,1)} {\tilde f}(z)
g(z) \nu (dz)$$ and the norm $\| f \| := (f,f)^{1/2}$, where $\nu $
is the Lebesgue measure on $\bf R^{2^{m+1}}$, which induces the
measure on $\bf K$ such that it is bounded on the ball of the unit
radius in $\bf K$ with hte centre at zero, $f$ and $g$ are
measurable functions with finite norms from $B({\bf K},0,1)$ into
$\bf K$, $m=3$ for $\bf H$, $m=7$ for $\bf O$.
\par Take on this Hilbert space $X$ the  operator
given by the following formula $(Af)(x) = xf(x)$. Then $A$ does not
have eigenvalues. If $Af=bf$ for some $b\in \bf K$, then $f=0$ on
$B({\bf K},0,1)\setminus \{ b \} $, consequently, $f=0$ $\nu
$-almost everywhere, that is $f=0\in X$. But, $sp_{L_q(X)}= B({\bf
K},0,1)$. For the proof of the latter equality consider the
characteristic functions $y_n := \chi _{B({\bf K},b,1/(2n))}$ of the
balls $B({\bf K},b,1/(2n))$. Take $x_n := y_n/ [\nu (B({\bf
K},b,1/(2n)))]^{1/2}  $. In view of $\nu (B({\bf K},z,r)) = (2\pi
)^{2^m}r^{2^{m+1}}/[(2^{m+1})!!]$ (see, for example, 3 \S XI.4.2 in
\cite{zorich}), where $m=3$ for $\bf H$, $m=7$ for $\bf O$, then $
\| x_n \| =1$ and for $B({\bf K},b,1/(2n))\subset B({\bf K},0,1)$
the following equality is accomplished  $$\| (A-bI)x_n \|^2=
\int_{B({\bf K},b,1/(2n))} |x-b|^2\nu (dx)/\nu (B({\bf
K},b,1/(2n)))=C n^{-2},$$ where $C=const>0$. If $B = (A-bI)^{-1}$,
then $1 = \| x_n \| = \| B(A-bI) x_n \| \le \| B \| \| (A-bI) x_n \|
= \| B \| C^{1/2}/n$. Thus, $ \| B \| \ge n C^{-1/2}$ for each
$n=1,2,3,...$, consequently, $B$ is unbounded, that is, $(A-bI)$
does not have a two sided inverse operator in $L_q(X)$ and
inevitably $B({\bf K},0,1)\subset sp_{L_q(X)}(A)$. If $b\notin
B({\bf K},0,1)$, then the function $f=(x-b)^{-1}$ for $|x|\le 1$ is
continuous on $B({\bf K},0,1)$. The multiplication on $f$ in $X$ is
the bounded operator, which is two sided inverse to $(A-bI)$:
$f(A-bI)=(A-bI)f=I$. Therefore, $b\notin sp_{L_q(X)}(A)$,
consequently, $sp_{L_q(X)}= B({\bf K},0,1)$.
\par {\bf (2).} Let $\{ e_n: n=1,2,3,... \} $ be an orthonoramal
basis in a separable Hilbert space $X$. Take the operator $A$ such
that $Ae_n=b_ne_n$, where $b_n \in \bf K$ for each $n$. If the
sequence $\{ b_n: n=1,2,3,.. \} $ is bounded, then $ \| A \| =
\sup_n |b_n| $ and the operator $A$ is normal. At the same time
$sp_{L_q(X)}(A)=cl \{ b_n: n=1,2,3,... \} $, since the spectrum of
the operator is closed in accordance with the general theorem. The
operator $A$ is selfadjoint if and only if $b_n\in \bf R$ for each
$n$. The unitarity of the operator $A$ is equivalent to $|b_n|=1$
for each $n$.
\par {\bf (3).} Let $(X,{\hat \mu })$ be a $\sigma $-finite
space with a (noncommutative) $\bf K$-valued measure. Consider a
$\hat \mu $-measurable function $f$ with the finite norm $\| f \|
_{L^{\infty }(X,{\hat \mu },{\bf K})}<\infty $, then define on the
space $L^2(X,{\hat \mu },{\bf K})=:X$ the multiplication operator
$M_f$ on the function $f$, that is, $M_f(g):=fg$. This operator is
bounded. Let a point $b$ belongs to the $\hat \mu $-essential domain
of the function $f$: $sp(f) := \{ z\in {\bf K}: \| {\hat \mu } \|
(f^{-1}(U))>0$ $\mbox{ for each open subset}$ $U\subset \bf K$ with
$z\in U \} $. Let $y_n:= \chi _{f^{-1}(U_n)}$, where $\| {\hat \mu }
\| (f^{-1}(U_n))=a_n>0$, where $U_n := \{ z\in {\bf K}: |z-b|<1/n \}
$, put $x_n := (a_n)^{-1/2}y_n$, consequently, $\| x_n \| =1$, also
$$ \| (M_f-bI)x_n \| ^2=\int_{f^{-1}(U_n)} |f(z)-b|^2 y_n(z){\hat \mu
}(dz)/a_n \le n^{-2}.$$ Therefore, there exists $ \lim_{n\to \infty
} \| (M_f-bI)x_n \| =0$ and $b\in sp_{L_q(X)}(M_f)$.

\section{Algebras of operators}
\par {\bf 1. Definition.} Let $V$ be a compact Hausdorff
space, also $C(V,{\bf K})$ be a quasicommutative algebra of all
continuous bounded functions $f: V\to \bf K$ with the supremum norm
$\| f \| :=\sup_{z\in V} |f(z)|$ and with the pointwise
multiplication, where ${\bf K} = \bf H$ or ${\bf K} = \bf O$. A
quasilinear functional $p\in L_q(C(V,{\bf K}),{\bf K})$ satisfying
the condition $p(fg)=p(f)p(g)$ for each $f, g \in C(V,{\bf K})$ is
called multiplicative.
\par {\bf 2. Theorem.} {\it If $J$ is a closed ideal in $C(V,{\bf K})$, then
there exists a closed subset $S$ in $V$ such that $J$ is the ideal
$C(V,{\bf K})_S$ of all functions from $C(V,{\bf K})$ equal to zero
on $S$ up to an automorphism of the algebra $\bf K$, where ${\bf
K}=\bf H$ or ${\bf K}=\bf O$. Vise versa, if $S$ is a closed subset
in $V$, then $C(V,{\bf K})_S$ is a closed ideal in $C(V,{\bf K})$.
An ideal $J$ in $C(V,{\bf K})$ is maximal if and only if the
corresponding closed subset $S=S(J)$ is a singleton.}
\par {\bf Proof.} The Banach algebra $C(V,{\bf K})$
is noncommutative, moreover, it is nonassociative for ${\bf K}=\bf
O$, where ${\bf K}=\bf H$ or ${\bf K}=\bf O$. If $f\in C(V,{\bf
K})$, then $f^{-1}(0)$ is a closed subset in $V$. Then the subset $S
:= \bigcap_{f\in J}f^{-1}(0)$ is closed in $V$. If $J=C(V,{\bf K})$,
then $S=\emptyset $. If $J$ is an ideal of the algebra $C(V,{\bf
K})$ and $\theta : {\bf K}\to \bf K$ is an automorphism of the
algebra $\bf K$ such that $\theta ({\bf R})=\bf R$, then $J$ and
$J\circ \theta := \{ g=f\circ \theta : f \in J \} $ are isomorphic
algebras. At the same time the group  of generators $ \{
i_0,i_1,...,i_m \} $ generates an isomorphic group $ \{ \theta
(i_0), \theta (i_1), ..., \theta (i_m ) \} $, where $m=3$ for $\bf
H$ and $m=7$ for $\bf O$ (in the latter case the group of generators
is nonassociative), $i_0:=1$, $i_1:=i$, $i_2:=j$, $i_3:=k$,
$i_4:=l$, $i_5=i_1l$, $i_6=i_2l$, $i_7=i_3l$, where $l$ is the
generator of the doubling procedure while construction of the
algebra $\bf O$ of octonions from the skew field of quaternions $\bf
H$ (see about generators of the algebras $\bf H$ and $\bf O$ in
\cite{baez,kurosh}). On the other hand, if two ideals $J_1$ and
$J_2$ are isomorphic, then $J_1\cap C(V,{\bf R})$ and $J_2\cap
C(V,{\bf R})$ are isomorphic algebras over $\bf R$, also the
extension of this isomorphism on $J_1$ and $J_2$ induces an
automorphism of the algebra $\bf K$, since for a compact set $V$ all
constants from $\bf K$ belong to $C(V,{\bf K})$, also $J_1$ and
$J_2$ are algebras over $\bf K$, that is, they are linear spaces
over $\bf K$, also in them there is defined the multiplication of
elements, which is distributive (associative in the particular case
of the skew field ${\bf K}=\bf H$ of quaternions).
\par Consider a function $f\in C(V,{\bf K})$ which is bounded on $V$
and equal to zero on $S$, $f|_S=0$. For each $\epsilon
>0$ the set $S_{f,\epsilon }:= \{ z\in V: |f(z)|\ge \epsilon \} $
is compact and $S_{f,\epsilon }\cap V = \emptyset $. If $z\in
S_{f,\epsilon }$, then there exists a function $f_z\in J$ such that
$f_z(z)\ne 0$, consequently, the function ${\tilde f}_zf_z$ is
positive on some neighborhood $U_z$ of the point $z$. Due to
compactness of the set $S_{f,\epsilon }$ there exists a finite set
of points $z_1,...,z_n$ such that $S_{f,\epsilon }\subset
\bigcup_{m=1}^nU_{z_m}$, consequently, $h_{\epsilon }:=\sum_{m=1}^n
{\tilde f}_{z_m}f_{z_m}>0$ on $S_{f,\epsilon }$ and attains the
minimum $u$ on this compact set at some point $t$ with $h_{\epsilon
}(t)>0$. Then the equation $w_{\epsilon }(z):= \max ( h_{\epsilon
}(z),u)$ defines a continuous function on $V$ satisfying the
inequalities $w_{\epsilon }(z)>0$ and $w_{\epsilon }(z)\ge
h_{\epsilon }(z)\ge 0$ for each $z\in V$. In view of $1/w_{\epsilon
}\in C(V,{\bf K})$ and $h_{\epsilon }\in J$ let us take the function
$g_{\epsilon }:=h_{\epsilon }w_{\epsilon }^{-1}$. It satisfies the
conditions $g_{\epsilon }\in  J$ and $0\le g_{\epsilon }(z)\le 1$
for each point $z\in V$, also $g_{\epsilon }(z)=1$ on $S_{f,\epsilon
}$. Then $fg_{\epsilon }\in J$ and $ \| f-fg_{\epsilon } \|  \le
\epsilon $, since $ \| 1- g_{\epsilon } \| \le 1$ and $1=g_{\epsilon
}(z)$ for each $z\in S_{f,\epsilon }$. In view of the fact that the
ideal $J$ is closed, it follows that $f\in J$.
\par Each compact Hausdorff space is normal.
If $S$ is a closed subset in $V$, also if $z\in V\setminus S$, then
due to the Tietze-Urysohn theorem (see Theorems 2.1.8, 3.1.7 and
3.1.9 in \cite{eng}) about extension of continuous functions there
exists a continuous function on $V$ such that $f|_S=0$ and $f(z)=1$.
Thus $f\in J$ and $f(z)\ne 0$, consequently, if $J=C(V,{\bf K})_T$,
then $S(J)=T$.
\par The last statement of this theorem follows now from the inclusion $C(V,{\bf
K})_T\supset C(V,{\bf K})_S$, if a closed set $T$ is contained in
the closed subset $S$. Certainly, the ideals here are defined up to
the  automorphism of the algebra $\bf K$, that gives isomorphisms of
ideals.
\par {\bf 3. Corollary.} {\it For each nonzero continuous
quasilinear multiplicative functional $p$ on a Banach algebra
$C(V,{\bf K})$ there exists a point $z_0\in V$ such that
$p(f)=f(z_0)$ for each function $f\in C(V,{\bf K})$ up to an
automorphism of the algebra $\bf K$, where ${\bf K}=\bf H$ or ${\bf
K}=\bf O$.}
\par {\bf Proof.} The Banach algebra $C (V,{\bf K})$
has a decomposition into the direct sum: $C (V,{\bf K}) = C (V,{\bf
R})\oplus i_1 C (V,{\bf R}) \oplus ... \oplus i_m C (V,{\bf R})$,
where $\{ i_0, i_1,..., i_m \} $ are standard generators of the
algebra $\bf K$, $m=3$ for the skew field of quaternions $\bf H$,
$m=7$ for the algebra of octonions $\bf O$, $i_0:=1$. $i_1:=i$,
$i_2:=j$, $i_3:=k$, $i_4:=l$, $i_5=i_1l$, $i_6=i_2l$, $i_7=i_3l$,
where $l$ is the generator of the doubling procedure while
construction of $\bf O$ from $\bf H$. Then $p(f) = p(f_0)p(i_0) +
p(f_1)p(i_1) + ... + p(f_m)p(i_m)$, where $f=f_0+f_1i_1+...+f_mi_m$,
$f_0, f_1,..., f_m \in C(V,{\bf R})$, $f\in C(V,{\bf K})$. From the
identities for generators and the multiplicativity of the functional
$p$ it follows, that $p(1)=p(1^2)=p(1)^2=1$, since all constants
belong to the algebra $C(V,{\bf K})$, $p(i_q^2)=-1$ for each
$q=1,2,...,m$, $p(ij)=-p(ji)=p(k)=p(i)p(j)=-p(j)p(i)$,
$p(jk)=-p(kj)=p(i)=p(j)p(k)=-p(k)p(j)$,
$p(ki)=-p(ik)=p(j)=p(k)p(i)=-p(i)p(k)$, $p(il)=-p(li)$,
$p((il)i)=-p(i(il))=-p(l)$ and so on, that is, $p$ gives the
automorphism of the algebra $\bf K$. The kernel $M:=p^{-1}(0)$ of
the functional $p$ is the proper maximal ideal in $C(V,{\bf K})$.
Let $f\in C(V,{\bf K})$, then $f-p(f)1\in M$. Due to Theorem 3.2 the
ideal $M$ is characterized uniquely up to an automorphism of the
algebra $\bf K$ by some point $z_0\in V$. Thus, $f(z_0)=p(f)$ for
each $f\in C(V,{\bf K})$ up to an isomorphism of the algebra and of
the ideals relative to automorphisms of the algebra $\bf K$.
\par {\bf 4. Theorem.} {\it A mapping $\phi $ of the algebra $C(V,{\bf K})$
onto the algebra $C(W,{\bf K})$ for compact Hausdorff spaces $V$
and $W$ is an isomorphism if and only if there exists a \\
homeomorphism $g: V\to W$ and an automorphism $\theta $ of the
algebra $\bf K$, such that $\phi (f)=\theta \circ f\circ g$ for each
$f\in C(V,{\bf K})$.}
\par {\bf Proof.} If there exists a homeomorphism
$g: V\to W$ and an automorphism $\theta $ of the algebra $\bf K$
such that $\phi (f)=\theta \circ f\circ g$ for each $f\in C(V,{\bf
K})$, then from the quasilinearity and multiplicativity of $\theta $
it follows the quasilinearity and multiplicativity of the mapping
$\phi $, since $(f_1f_2)\circ g(z) = f_1(g(z)) f_2(g(z))$ and
$(af_1+bf_2)\circ g(z)=af_1(g(z))+bf_2(g(z))$ for each $f_1, f_2 \in
C(V,{\bf K})$, $a, b\in \bf K$ and each $z\in V$.
\par Let now $\phi $ be an algebraic isomorphism of the algebra
$C(V,{\bf K})$ onto the algebra $C(W,{\bf K})$, then for a maximal
ideal $M$ in $C(W,{\bf K})$ its inverse image $\phi ^{-1}(M)$ is the
maximal ideal in $C(V,{\bf K})$. In view of the fact that $V$ and
$W$ are compact, it follows that constants from $\bf K$ belong to
these algebras $C(V,{\bf K})$ and $C(W,{\bf K})$, consequently, the
restriction of $\phi $ on the subalgebra isomorphic with $\bf K$ in
$C(V,{\bf K})$ gives some automorphism $\theta $ of the algebra $\bf
K$ such that $\theta ^{-1}\circ \phi =: \phi _1$ is the
 identity mapping on $\bf K$. If $z\in W$ is a point corresponding to $M$,
 then denote by $q(z)$ a point in $V$
corresponding to $\phi _1^{-1}(M)$. In view of the fact that each
maximal ideal $M_0$ in $C(V,{\bf K})$ has the form $\phi
_1^{-1}(\phi _1(M_0))$ with the maximal ideal $\phi _1(M_0)$ in
$C(W,{\bf K})$, then $q: W\to V$ is the uniquely defined mapping. In
view of the fact that $f-f(q(z))1$ is equal to zero at the point
$q(z)$ for each function $f\in C(V,{\bf K})$, then $\phi
_1(f-f(q(z))1)$ is equal to zero at the point $q(z)$. Since a
compact Hausdorff space is completely regular, then the set of the
form $\phi ^{-1}(U)$, where $U$ is an open subset in $\bf K$ and
$f\in C(V,{\bf K})$, form a (sub)base of the topology in $V$. From
$q^{-1}(f^{-1}(U))=\phi _1(f)^{-1}(U)$ and the continuity of the
function $\phi _1(f): W\to \bf K$ for the continuous function $f:
V\to \bf K$ it follows, that the mapping $q$ is continuous, since
the inverse image $q^{-1}$ of each open subset of the subbase in $V$
is an open set of the subbase of the topology in $W$. Symmetrically
$q^{-1}$ is also continuous, therefore $q: W\to V$ - is the
homeomorphism.
\par {\bf 5. Definitions.} Let be given a space $C(V,{\bf R})$
of all continuous functions from a Tychonoff topological space $V$
into $\bf R$, also let be given its subset $\cal P$ satisfying the
following properties:
\par $(i)$ if $f$ and $-f\in \cal P$, then $f=0$;   \par $(ii)$ if $a>0$
is a positive constant and $f\in \cal P$, then $af\in \cal P$;
\par $(iii)$ if $f, g \in \cal P$, then $f+g\in \cal P$; so such $\cal P$
is called the cone in $C(V,{\bf R})$.
\par If $X_0$ is a linear over $\bf R$ subspace in $C(V,{\bf
R})$ and $\cal P$ is a cone in it, then in $X_0$ there is given a
partial ordering $f\le g$ if and only if $g-f\in \cal P$. A constant
function $I$ equal to $1$ on $V$ is called the ordering unit in
$X_0$, if for each $f\in X_0$ there exists a constant $a>0$ such
that $-aI\le f\le aI$.
\par Let $X_0$ be a partially ordered subspace in $C(V,{\bf
R})$ with the ordering unit $I$, also $p$ be a quasilinear
functional on a $\bf K$-linear subspace $X$ in $C(V,{\bf K})$ such
that $X_0=X\cap C(V,{\bf R})$. If a functional $p$ satisfies the
following conditions:
\par  $(iv)$ $p(f)\ge 0$ for each $f\in X_0$ with $f\ge 0$;
\par $(v)$ a bounded functional $p$ on $\bf K$ is an automorphism of the algebra
$\bf K$, if $p\ne 0$, then $p$ is called positive; if in addition
\par $(vi)$ $p(I)=1$, then $p$ is called a state on $X$.
\par If $p$ is an extreme point of a (convex) family ${\cal L}(X)$
of states for $X$, then $p$ is called a pure state on a $\bf
K$-linear (sub)space $X$.
\par {\bf 6. Theorem.} {\it A nonzero quasilinear functional $p$ on
$C(V,{\bf K})$, where ${\bf K}=\bf H$ or ${\bf K}=\bf O$, for a
compact Hausdorff topological space $V$ is a pure state on $C(V,{\bf
K})$ if and only if it is multiplicative.}
\par {\bf Proof.} Consider at first a bounded functional $p_0$
of a functional $p$ on $X_0$. Let $p_0$ be satisfying the property:
for each positive functional $q$ on $X_0$ such that $q\le p$ there
exists a constant $a\ge 0$, for which $q = ap_0$. Suppose that
$p_0=aq_1+(1-a)q_2$, where $q_1$ and $q_2$ are states on $X_0$,
$0<a<1$. Then $0\le aq_1\le p_0$, therefore, $aq_1=bp_0$ for some
constant $b>0$. From $p_0(I)=q_1(I)=1$ it follows, that $a=b$ and
$q_1=p_0$. Analogously $q_2=p_0$, consequently, $p_0$ is the pure
state. \par If $p_0$ is a pure state and $0\le q\le p_0$ on $X_0$,
then $0\le q(I)\le p_0(I)=1$. In the case $q(I)=0$ for each $f\in
X_0$ are accomplished the inequalities $0=q(-aI)\le q(f)\le q(aI)=0$
for some constant $a$, consequently, $q(f)=0$ and $q=0$. If
$q(I)=1=p_0(I)$, then the proof  above shows, that the positive
functional $p_0-q$ is equal to zero, that is, $q=p_0$. If
$0<q(I)<1$, then $q=(1-b)q_1+bq_2$, where $b=q(I)$, also $q_1$ and
$q_2$ are states defined by the equalities $q_1=(1-b)^{-1}(p_0-q)$,
$q_2=b^{-1}q$. In view of the fact that $p_0$ is the pure state, it
follows that $q_2=p_0$ and $q=bp_0$. Thus, $p_0$ is the pure state
on $X_0$ (without Condition $(v)$) with the unit of the ordering $I$
if and only if for each positive functional $q$ on $X_0$ such that
$q\le p_0$, there exists a constant $a\ge 0$, for which $q=ap_0$.
\par Due to Corollary 3.3 for a multiplicative functional $p$
there exists a point $z_0\in V$ and an automorphism $\theta $ of the
algebra $\bf K$ such that $p_1(f)=f(z_0)$ for each function $f\in
C(V,{\bf K})$, where $p_1=\theta \circ p$. Each continuous function
equal to zero at the point $z_0$ is the linear combination of two
nonnegative functions equal to zero at the point $z_0$. Therefore,
from $0<q\le p$ for some quasilinear functional $q$ on $C(V,{\bf
K})$, it follows, that $q(f)=0$ for $f(z_0)=0$. Thus, the
functionals $q$ and $p$ have the same subspaces of functions in
$C(V,{\bf K})$, on which they are equal to zero. Therefore, there
exists an automorphism $\theta _{p,q}$ of the algebra $\bf K$, for
which $q=a\theta _{p,q}\circ p$ for some positive constant $a$. Due
to the first part of the proof $p$ is the pure state on $C(V,{\bf
K})$, since $p(I)=1$.
\par Suppose now, that $p$ is the pure state on
$C(V,{\bf K})$. If $0\le f\le 1$ and $q(g)=p(fg)$, then $q$ is the
quasilinear functional on $C(V,{\bf K})$ such that $q\le p$. Thus,
$q=ap$ for some constant $a\ge 0$. If $p(h)=0$, then
$p(fh)=q(h)=ap(h)=0$. In view of the fact that each $g$ from
$C(V,{\bf K})$ is a $\bf K$-linear combination of functions with
nonnegative values between $0$ and $1$, then $p(gh)=0$ for each
continuous function $g\in C(V,{\bf K})$. Thus, the kernel of the
functional $p$ is the maximal ideal in $C(V,{\bf K})$. From $p(I)=1$
it follows the multiplicativity of the functional $p$.
\par {\bf 7. Definition.} Let $V$ be a Tychonoff topological space,
also ${\bf K}=\bf H$ or ${\bf K}=\bf O$. A quasilinear functional
$p$ on $C(V,{\bf K})$ is called Hermitian, if $p({\tilde
f})=(p(f))^{\tilde .}$ for each continuous function $f\in C(V,{\bf
K})$, also $p(fi_q)=p(f)p(i_q)$ for each continuous real valued
function $f\in C(V,{\bf R})$ and the standard generator $i_q$ of the
algebra $\bf K$, $q=0,1,...,m$, $m=3$ for $\bf H$, $m=7$ for $\bf
O$, where a number $b\in \bf O$ is considered also as the constant
function on $V$.
\par {\bf 8. Proposition.} {\it A family $P$ of Hermitian quasilinear functionals
$p$ on $C(V,{\bf K})$, having the property $p_1(i_q)=p_2(i_q)$ for
each $p_1, p_2\in P$ and each $q=0,1,...,m$, in the normed dual
space generates a lattice.}
\par {\bf Proof.} For Hermitian functionals $p_1$ and $p_2$
and a positive continuous function $f\in C(V,{\bf K})$ we define
$$(i)\quad (p_1\vee p_2)(f):=\sup \{ p_1(f_1)+p_2(f_2): f_1, f_2
\in C(V,{\bf K}), 0\le f_1, 0\le f_2, f=f_1+f_2 \} .$$
 Then $|p_1(f_1)+p_2(f_2)| \le \| p_1 \| \| f_1
\| + \| p_2 \| \| f_2 \| \le (\| p_1 \| + \| p_2 \| ) \| f \| $,
consequently, $(p_1\vee p_2)(f)$ is finite. In view of $f+0=0+f=f$,
we have $p_1\le p_1\vee p_2$ and $p_2\le p_1\vee p_2$. If $q$ is a
Hermitian quasilinear functional such that $p_1\le q$ and $p_2\le
q$, then $p_1(f_1)+p_2(f_2)\le q(f_1)+q(f_2)=q(f)$, that is,
$p_1\vee p_2\le q$. \par We show that $p_1\vee p_2$ is the $\bf
R$-linear on the cone of positive elements. Consider now nonnegative
continuous functions on $V$. Choose $g_1$ and $g_2$ such that
$g=g_1+g_2$, also $h_1$ and $h_2$ such that $h=h_1+h_2$, then
$g_1+h_1+g_2+h_2=g+h$ and \par
$p_1(g_1+h_1)+p_2(g_2+h_2)=p_1(g_1)+p_1(h_1)+p_2(g_2)+p_2(h_2)\le
(p_1\vee p_2)(g+h)$. Thus,
$$(ii)\quad (p_1\vee p_2)(g)+(p_1\vee p_2)(h)\le (p_1\vee p_2)(g+h).$$
 Let $g+h=f_1+f_2$ and
$p_1(f_1)+p_2(f_2)$ approximates $(p_1\vee p_2)(g+h)$. Use the
decompositions $f_1=f_{11}+f_{12}$, $f_2=f_{21}+f_{22}$ with
$f_{11}+f_{21}=g$, $f_{12}+f_{22}=h$. This is possible, since
$f_1\le g+h$, $f_{11}\le g$ and $f_{12}\le h$,
$f_2=g+h-f_1=g-f_{11}+h-f_{12}$, putting $f_{12}=g-f_{11}$ and
$f_{22}=h-f_{12}$. Therefore,
\par $p_1(f_1)+p_2(f_2)=p_1(f_{11})+p_1(f_{12})+p_2(f_{21})+p_2(f_{22})$
$\le (p_1\vee p_2)(g)+(p_1\vee p_2)(h)$. Thus, $$(iii)\quad (p_1\vee
p_2)(g+h)\le (p_1\vee p_2)(g)+(p_1\vee p_2)(h).$$
 Then from the
inequality $(ii,iii)$ it follows the additivity of $p_1\vee p_2$ on
the cone of positive functions. If $0<a$, $f_1+f_2=f$ and
$p_1(f_1)+p_2(f_2)$ approximates $(p_1\vee p_2)(f)$, then
$a(p_1(f_1)+p_2(f_2))\le (p_1\vee p_2)(af)$, that is, $a(p_1\vee
p_2)(f)\le (p_1\vee p_2)(af)$, consequently, $a^{-1}(p_1\vee
p_2)(af)\le (p_1\vee p_2)(a^{-1}af)=(p_1\vee p_2)(f)$. Thus,
$a(p_1\vee p_2)=(p_1\vee p_2)(af)$.
\par Let now $f\in C(V,{\bf R})$. Then there exist
nonnegative continuous functions $f_1$ and $f_2$ such that
$f=f_1-f_2$. Then define $(p_1\vee p_2)(f):=(p_1\vee p_2)(f_1)-
(p_1\vee p_2)(f_2)$.  If $f=g_1-g_2$ with nonnegative functions
$g_1$ and $g_2$, then $f_1+g_2=g_1+f_2$, therefore, $(p_1\vee
p_2)(f_1+g_2)= (p_1\vee p_2)(f_1)+(p_1\vee p_2)(g_2)=(p_1\vee
p_2)(g_1+f_2)=(p_1\vee p_2)(g_1)+(p_1\vee p_2)(f_2)$, consequently,
$(p_1\vee p_2)(f_1)-(p_1\vee p_2)(f_2)=(p_1\vee p_2)(g_1)-(p_1\vee
p_2)(g_2)$. Thus, $(p_1\vee p_2)$ has the extension onto $C(V,{\bf
R})$ uniquely. The additivity and positive homogeneousity $(p_1\vee
p_2)$ on the cone of nonnegative functions give the additivity and
real homogeneousity of the functional $(p_1\vee p_2)$ on $C(V,{\bf
R})$. \par In view of $p_1(i_q)=p_2(i_q)$ for each $q=0,1,...,m$ and
each $p_1, p_2\in P$, for arbitrary continuous function $f\in
C(V,{\bf K})$ and $p:=p_1\vee p_2$ put
$p(f)=p(f_0)p(1)+p(f_1)p(i_1)+...p(f_m)p(i_m)$, where
$p(i_q):=p_1(i_q)$ for each $q=0,1,...,m$, the function $f$ is
written in the form $f=f_0+f_1i_1+...+f_mi_m$ with real continuous
functions $f_0,f_1,...,f_m\in C(V,{\bf R})$ and standard generators
$\{ i_0,i_1,...,i_m \} $, $i_0=1$, $m=3$ for $\bf H$, $m=7$ for $\bf
O$. From $p_s({\tilde f})=(p_s(f))^{\tilde .}$ for $s=1$ and $s=2$
it follows, that $p({\tilde
f})=p(f_0)-p(f_1)p(i_1)-...-p(f_m)p(i_m)=(p(f))^{\tilde .}$.
\par We define also $p_1\wedge p_2 := -((-p_1)\vee (-p_2))$.
Therefore, $p_1\vee p_2$ and $p_1\wedge p_2$ have finite norms.
\par {\bf 9. Proposition.} {\it If $p$ is a bounded quasilinear
functional on $C(V,{\bf K})$, $p^+ := p\vee 0$, $p^- := -(p\wedge
0)$, then $p=p^+ - p^-$, $p^+\wedge p^- =0$ and $ \| p \| = \| p^+
\| + \| p^- \| =p^+ (1) + p^- (1)$.}
\par {\bf Proof.} For a positive continuous function
$f\in C(V,{\bf R})$ it is accomplished the equality:
\par $p^+(f)=\sup \{ p(f_1): 0\le f_1\le f, f_1\in C(V,{\bf R}) \}
$,
\par $-p^- (f)= (p\wedge 0)(f)=-((-p)\vee 0)(f)=-\sup \{ - p(g):
0\le g\le f, g\in C(V,{\bf R}) \} $
\par $= \inf \{ p(g): 0\le g\le f, g\in C(V,{\bf R}) =
\inf \{ p(f-f_1): 0\le f_1\le f, f_1\in C(V,{\bf R}) \} $
\par $=p(f)-\sup \{ p(f_1): 0\le f_1\le f , f_1\in C(V,{\bf R}) \}
=p(f)-p^+ (f)$. \\
Thus, $p=p^+ - p^-$, since $0(f)=0$ for each $f\in C(V,{\bf K})$,
also putting $p^+(i_q)=p^-(i_q)=p(i_q)$ for each $q=0,1,...,m$. For
a nonnegative function $f\in C(V,{\bf R})$
\par $(p^+ \wedge p^-)(f)=\inf \{ p^+(f_1)+p^-(f-f_1): 0\le f_1\le
f, f_1\in C(V,{\bf R}) \} $  $ = \inf \{ p(f_1) + p^-(f): 0\le
f_1\le f, f_1\in C(V,{\bf R}) \} $  $= p^- (f) + \inf \{ p(f_1):
0\le f_1\le f, f_1\in C(V,{\bf R}) \} =0$.
\par If $v$ is a nonnegative quasilinear functional and $v(1)=0$,
then $v=0$, since $I$ is the unity ordering. Therefore, $ \| v \| =
v(1) =0$. If $v(1)\ne 0$, then $v(1)^{-1}v$ is the state on
$C(V,{\bf K})$ and $ \| v(1)^{-1}v \| =1$, consequently, $ \| v \| =
v(1)$. Thus,
\par $ \| p \| \le \| p^+ \| + \| p^- \| = p^+(1) + p^- (1) $
$ = \sup \{ p(f): 0\le f\le 1, f\in C(V,{\bf R}) \} $ $- \inf \{
p(g): 0\le g\le 1, g\in C(V,{\bf R}) \} $. \par For suitable $f,
g\in C(V,{\bf R})$ satisfying inequalities $0\le f\le 1$ and $0\le
g\le 1$, $p(f)-p(g)$ approximates $p^+(1)+p^-(1)$ and $p(f)-p(g)\le
\| p \| \| f-g \| \le \| p \| .$ Therefore, $p^+(1)+p^-(1) \le \| p
\| $ and this togehter with the inequalities gives: $ \| p \| = \|
p^+ \| + \| p^- \| =p^+(1) + p^- (1)$.
\par {\bf 10. Theorem.} {\it If $\cal A$ is a closed relative to the norm
subalgebra in $C(V,{\bf K})$ for a Hausdorff compact topological
space $V$, ${\bf K}=\bf H$ or ${\bf K}=\bf O$, $1\in \cal A$,
${\tilde f}\in \cal A$ for each function $f\in \cal A$, also for
each pair of different points $z_0, z_1 \in V$ there exists a
function $f\in \cal A$ such that $f(z_0)\ne f(z_1)$, then ${\cal
A}=C(V,{\bf K})$.}
\par {\bf Proof.} An algebra ${\cal A}$ has the decomposition
${\cal A}={\cal A}_0\oplus {\cal A}_1i_1\oplus ... \oplus {\cal
A}_mi_m$, where ${\cal A}_0={\cal A}\cap C(V,{\bf R})$,
$i_0,i_1,...,i_m$ are generators of the algebra $\bf K$, ${\cal
A}_q\subset C(V,{\bf R})$ for each $q$. At the same time the
algebras ${\cal A}_q$ over $\bf R$ are contained in $C(V,{\bf R})$,
also ${\cal A}_q$ is isomorphic with ${\cal A}_0$ for each
$q=1,...,m$, since $\cal A$ is the algebra over $\bf K$, in
particular, it is invariant relative to the multiplication on each
$0\ne b\in \bf K$, where $i_0=1$. This theorem is accomplished for
the algebra over $\bf R$, that is, ${\cal A}_0=C(V,{\bf R})$, that
composes the context of the Stone-Weierstrass theorem, consequently,
${\cal A}=C(V,{\bf K})$.
\par {\bf 11. Theorem.} {\it An element $A$ of a Banach algebra $\cal C$ over $\bf K$
has a spectral radius $r(A)$ given by the formula: $r(A)=\lim_{n\to
\infty } \| A^n \| ^{1/n}$.}
\par {\bf Proof.} Due to Lemma 2.14 there is accomplished the inequality
$ r(A)\le \| A \| .$ If $a\in sp (A)$, then due to Theorem 2.22
$a^n\in sp (A^n)$, since the skew field of quaternions $\bf H$ is
associative, also the algebra of octonions is alternative,
therefore, $|a^n|\le \| A^n \| $. Then $$ |a| \le {\underline {\lim
}} \| A^n \| ^{1/n}\mbox{ and }r(A) \le {\underline {\lim }} \| A^n
\| ^{1/n}.$$
 The function $R(z;A)$ by the variable $z$ is $\bf
K$-holomorphic on $\rho (A)$. The power series corresponding to
$R(z;A)$ diverges for $|z|> ( {\overline {\lim }} \| A^n \|
^{1/n})^{-1} $. That is, if $0\le a'< {\overline {\lim }} \| A^n \|
^{1/n}$, then there exists $a\in \bf K$, such that $a'< |a|$, for
which $I-a^{-1}A$, hence $A- a I$ also has not an inverse element in
$\cal C$. Therefore, $a \in sp (A)$ and $a'<r(A)$. In view of
arbitrariness of a nonnegative number $a'$ less than ${\overline
{\lim }} \| A^n \| ^{1/n}$ we have $${\overline {\lim }} \| A^n \|
^{1/n} \le r(A)\le {\underline {\lim }} \| A^n \| ^{1/n}.$$
 On the other hand, it is always accomplished the inequality ${\underline
{\lim }} \| A^n \| ^{1/n}\le {\overline {\lim }} \| A^n \| ^{1/n} $,
consequently, there exists $r(A)=\lim_{n\to \infty } \| A^n \|
^{1/n}$.
\par {\bf 12. Corollary.} {\it If $A$ and $B$ are commuting elements in
a Banach algebra $\cal C$, then $r(AB)\le r(A)r(B)$ and $r(A+B)\le
r(A)+r(B)$.}
\par {\bf Proof.} Due to Theorem 3.11 $r_{\cal
B}(A)=r_{\cal C}(A)$, when $A$ belongs to the Banach subalgebra
$\cal B$ in $\cal C$. If $A$ and $B$ are commuting elements of the
algebra $\cal C$, then due to Theorem 2.22 for quasicommutative
subalgebra $\cal B$, generated by $A$ and $B$, from $$ \| (A B)^n \|
= \| A^nB^n \| \le \| A ^n \|  \| B^n \| \mbox{ and } \| (A+B)^n \|
\le \sum_{m=0}^n \| A^m \| \| B^{n-m} \| $$ it follows the statement
of this corollary.
\par {\bf 13. Proposition.} {\it Let $A$ be an element of a $C^*$-algebra $\cal C$.
\par $(i)$ If $A$ is normal, then $r(A) = \| A \| $.
\par $(ii)$ If $A$ is selfadjoint, then $sp (A)$ is a compact
subset in $\bf R$, also its spectrum contains at least one of the
numbers $ - \| A \| $ or $ \| A \| $.
\par $(iii)$ If $A$ is a unitary element in $\cal A$, then $ \| A \| =1$
and $sp (A)$ is a compact subset in the sphere $ \{ z\in {\bf K}:
|z|=1 \} $.}
\par {\bf Proof.} $(i)$ For a selfadjoint operator $B$
in $\cal C$ and a natural number $n>0$ it is accomplished the
equality $ \| (B)^{2n} \| = \| (B^n)^*B^n \| = \| B^n \| ^2$. By
induction by $q$ we get $ \| B ^q \| = \| B \| ^q$ for $q=2^m$, also
due to Theorem 3.11 $r(B)=\lim_{q\to \infty } \| B^q \| ^{1/q} = \|
B \| $. For a normal element $A$ and selfadjoint $B=AA^*$ from the
preceding discussion, Corollary 3.12 and $C^*$ property of the norm
it follows, that $$ \| A^2 \| = \| A^*A \| =r(A^*A)\le
r(A^*)r(A)=r(A)^2\le \| A \| ^2\mbox{, consequently, }r(A) = \| A \|
.$$
\par $(ii)$ for a selfadjoint element $A$ in the algebra $\cal C$
and its spectrum $sp (A)$ is compact in accordance with Lemma 2.14
and therefore contains a scalar with the absolute value $r(A)$, also
$r(A) = \| A \| $ due to Part $(i)$ of this proposition. Therefore,
it is sufficient to show that $sp (A)\subset \bf R$. Let $t\in sp
(A)$, $t = a + M b$, where $a, b\in \bf R$, $M + {\tilde M}=0$,
$|M|=1$, $M\in \bf K$.  Put $B_n=A-aI+MnbI$, then
$M(n+1)b=a+Mb-a+Mnb\in sp (B_n)$. At the same time $$(n^2+2n+1)b^2 =
|M(n+1)b|^2 \le [r(B_n)]^2\le \| B_n \| ^2=$$ $$\| B_n^*B_n \| = \|
(A-aI-MnbI)(A-aI+MnbI) \| $$  $$= \| (A-aI)^2+n^2b^2I +AM-MA\| \le
\| A -aI \| + 2 \| A \| + n^2b^2,$$ since $B_n^*B_n$ is selfadjoint,
$(MA)^* = A^*M^* = -AM$, $A^*M+M^*A = AM-MA$. Thus, $(2n+1)b^2\le \|
A -aI \| ^2+ 2 \| A \| $ for each $n=1,2,...$, therefore, $b=0$ and
$t=a\in \bf R$.
\par $(iii)$ For a unitary operator due to property of the norm in the
$C^*$-algebra we have: $ \| A \| ^2= \| A^*A \| = \| I \| =1$,
Therefore, $ \| A \| =1$. If $z\in sp (A)$, then $a^{-1}\in sp
(A^{-1})=sp (A^*)$, consequently, $|a|\le \| A \| =1$, $ |a|^{-1}\le
\| A^* \| =1$, Thus, $|a|=1$.
\par {\bf 14. Corollary.} {\it If $A$ is a normal element
in a $C^*$-algebra $\cal C$ and $A^k=0$ for some positive number
$k$, then $A=0$.}
\par {\bf Proof.} In view of $A^n=0$ for $n\ge k$, then from
Proposition 3.13 it follows, that $$ \| A \| =r(A) =\lim_{n\to
\infty } \| A ^n \| ^{1/n}=0.$$
\par {\bf 15. Theorem.} {\it If $A$ is a selfadjoint element of a $C^*$-algebra
$\cal C$ over ${\bf K}=\bf H$ or ${\bf K}=\bf O$, then there exists
a unique continuous mapping $f: C(sp(A),{\bf K})\to \cal C$ such
that \par $(i)$ $f(A)$ has an elementary value, if $f$ is a
polynomial.
\par Moreover, if $f, g \in C(sp(A),{\bf K})$, $a, b\in \bf K$, then
\par $(ii)$  $ \| f(A) \| = \| f \| $;
\par $(iii)$ $(af+bg)(A)=af(A)+bg(A)$;
\par $(iv)$ $(fg)(A)=f(A)g(A)$;
\par $(v)$ $f(A)$ is normal;
\par $(vi)$ ${\tilde f}(A)=[f(A)]^*$, in particular, $f(A)$ is selfadjoint
if and only if $f(sp(A))\subset \bf R$.}
\par {\bf Proof.} Due to Proposition 3.13 $sp (A)$ is the
compact subset in $\bf R$. Theorem 3.10 states that the set of all
polynomials with coefficients in $\bf K$ (with all possible orders
of constants and variables in each term for $\bf H$ or $\bf O$ and
possibly different orders of associated products in the case $\bf
O$) is dense in $C(sp(A),{\bf K})$. Consider a subalgebra $\cal A$
in $\cal C$, generated by an operator $A$, then due to Theorem 2.22
the element $f(A)$ of the algebra $\cal C$ is normal, since in
$C(\Lambda ,{\bf K})$ each element $z$ commutes with its adjoint.
Due to the same theorem there exists a mapping $f\mapsto f(A)$,
which is isometric.
\par {\bf 16. Note.} Apart from the complex case in general from the noncommutativity of
of the algebra $\bf K$ the property $f(A)B=Bf(A)$ for each element
$B\in \cal C$ commuting with $A$ is not accomplished for arbitrary
function $f\in C(sp(A),{\bf K})$, since $bI$ is not necessarily
commuting with $B$ for each $b\in \bf K$.
\par {\bf 17. Proposition.} {\it If $A$ is a selfadjoint element
of a $C^*$-algebra $\cal C$, then the set $\{ f(A): f\in
C(sp(A),{\bf K}) \} $ is the least quasicommutative $C^*$-subalgebra
${\cal C}(A)$ in $\cal C$ containing $A$ and $I$, also each element
from ${\cal C}(A)$ is a limit of a sequence of polynomials of $A$,
where ${\bf K}=\bf H$ or ${\bf K}=\bf O$.}
\par {\bf Proof.} Due to Theorems 2.22 and 3.15 the mapping
$C(sp(A),{\bf K})\ni f\mapsto f(A)\in \cal C$ is the isometric
$*$-isomorphism, $id(A)=A$, where $id(z)=z$ is the identity function
on $sp(A)$, $z\in sp(A)$. From the completeness of the metric space
$C(sp(A),{\bf K})$ it follows, that its image $\{ f(A): f\in
C(sp(A),{\bf K}) \} $ is the complete metric space, consequently, it
is closed in $\cal C$. Thus, $I, A\in {\cal C}(A)$, each element
from ${\cal C}(A)$ is the limit of polynomials, but ${\cal C}(A)$
contains all polynomials of $I, A$.
\par {\bf 18. Proposition.} {\it If $\cal C$ is a $C^*$-algebra
over $\bf K$, also $\cal E$ is a $C^*$-subalgebra in $\cal C$, $B\in
\cal E$, where ${\bf K}=\bf H$ or ${\bf K}=\bf O$, then $sp_{\cal
C}(B)=sp_{\cal E}(B)$.}
\par {\bf Proof.} If $u\in sp_{\cal C}(B)$, then $uI-B$ has not
an inverse in $\cal C$, consequently, $uI-B$ has not an inverse in
$\cal E$, therefore, $u\in sp_{\cal E}(B)$. Thus, $sp_{\cal
E}(B)\subset sp_{\cal C}(B)$.
\par For the proof of the reverse statement it is sufficient to verify, that
if $A\in \cal E$ and $A$ has an inverse element $A^{-1}$ in $\cal
C$, then $A^{-1}\in \cal E$. Let at first $B$ be selfadjoint, then
in view of $0\notin sp_{\cal C}(A)$, the equation $f(t)=1/t$ defines
a continuous function on $sp_{\cal C}(A)$. In accordance with
Theorem 3.15 applied to $\cal C$, there exists $f(A)\in \cal C$. Due
to Proposition 3.17 there exists $f(A)\in \cal E$. In view of $t
f(t)=1$ for each $t\in sp(A)$, then in accordance with Theorem 3.15
$A f(A)=I$, that is, $A^{-1}=f(A)\in \cal E$.
\par Consider now not necessarily selfadjoint element
$A\in \cal C$, having an inverse element $Y\in \cal C$. Then $A^*\in
\cal E$ and it has an inverse element $Y^*\in \cal C$. From the fact
that the element $AA^*\in \cal C$ is selfadjoint and has the inverse
element $YY^*\in \cal C$, then due to the proof given above $YY^*\in
\cal E$. Thus, $A^{-1}=Y=(YY^*)A^*\in \cal E$.
\par {\bf 19. Note.} In accordance with Proposition 3.18 it can be omitted the index
$\cal C$ or $\cal E$ and it can be written simply $sp(B)$ for the
spectrum of the element $B$ independently of the containing it
$C^*$-algebra.
\par {\bf 20. Theorem.} {\it If $A$ is a selfadjoint element
of a $C^*$-algebra $\cal C$ over $\bf K$, where ${\bf K}=\bf H$ or
${\bf K}=\bf O$, also $f\in C(sp(A),{\bf K})$, then $sp(f(A))= \{
f(t): t\in sp (A) \} $.}
\par {\bf Proof.} Due to Propositions 3.17 and 3.18
it can be considered the spectrum $sp(f(A))$ relative to the
$C^*$-subalgebra ${\cal C}(A)$ in $\cal C$. On the other hand, the
mapping \par $C(sp(A),{\bf K})\ni f\mapsto f(A)\in {\cal
C}(A)\subset C(sp(A),{\bf K})$ \\
is the isometric $*$-isomorphism, therefore, the spectrum of the
element $f(A)$ is such that: $sp(f(A))= \{ f(t): t\in sp (A) \} $.
\par {\bf 21. Theorem.} {\it Each element $A$ from a $C^*$-algebra $\cal C$
over $\bf K$, where ${\bf K}=\bf H$ or ${\bf K}=\bf O$, is a finite
linear combination of unitary elements.}
\par {\bf Proof.} In view of the fact that the algebra $\cal C$ is also the vector
space over $\bf K$ it follows, that $\cal C$ is isomorphic with the
direct sum ${\cal C}={\cal C}_0i_0\oplus {\cal C}_1i_1\oplus ...
\oplus {\cal C}_mi_m$, where $ \{ i_0,i_1,...,i_m \} $ are
generators of the algebra $\bf K$, also the linear spaces ${\cal
C}_0, {\cal C}_1,...,{\cal C}_m$ over $\bf R$ are pairvise
isomorphic. At the same time each element $A_p\in {\cal C}_p$ is
selfadjoint due to Theorem 2.22 and Propositions 3.13, 3.17, where
$p=0,1,...,m$.
\par Thus, it is sufficient to consider the case, when an element $A\in {\cal
C}$ is selfadjoint. The renormalization  $A\mapsto A/\| A \| $ for
$A\ne 0$ leads to the consideration of the case $\| A \| \le 1$. In
this case $sp(A)\subset [-1,1]$ and it can be prescribed a function
$f\in C(sp (A),{\bf K})$ by the formula $f(t)=t+M(1-t^2)^{1/2}$,
where $M+{\tilde M}=0$, $|M|=1$, $M\in \bf K$. In view of
$t=[f(t)+{\tilde f}(t)]/2$, $f(t){\tilde f}(t)={\tilde f}(t)f(t)=1$
for each $t\in sp(A)$, then the operator $f(A)=:U\in \cal C$
satisfies the equalities $A=(U+U^*)/2$, $UU^*=U^*U=I$.
\par {\bf 22. Theorem.} {\it Let $\cal C$ and $\cal E$ be
$C^*$-algebras over $\bf K$, where ${\bf K}=\bf H$ or ${\bf K}=\bf
O$, also $\phi : {\cal C}\to \cal E$ be a $C^*$-isomorphism. Then
there are accomplished the following statements:
\par $(i)$ $sp(\phi (A))\subset sp(A)$ and $ \| \phi (A) \| \le \| A
\| $ for each $A\in \cal C$, in particular, $\phi $ is continuous.
\par $(ii)$ If $A$ is a selfadjoint element in $\cal C$ and
$f\in C(sp(A),{\bf K})$, then $\phi (f(A))=f(\phi (A))$.
\par $(iii)$ If $\phi $ is a $*$-isomorphism, then $ \| \phi
(A)  \| = \| A \| $ and $sp (\phi (A))=sp (A)$ for each $A\in \cal
C$, also $\phi (A)$ is a $C^*$-subalgebra in $\cal E$.}
\par {\bf Proof.} $(i)$. If $b\notin sp(A)$, then $bI-A$
has an inverse element $Y\in \cal C$. In view of $\phi (I)=I$,
$bI-\phi (A)$ has an inverse element $\phi (S)$ in $\cal E$.
Therefore, $b\notin sp(\phi (A))$, consequently, $sp(\phi
(A))\subset sp(A)$. Due to Proposition 3.13 there are accomplished
the equalities $\| A \| ^2= \| A^*A \| =r(A^*A)$ and  $ \| \phi (A)
\| ^2= \| \phi (A)^*\phi (A) \| = \| \phi (A^*A) \| =r(\phi (A^*A)$.
In view of $sp(\phi (A^*A))\subset sp (A^*A)$, then $r(\phi
(A^*A))\le r(A^*A)$, therefore, $ \| \phi (A) \| \le \| A \| $.
\par $(ii)$. If $ \{ p_n : n \} $ is a sequence of polynomials,
converging to $f$ uniformly on $sp(A)$, then due to $(i)$ and also
on $sp(\phi A))$, then there exist $\lim_{n\to \infty }\phi
(p_n(A))= \phi (f(A))$ and $\lim_{n\to \infty }p_n(\phi A)= f(\phi
(A))$. Therefore, from $\phi (p_n(A))=p_n(\phi (A))$ for each $n$,
since $\phi $ is the homomorphism, it follows Statement  $(ii)$.
\par $(iii)$. For a selfadjoint element $B$ from $\cal C$
due to Statement $(i)$ it is accomplished the inclusion $sp (\phi
(B))\subset sp (B)$. If these sets are different, then there exists
a function $0\ne f\in C(sp(B),{\bf K})$ such that $f|_{sp(B)}=0$.
Due to Part $(ii)$ of this theorem we get: $f(B)\ne 0$, $\phi
(f(B))=f(\phi (B))=0$, that contradicts to the bijectivity of the
isomorphism $\phi $. Then $sp (\phi (B))=sp (B)$ and $r(\phi
(B))=r(B)$ for each selfadjoint element $B\in \cal C$. If $A\in \cal
C$ and $B=A^*A$, then due to the proof given above $ \| A \|
^2=r(A^*A)=r(\phi (A^*A))= \| \phi (A) \| ^2$, $ \| \phi (A) \| = \|
A \| $. In view of the fact that $\cal C$ is complete as the metric
space, also $\phi : {\cal C}\to \cal E$ is the isometry, then $\phi
({\cal C})$ is closed in $\cal E$, it contains the unitу,
consequently, it is the $C^*$-subalgebra in $\cal E$. Due to
Proposition 3.18 the spectrum $\phi (A)$ in $\cal E$ is the same, as
in $\phi ({\cal C})$, also $sp(A)=sp_{\phi ({\cal C})} (\phi (A))$,
since $\phi : {\cal C}\to \phi ({\cal C})$ is the isomorphism.
\par {\bf 23. Theorem.} {\it If $\cal C$ and $\cal E$ are $C^*$-algebras
over $\bf K$, where ${\bf K}=\bf H$ or ${\bf K}=\bf O$, also $\phi $
are $*$-homomorphisms from $\cal C$ into $\cal E$, then $\phi ({\cal
C})$ is the $C^*$-subalgebra in $\cal E$.}
\par {\bf Proof.} In view of the fact that $\phi ({\cal C})$ is
the $*$-subalgebra in $\cal E$ it is sufficient to prove that $\phi
({\cal C})$ is closed in $\cal E$. Consider $B\in \cal E$ and a
sequence $ \{ A_n : n, A_n\in {\cal C} \} $ with $\lim_{n\to \infty
} \| B -\phi (A_n) \| =0$. Due to Theorem 3.21 it is sufficient to
consider the case of selfadjoint operators $B, A_1, A_2,...$.
Without restriction of the generality it can be taken a subsequence
and in the same notation we suppose that $ \| \phi (A_{n+1})-\phi
(A_n) \| <2^{-n}$, $n=1,2,...$. Let $f_n: {\bf R}\to
[-2^{-n},2^{-n}]$, such that $f_n(t)=t$ while $|t|\le 2^{-n}$. From
Theorem 3.22 it follows, that $\phi (A_{n+1})-\phi (A_n)=f_n(\phi
(A_{n+1}-A_n))=\phi (f_n(A_{n+1}-A_n))$. In view of $ \|
f_n(A_{n+1}-A_n) \| \le 2^{-n}$ the series $A_1+\sum_{n=1}^{\infty }
f_n(A_{n+1}-A_n)$ converges to the element $A\in \cal C$, also due
to continuity of $\phi $ there exists $$\phi (A)=\lim_{m\to \infty }
\{ \phi (A_1) + \sum_{n=1}^{m-1} \phi (f_n(A_{n+1}-A_n)) \} =$$
$$\lim_{m\to \infty } \{ \phi (A_1)+ \sum_{n=1}^{m-1}[\phi
(A_{n+1})-\phi (A_n)] \} =\lim_{n\to \infty } \phi (A_m),$$
consequently, $B\in \phi ({\cal C})$, that is, $\phi ({\cal C})$ is
complete and hence is the $C^*$-algebra.
\par {\bf 24. Definition.} An element $A$ from a $C^*$-algebra $\cal C$
over $\bf K$, where ${\bf K}=\bf H$ or ${\bf K}=\bf O$, is called
nonnegative, if $A$ is selfadjoint and $sp(A)\subset [0,\infty
)\subset \bf R$. Denote by ${\cal C}^+$ the set of all nonnegative
elements in $\cal C$.
\par {\bf 25. Lemma.} {\it If $A$ is a selfadjoint element
of a $C^*$-algebra $\cal C$ over $\bf K$, where ${\bf K}=\bf H$ or
${\bf K}=\bf O$, $b\in \bf R$, $b\ge \| A \| $, then $A\in {\cal
C}^+$ if and only if $ \| A -bI \| \le b$.} \par {\bf Proof.} Due to
Proposition 3.13 $sp (A)\subset [-b,b]$ and also $$ \| A -bI \|
=r(A-bI)=\sup_{t\in sp (A)} |t-b| =\sup_{t\in sp(A)} (b-t).$$
Therefore, $ \| A -bI \| \le b$ if and only if $sp (A)\subset
[0,\infty )$.
\par {\bf 26. Theorem.} {\it Let $\cal C$ be a $C^*$-algebra
over $\bf K$, where ${\bf K}=\bf H$ or ${\bf K}=\bf O$, then
\par $(i)$ ${\cal C}^+$ is closed in $\cal C$;
\par $(ii)$ $bA\in {\cal C}^+$, if $A\in {\cal C}^+$ and $b\in [0,\infty )$;
\par $(iii)$ $A+B\in {\cal C}^+$, if $A, B\in {\cal C}^+$;
\par $(iv)$ $AB\in {\cal C}^+$, if $A, B\in {\cal C}^+$ also $AB=BA$;
\par $(v)$ $A=0$, if $A\in {\cal C}^+$ and $-A\in {\cal C}^+$.}
\par {\bf Proof.} $(i)$. Due to Lemma 3.25 ${\cal C}^+=
\{ A\in {\cal C}: A=A^* \mbox{ and } \| A - \| A \| I \| \le \| A \|
\} $, consequently, ${\cal C}^+$ is closed in $\cal C$.
\par $(ii)$. If $A\in {\cal C}^+$ and $b\in [0,\infty )$, then $bA$
is selfadjoint and $sp (bA)= \{ bt: t\in sp(A) \} \subset [0,\infty
)$.
\par $(iii)$. Due to Lemma 3.25 $ \| A - \| A \| I \| \le \| A \| $
and $ \| B - \| B \| I \| \le \| B \| $. Therefore, $ \| A +B - ( \|
A \| + \| B \| )I \| \le \| A \| + \| B \| $ and from this Lemma for
$b=\| A \| + \| B \| $ it follows, that $A+B\in {\cal C}^+$.
\par $(iv)$. Each element $A, B, AB$ has the same spectrum in $\cal
C$ as in the quasicommutative $C^*$-subalgebra generated by $ \{
I,A,B \} $. Due to Theorem 3.20 $sp(AB)\subset \{ st: s\in sp(A),
t\in sp(B) \} \subset [0,\infty )$.
\par $(v)$. If $A, -A\in {\cal C}^+$, then $A$ is selfadjoint and
$sp(A)\subset [0,+\infty )\cap (-\infty ,0] = \{ 0 \} $, therefore,
$ \| A \| =r(A)=0$.
\par {\bf 27. Proposition.} {\it Suppose that $A$ is a selfadjoint element
of a $C^*$-algebra $\cal C$ over $\bf K$, where ${\bf K}=\bf H$ or
${\bf K}=\bf O$, also $f\in C(sp(A),{\bf K})$, then \par $(i)$
$f(A)\in {\cal C}^+$ if and only if $f(t)\ge 0$ for each $t\in
sp(A)$;
\par $(ii)$ $ \| A \| I+sA\in {\cal C}^+$ for $s=1$ and $s=-1$;
\par $(iii)$ $A$ can be taken in the form $A^+-A^-$, where $A^+, A^-\in
{ \cal C}^+$ and $A^+A^-=A^-A^+=0$, moreover, these conditions
define $A^+$ and $A^-$ uniquely, also $ \| A \| = \max ( \| A^+ \| ,
\| A^- \| )$.}
\par {\bf Proof.} $(i)$. Due to Theorem 3.20 $f(A)$ has
the spectrum $ \{ f(t): t\in sp (A) \} $, therefore, $f$ is
nonnegative on $sp(A)$, if $f(A)\in {\cal C}^+$. Vise versa, if
$f(t)\ge 0$ for each $t\in sp(A)$, then $f(A)$ is selfadjoint, since
$sp(f(A)) \subset [0,\infty )$.
\par $(ii)$. A function $f\in C(sp(A),{\bf K})$ given by the formula
$f(t):= \| A \| +st$ takes nonnegative values, then due to $(i)$
$f(A)\in {\cal C}^+$, that is, $ \| A \| I+sA\in {\cal C}^+$.
\par $(iii)$. For the continuous function $f(t)=t$, $f^+(t):=
\max (t,0)$, $f^-(t):=\max (-t,0)$ there are accomplished equalities
$f=f^+-f^-$, $f^+f^-=f^-f^+=0$. Then $A=A^+-A^-$, $A^+A^-=A^-A^+=0$,
where $A^+=f^+(A)$, $A^-=f^-(A)$. In accordance with Part $(i)$ the
elements $A^+$ and $A^-$ belong to ${\cal C}^+$. From the equality
$\| f \|= \max ( \| f^+ \| , \| f^- \| )$ it follows, that $ \| A \|
= \max ( \| A^+ \| , \| A^- \| )$. For the proof of the uniquiness
consider the decomposition $A=B-C$, where $B, C\in {\cal C}^+$ and
$BC=CB=0$, then $A^n=B^n=(-C)^n$ for each $n=1,2,...$, therefore,
$p(A)=p(B)+p(-C)$, when $p$ is a polynomial with zero constant term.
There exists a sequence $ \{ p_n: n \} $ of polynomials, which
converges to $f^+$ uniformly on $sp(A)\cup sp(B)\cup sp(-C)$,
moreover, $$f^+(A)=\lim_{n\to \infty } p_n(A)=\lim_{n\to \infty }
[p_n(B)+p_n(-C)]=f^+(B)-f^+(-C).$$ In view of $f^+(B)=B$,
$f^-(-C)=0$, then $B=f^+(A)=A^+$, $C=B-A=A^+-A=A^-$.
\par {\bf 28. Corollary.} {\it If $A\in {\cal C}$,
where $\cal C$ is a $C^*$-algebra over ${\bf K}=\bf H$ or ${\bf
K}=\bf O$, then $A$ is a $\bf K$-linear combination of at most, than
$2p$ terms from ${\cal C}^+$, where $p=4$ for $\bf H$ and $p=8$ for
$\bf O$.}
\par {\bf Proof.} Due to \S 3.21 $A$ is a $\bf
K$-linear combination of at most, than $p$ selfadjoint elements from
$\cal C$, also due to Proposition 3.27 each selfadjoint element is
the difference of two nonnegative elements.
\par {\bf 29. Lemma.} {\it  If $\cal C$ is a $C^*$-algebra over ${\bf K}=
\bf H$ or ${\bf K}=\bf O$, $A\in \cal C$, $-A^*A\in {\cal C}^+$,
then $A=0$.}
\par {\bf Proof.} An element $A$ has a decomposition in the form of the sum
$A=\sum_{s=0}^mi_sA_s$, where $ \{ i_0,i_1,...,i_m \} $ are
generators of the algebra $\bf K$, each $A_s$ is a selfadjoint
element from $\cal C$ for each $s$ (see \S 3.21). From
$sp(A_s)\subset \bf R$ it follows, that $sp(A_s^2)\subset [0,\infty
)$. In view of the fact that $AA^*$ is selfadjoint, it follows that
$sp(AA^*)\subset \bf R$. Due to Conditions of this lemma and
Definition 3.24 $sp (-A^*A)\subset [0,\infty )$. If $b\ne 0$ and
$b\in sp(A^*A)$, then $b\in \bf R$, also $\bf R$ is the centre of
the algebra $\bf K$, then $A^*A-bI$ and, consequently,
$b^{-1}A^*A-I$ is noninvertible. Show that $I-A^*A$ is invertible if
and only if $I-AA^*$ is invertible. Let $(I-A^*A)$ be invertible.
Verify that then $(I-AA^*)$ is invertible:
$(I-A^*A)[A^*((I-A^*A)^{-1}A)+I]=A^*((I-A^*A)^{-1}A)+I-(A^*A)
(A^*(I-A^*A)^{-1}A)-A^*A$. Then due to the alternativity of the
algebra $\bf K$ and Theorem 2.22 it is accomplished the equality
$(I-A^*A)[A^*((I-A^*A)^{-1}A)+I] = A^*[
(I-A^*A)^{-1}-(AA^*)(I-AA^*)^{-1}]A+I-A^*A=I$ and the analogous
equality is accomplished for the right multiplication on $(I-A^*A)$.
\par Thus, $sp (-AA^*)\subset sp
(-A^*A)\cup \{ 0 \} \subset [0,\infty )$, that is, the element
$-AA^*$ is nonnegative. Then $A^*A=(A_0 - \sum_{s=1}^m i_sA_s) (A_0+
\sum_{p=1}^mi_sA_s)$, since to each $A_s$ there corresponds a
selfadjoint function $f_s={\tilde f}_s$ from $C(sp(A),{\bf K})$ due
to Theorem 2.22. Therefore $A^*A+AA^*=2\sum_{s=0}^mA_s^2$,
consequently, $A^*A=2\sum_{s=0}^mA_s^2+(-AA^*)$, where all terms in
the right part of the equality are nonnegative. Thus, $A^*A$ and
$(-A^*A)$ are nonnegative, therefore, $A^*A=0$ due to Theorem 3.26.
Thus, $\| A \| ^2=\| A^*A \| =0$, consequently, $A=0$.
\par {\bf 30. Theorem.} {\it If $\cal C$ is a $C^*$-algebra over
${\bf K}=\bf H$ or ${\bf K}=\bf O$, then the following conditions
are equivalent: \par $(i)$ $A\in {\cal C}^+$; \par $(ii)$ $A=H^2$
for some $H\in {\cal C}^+$;  \par $(iii)$ $A=B^*B$ for some $B\in
\cal C$. At the same time under satisfaction of these conditions $H$
is unique. If $X$ is a Hilbert space over $\bf K$, also ${\cal C}$
is a subalgebra in $L_q(X)$, then these conditions are equivalent to
the following condition:
\par $(iv)$ $<Ax,x> \ge 0$ for each $x\in X$.}
\par {\bf Proof.} In the algebra $\bf K$ each polynomial has
a root (see Theorem 3.17 in \cite{luoystoc}). If $A\in {\cal C}^+$,
then the equation $f(t)=t^{1/2}$ defines a nonnegative continuous
function $f$ on $sp(A)\subset [0,\infty )$. For $H=f(A)$ it is
evident that $H\in {\cal C}^+$ and $H^2=A$ due to Theorem 2.22. That
is, from $(i)$ it follows $(ii)$, also from $(ii)$ it follows
$(iii)$.
\par Suppose that $A=B^*B$, where $B\in \cal C$. From
selfadjointness of $A$ it follows, that $A=A^+-A^-$ in accordance
with Proposition 3.27. Let $C:=BA^-$, then $C^*C=(A^-B^*)(BA^-)$.
From Theorem 2.22 applied to the $C^*$ subalgebra generated by $B$,
using $C(sp(B),{\bf K})$ and the alternativity of the algebra $\bf
K$, we get, that $C^*C=A^-(A^+-A^-)A^-=-(A^-)^3$. In view of $A^-\in
{\cal C}^+$ and $(A^-)^3$ has the spectrum $ \{ t^3: t\in sp(A^-) \}
$, then $-C^*C=(A^-)^3\in {\cal C}^+$. Due to Lemma 3.29 $C=0$,
consequently, $A^-=0$ in accordance with corollary 3.14, since $A^-$
is selfadjoint. Thus, $A=A^+\in {\cal C}^+$, hence from $(iii)$ it
follows $(i)$.
\par Verify now the uniquiness of $H$ in $(ii)$ for $A\in {\cal
C}^+$. Let $G\in {\cal C}^+$, $G^2=A$, also $H=f(A)$ (see above).
There exists a sequence of polynomials $p_n$, converging to $f$ on
$sp(A)$. Put $q_n(t)=p_n(t^2)$. From $sp(A)=sp(G^2)=\{ t^2: t\in
sp(G) \} $ it follows, that $\lim_{n\to \infty }q_n(t)=\lim_{n\to
\infty }p_n(t^2)=f(t^2)=t$ uniformly on $sp(G)$, consequently,
$$G=\lim_{n\to \infty }q_n(G)=\lim_{n\to \infty }p_n(G^2)=\lim_{n\to
\infty }p_n(A)=f(A)=H,$$ that is, $H$ is unique.
\par If ${\cal C}$ is a subalgebra in $L_q(X)$, then ${\cal
C}^+={\cal C}\cap (L_q(X))^+$. Therefore, for the proof of the
equivalence of $(i)$ and $(iv)$ it is sufficient to consider ${\cal
C}=L_q(X)$. Then the statement about equivalence follows from the
fact that, if $A\in L_q(X)$, then $<Ax,x> \ge 0$ for each $x\in X$
if and only if $A=A^*$ and $sp(A)\subset [0,\infty )$ (see Lemma
2.17 and Theorem 2.31).
\par {\bf 31. Note.} An element $H$ in Theorem 3.30 is called
a nonnegative root from $A\in {\cal C}^+$ and it is denoted by
$A^{1/2}$. An analogous procedure can be applied for producing
elements $A^b$ for others real values $b>0$ with the help of the
function $f_b(t):=t^b$, moreover, for an invertible $A$ it can be
taken arbitrary $b\in \bf R$ with $f_0(t):=1$, such that
$A^bA^c=A^{c+b}$ for each $b, c \in \bf R$ for an invertible $A\in
{\cal C}^+$.
\par A linear over $\bf R$ subspace ${\cal C}_h$ of all
Hermitian elements of a $C^*$-algebra $\cal C$ is closed in $\cal
C$, consequently, it is the Banach space. At the same time ${\cal
C}_h$ is partially ordered by the relation $A\le B$ if and only if
$B-A\in {\cal C}^+$. Then ${\cal C}^+ = \{ A\in {\cal C}_h: A\ge 0
\} $.
\par {\bf 32. Definitions and Notations.} Let ${\cal M}\subset
{\cal C}$, such that ${\cal M}^* = \cal M$, that is, $\cal M$ is a
selfadjoint algebra, also $\cal M$ is linear over $\bf K$ (that is,
has the structure of a vector space over $\bf K$) with $I\in \cal
M$, where ${\bf K}=\bf H$ or ${\bf K}=\bf O$. By ${\cal M}^+$ we
denote ${\cal M}\cap {\cal C}^+$. Due to Corollary 3.28 each element
from $\cal M$ is a $\bf K$-linear combination of elements from
${\cal M}^+$.
\par A functional $\rho : {\cal M}\to \bf K$ we call Hermitian, if
$\rho \in L_q({\cal M},{\bf K})$ and $\rho = \rho ^*$, that is,
$\rho (A^*)=(\rho (A))^{\tilde .}$.
\par A functional $\rho \in L_q({\cal M},{\bf K})$ we call
nonnegative, if $\rho (A)\ge 0$ for each $A\in {\cal M}^+$. If also
$\rho (I)=1$, then a nonnegative functional is called a state.
\par Each functional $\rho \in L_q({\cal M},{\bf K})$ we write in the form:
$\rho (A) = \sum_{s\in \{ i_0, i_1,..., i_m \} } s \rho _s(A)$,
where $\rho _s(A)\in \bf R$ for each $A\in \cal M$ and each $s\in \{
i_0, i_1,..., i_m \}$, $ \{ i_0, i_1,..., i_m \} $ is the set of
standard generators of the algebra $\bf K$. Put also $\rho
_{bs}(A):=b\rho _s(A)$ for each $b\in \bf R$.
\par {\bf 33. Lemma.} {\it Let $\cal M$ and $\bf K$ be the same as in
\S 3.32. \par $(i)$. If $\rho \in L_l({\cal M},{\bf K})$, then $\rho
_s(A)=\rho _{qs}(qA)$ for each $s, q \in \{ i_0,...,i_m \} $, $A\in
\cal M$.
\par $(ii)$. If $\rho \in L_r({\cal M},{\bf K})$, then $\rho
_s(A)=\rho _{sq}(Aq)$ for each $s, q \in \{ i_0,...,i_m \} $, $A\in
\cal M$.
\par $(iii)$. If a functional $\rho \in L_q({\cal M},{\bf K})$
is Hermitian, then $\rho _s(A)=0$ for each $A = A^* \in \cal M$ and
each $s\in \{ i_1,i_2,...i_m \} $.}
\par {\bf Proof.} $(i,ii)$. Let $\rho \in L_l({\cal M},{\bf K})$,
Then $\rho (aA)=a\rho (A)$ for each $a\in \bf K$ and $A\in \cal M$.
In particular, while $a\in \{ i_0, i_1,..., i_m \} $ we get $\rho
(aA) =\sum_{q\in \{ i_0,...,i_m \} }q \rho _q(aA) = \sum_{s\in \{
i_0,...,i_m \} } as\rho _s(A)$, consequently, $\rho _{as}(aA)=\rho
_s(A)$ for each $s\in \{ i_0,...,i_m \} $. In the case $L_r$ the
proof is analogous.
\par $(iii)$. From the Hermiticity of a functional it follows the equality
\par $\rho (A)=\rho _{i_0}(A_{i_0})+ \sum_{s\in \{ i_1,...,i_m \} } s
\rho _s(A_s)=\rho (A^*)=\rho _{i_0}(A_{i_0}) - \sum _{s\in \{
i_1,...,i_m \} } s \rho _s(A_s)$ \\
for each $A=A^*$, consequently, $\rho _s(A)=0$ for each $s\in \{
i_1,...,i_m \} $.
\par {\bf 34. Lemma.} {\it If $\rho \in L_l({\cal M},{\bf K})$
or $\rho \in  L_r({\cal M},{\bf K})$, then $\| \rho \| \le (p)^{1/2}
\| \rho _e \| $, where $p=4$ for ${\bf K}=\bf H$ and $p=8$ for ${\bf
K}=\bf O$. If in addition a functional $\rho $ is Hermitian, then
$\| \rho \| = \| \rho _e \| =\sup \{ \rho (A): A=A^*\in {\cal M} \}
$.}
\par {\bf Proof.} Due to Lemma 3.33 $\rho _s(A)=\rho _e(sA)$
for $\rho \in L_l({\cal M},{\bf K})$ or $\rho _s(A)=\rho _s(As)$ for
$\rho \in L_r({\cal M},{\bf K})$ for each $A\in \cal M$. Therefore,
$\| \rho _s \| = \| \rho _e \| $ for each $s\in \{ i_0,...,i_m \} $,
where $ \| \rho \| := \sup_{0\ne A\in \cal M} |\rho (A)| / \| A \|
$, but $ |\rho (A)|^2=\sum_{s\in \{ i_0,...,i_m \} } \rho _s^2(A)\le
p\max_s \rho _s^2(A).$ If a functional $\rho $ is Hermitian, then
$\rho _s(A)=0$ for each $A=A^*$ and each $s\in \{ i_1,...,i_m \} $,
consequently, $|\rho (A)|=|\rho _e(A)| $ for $A=A^*$. It remains to
prove, that $ \| \rho \| =\sup \{ \rho (A): A=A^*\in {\cal M} \} $.
For each $\epsilon >0$ there exists $A\in {\cal M}$ with $ \| A \|
=1$, such that $|\rho (A)|> \| \rho \| - \epsilon $. For a suitable
number $b\in \bf K$ with $|b|=1$ we get: $\| \rho \| - \epsilon < \|
\rho \| = \rho (bA) = (\rho (bA))^{\tilde .} =\rho ((bA)^*)$. Take
the real part $Re(bA)=:(bA)_0$ of the element $bA$, then $ |Re
(bA)|\le 1$ and $\rho ((bA)_0)> \| \rho \| - \epsilon $. Thus, $\|
\rho \| \le \sup \{ \rho (A): A=A^*\in {\cal M}, \| A \| \le 1 \} $,
the inverse inequality is evident.
\par {\bf 35. Lemma.} {\it If a functional $\rho $ is nonnegative, then
$\rho $ is Hermitian.}
\par {\bf Proof.} If a functional $\rho $ is nonnegative,
also $A=A^*\in {\cal M}$, then $\rho ( \| A \| I +cA)\ge 0$ for
$c=1$ and $c=-1$, since $\| A \| I + cA \in {\cal M}^+$, also $\rho
(A)\in \bf R$, since $ \rho (A) = [\rho (\| A \| I+A) + \rho ( \| A
\| I-A)]/2$, consequently, $\rho $ is the Hermitian functional.
\par {\bf 36. Note.} A real vector space
${\cal M}_h$ consisting of all selfadjoint elements from $\cal M$ is
a partially ordered vector space with a positive cone ${\cal M}^+$
and the unit ordering $I$. If a functional $\rho $ belongs to
$L_l({\cal M},{\bf K})$ or $L_r({\cal M},{\bf K})$, then it is
Hermitian if and only if $\rho |_{{\cal M}_h}$ is $\bf R$-linear,
also each $\bf R$-linear functional on ${\cal M}_h$ has a unique
extension up to a Hermitian left or right $\bf K$-linear functional
on $\cal M$ with the help of relations of Lemma 3.33.
\par  All nonnegative left (or
right) $\bf K$-linear functionals on $\cal M$ form a cone ${\cal
P}_l$ (or ${\cal P}_r$ respectively) in the real vector space,
consisting of all Hermitian functionals from $L_l({\cal M},{\bf K})$
or $L_r({\cal M},{\bf K})$ respectively. At the same time ${\cal
P}_s\cap (-{\cal P}_s)=\{ 0 \} $, where $s=l$ or $s=r$, since $\cal
M$ is the $\bf K$-linear span of the set ${\cal M}^+$. The set of
all Hermitian functionals is partially ordered: $\rho _1\le \rho _2$
if and only if $\rho _2-\rho _1\in {\cal P}_s$, where $s=l$ or
$s=r$.
\par If $X$ is a Hilbert space over $\bf K$, then
for $x\in X$ the formula $w_x(A):=<x,A(x)>$ defines the right $\bf
K$-linear functional $w_x\in L_r({\cal M},{\bf K})$ on the $\bf
K$-linear subspace ${\cal M}$ in $L_q(X)$. At the same time
$w_x(I)=\| x \| ^2$, consequently, $w_x$ is nonnegative. A
functional $w_x$ is a state, if $ \| x \| =1$.
\par It is necessary to note, that if $\rho \in L_l({\cal M},{\bf K})\cap
L_r({\cal M},{\bf K})$ or $A\in L_l({\cal M},{\cal N})\cap L_r({\cal
M},{\cal N})$, then it does not mean a $\bf K$-linearity of a
functional $\rho $ or an operator $A$ in the same sense, that as in
the usual case of the commutative field of complex numbers $\bf C$
due to Lemma 3.33, since ${\cal M}$ is a $\bf K$-linear span of the
cone ${\cal M}^+$ of nonnegative elements, where $\cal M$ and $\cal
N$ are algebras over $\bf K$. Indeed, the equality $\rho
(a(bx))=a\rho (bx)=a\rho (xb)=a\rho (x)b$ is not guaranteed, as well
as $\rho (a(bx))=(ab)\rho (x)$ for arbitrary elements $a, b \in \bf
K$ and a vector $x$. By a right or left $\bf K$-linear functional or
operator over $\bf K$ we shall mean a functional or an operator $A$,
which is simultaneously right or left $\bf K$-linear, for example,
the operator $A$ on $l_2({\bf K})$ such that $A(e_j) = \sum_i
a_{j,i}e_i$ with real coefficients $a_{j,i}$, where $\{ e_j:
j=1,2,... \} $ is the standard orthonormal basis of a separable
Hilbert space over $\bf K$, $e_j=(0,...,0,1,0,...)$ with the unity
on the $j$-th place and others coordinate equal to zero.
\par If $\cal C$ is a $C^*$-subalgebra in $L_q(X)$, also $\cal M$ is a
selfadjoint subalgebra in $\cal C$, $I\in \cal M$, then $w_x|_{\cal
M}$ is a nonnegative functional on $\cal M$. States of an algebra
$\cal M$ arising by this way from the unit vectors in $X$ are called
vector states of the algebra $\cal M$.
\par Due to Lemma 3.33 the family of functionals belonging to $L_l({\cal M},{\bf
K})\cap L_r({\cal M},{\bf K})$ is nonvoid and is characterized by
the conditions $\rho _s(A_si_s)=\rho _s(i_sA_s)=i_s\rho _1(A_s)$ for
each $s\in \{ 0, 1,..., m \} $, where $A=\sum_{s=0}^m A_si_s$ is the
selfadjoint decomposition of the element $A\in {\cal M}$, $A_s\in
{\cal M}_0$ for each $s\in \{ 0, 1,..., m \} $, $\rho _1(A_s)\in \bf
R$, $A_s^*=A_s$.
\par {\bf 37. Proposition.} {\it If $\rho $ is a nonnegative
right or left $\bf K$-linear functional on a $C^*$-algebra $\cal C$
over $\bf K$, then $|\rho (B^*A)|^2 \le \rho (A^*A) \rho (B^*B)$ for
each $A, B\in {\cal C}$.}
\par {\bf Proof.} If $A\in \cal C$, then $A^*A\in {\cal
C}^+$, consequently, $\rho (A^*A)\ge 0$. The equation $<A,B>:=\rho
(A^*B)$ for each $A, B\in \cal C$ defines the scalar product (may be
without Property $(2)$ from \S 2.16). Then the statement of this
proposition follows from the Cauchy-Schwarz-Bunyakovskii inequality
for the scalar product $<A,B>$ on $\cal C$. That is, it remains to
prove the Cauchy-Schwarz-Bunyakovskii inequality for a vector space
$X$ over $\bf K$ with a function $<x,y>$, satisfying Conditions
2.16(1,3-7) for each $x, y \in X$. \par From Conditions 2.16(1,3-7)
it follows, that $<x,y>=\sum_{p,q}<x_p,y_q>i_p^*i_q$. Therefore, due
to the alternativity of the algebra $\bf K$, it is accomplished the
identity $<xa,xa>=<x,x>$ for each $x\in X$ and $a\in \bf K$ with
$|a|=1$. On the other hand, for given $x, y\in X$ there exists $a\in
\bf K$ with $|a|=1$, such that $<xa,y>=|<x,y>|\in [0,\infty )$. Then
without a restriction of the generality we can suppose, that
$<x,y>\ge 0$. If $<x,y>=0$, then the inequality $|<x,y>|^2\le <x,x>
<y,y>$ is trivial. Consider the case $<x,y> >0$ and the function
$f(b):=<xb+y,xb+y>$, where $b\in \bf R$. Due to Properties
2.16(1,3-7) $f(b)=b^2<x,x> + 2b<x,y>+ <y,y>$ for $<x,y>\in \bf R$.
In view of $<z,z>\ge 0$ for each $z\in X$, then this function $f$ is
nonnegative for each real numbers $b$ if and only if a discriminant
is nonpositive, that is, $<x,y>^2 -<x,x> <y,y>\le 0$. Thus,
$|<x,y>|^2\le <x,x> <y,y>$ for each $x, y \in X$.
\par {\bf 38. Theorem.} {\it If $\cal M$ is a selfadjoint
subspace over $\bf K$ in a $C^*$-algebra $\cal C$ and $I\in \cal M$,
where ${\bf K}=\bf H$ or ${\bf K}=\bf O$, then a left or right $\bf
K$-linear functional $\rho $ is nonnegative if and only if $\rho $
is bounded and $\| \rho \| = \rho (I)$.}
\par {\bf Proof.} Suppose at first, that $\rho $
is nonnegative, consequently, it is Hermitian. Due to the right or
left $\bf K$-linearity of the functional there exists $a\in \bf K$
such that $a\rho (A)\ge 0$ and let $H = Re(aA) := (aA)_0\in {\cal
C}_0$. Then $\| H \| \le \| A \| $, $H\le \| H \| I\le \| A \| I$,
$\| A \| \rho (I)-\rho (H)=\rho (\| A \| I-H) \ge 0$, therefore,
$|\rho (A)|=\rho (aA)=(\rho (aA))^{\tilde .}=\rho (A^*{\tilde
a})=\rho ((aA+A^*{\tilde a})/2)=\rho (H)\le \rho (I) \| A \| $. This
shows, that $\rho $ is bounded and $\| \rho \| \le \rho (I)$, also
the inverse inequality is evident.
\par Vise versa suppose that $\rho $ is bounded and $\| \rho \| =\rho
(I)$. It is sufficient to consider the case, in which $ \| \rho \|
=\rho (I)=1$. For $A\in {\cal M}^+$  consider $\rho (A) = a+Mb$,
where $a, b\in \bf R$, $M + {\tilde M} = 0$, $M\in \bf K$, $|M|=1$.
It is necessary to verify, that $a\ge 0$ and $b=0$. For a small
positive number $c$ there is the inclusion: $sp (I-cA)= \{ 1-ct:
t\in sp (A) \} \subset [0,1]$, since $sp(A)\subset [0,\infty )$,
Therefore, $\| I-cA \| =r(I-cA)\le 1$, consequently, $1-ca\le
|1-c(a+Mb)|=|\rho (I-cA)|\le 1$, Thus, $a\ge 0$. Consider
$B_n:=A-aI+MnbI$ for each natural number $n=1,2,...$, then $\| B_n
\| ^2=\| B_n^*B_n \| = \| (A-aI)^2+n^2b^2I \| \le \| A-aI \|
^2+n^2b^2$, consequently, $(n^2+2n+1)b^2=|\rho (B_n)|^2\le \| A-aI
\| ^2+n^2b^2$, $n=1,2,...$, thus, $b=0$.
\par {\bf  39. Note.} Due to Theorem 3.38 right and left $\bf K$-linear
states $\rho $ form a set ${\cal L}({\cal M})$ contained in the
sphere of the unit radius with the centre at zero in the dual space
${\cal M}^{\star }$ of all right and left $\bf K$-linear
functionals. At the same time ${\cal M}^{\star }$ is the $\bf
R$-linear space. The base of the weak $*$ topology in ${\cal
M}^{\star }$ form subsets of the form \par $ W_c(\eta
;A_1,...,A_n):= \{ |\rho (A_p)-\eta (A_p)|<c, p=1,...,n \} $, \\
where $c>0$, $\eta \in {\cal M}^{\star }$. Relative to the weak $*$
topology ${\cal L}({\cal M})$ is the compact topological Tychonoff
space, moreover, it is $\bf R$-convex, since \par ${\cal L}({\cal
M})= \{ \rho \in {\cal M}^{\star }: \rho (I)=1, \rho (A)\ge 0
\quad \forall A\in {\cal M}^+ \} $, \\
also $B({\bf K},0,1)^{\omega _0}$ is
compact in the Tychonoff product topology, where $B({\bf K},z,r)$ is
the ball in $\bf K$ of the radius $r$ with the centre at $z$,
$\omega _0$ is the first countable ordinal.
\par {\bf 40. Theorem.} {\it If $\cal C$ is a $C^*$-algebra
over ${\bf K}=\bf H$ or ${\bf K}=\bf O$, $\cal M$ is a selfadjoint
subspace over $\bf K$ in $\cal C$ with the unit $I\in \cal M$. If
$A\in \cal M$, $a\in sp(A)$, then there exists a right and left $\bf
K$-linear state $\rho $ on $\cal M$ such that $\rho (A)=a$.}
\par {\bf Proof.} For each $b, c, q\in \bf K$ there are accomplished
inclusions: $ab+c\in sp(Ab+cI)$, $qa+c\in sp(qA+cI)$, consequently,
$|ab+c|\le \| Ab+cI \| $, $|qa+c|\le \| qA+cI \| $. Then equations
$\rho ^0(Ab+cI)=ab+c$ and $\rho ^0(qA+cI)=qa+c$ uniquely define a
right and left $\bf K$-linear functional $\rho ^0$ on $\bf K$-vector
subspace $Y$ over $\bf K$, generated by vectors of the form $ \{
qA+cI, Ab+cI: b, c, q\in {\bf K} \} $, $\rho ^0(A)=a$, $\rho
^0(I)=1$, $\| \rho ^0 \| =1$. Due to Lemma 3.33 and Note 3.36 the
functional $\rho ^0_0$ on ${\cal M}_0$ is uniquely reconstructed
from the right and left $\bf K$-linear functional $\rho ^0$ on $\cal
M$. On the other hand, the functional $\rho ^0_0$ is realvalued and
$\bf R$-linear on $\bf R$-linear space $Y_0$, by the Hahn-Banach
theorem the functional $\rho ^0_0$ has the extension from $Y_0$ up
to the $\bf R$-linear functional $\rho _0$ on ${\cal M}_0$,
consequently, the desired right and left $\bf K$-linear state $\rho
$ on $\cal M$ exists.
\par {\bf 41. Theorem.} {\it Let $\cal C$ be a $C^*$-algebra over
${\bf K}=\bf H$ or ${\bf K}=\bf O$ with the unit $I\in \cal C$, also
$\cal M$ be a selfadjoint subspace in $\cal C$, $I\in \cal M$, $A\in
\cal M$.
\par $(i)$. If $\rho (A)=0$ for each right and left $\bf K$-linear state
$\rho $ on $\cal M$, then $A=0$.
\par $(ii)$. If $\rho (A)\in \bf R$ for each right or left $\bf K$-linear
state $\rho $ on $\cal M$, then $A=A^*$. \par $(iii)$. If $\rho
(A)\ge 0$ for each right or left $\bf K$-linear state $\rho $ on
$\cal M$, then $A\in {\cal M}^+$. \par $(iv)$. If $A$ is a normal
element from $\cal M$, then there exists such right and left $\bf
K$-linear state $\rho $ on $\cal M$, that $|\rho (A)|= \| A \| $.}
\par {\bf Proof.} $(i)$. Due to Theorem 3.40, $sp(A)= \{ 0
\} $, if $A$ is selfadjoint and $\rho (A)=0$ for each state $\rho $
on ${\cal M}$. For an arbitrary $A$ we can use the decomposition $A
= \sum_{s=0}^mA_si_s$, where each element $A_s$ is selfadjoint. Then
$\rho (A) = \sum_{s=0}^m \rho (A_s)i_s$ and $\rho (A_s)\in \bf R$
for each $s$, but then $\rho (A_s)=0$ for each $s$. Also due to the
beginning of the proof $A_s=0$ for each $s$, consequently, $A=0$.
\par $(ii)$. If $\rho (A)\in \bf R$ for each state $\rho $ on ${\cal M}$,
then $\rho (A-A^*)=\rho (A) - (\rho (A))^{\tilde .}=0$, that is,
$A-A^*=0$.
\par $(iii)$. If $\rho (A)\ge 0$ for each state $\rho $ on ${\cal M}$,
then $A$ is selfadjoint due to Part $(ii)$, also $sp(A) \subset \bf
R^+$ in accordance with Theorem 3.40, therefore, $A\in {\cal M}^+$.
\par $(iv)$. If $A$ is normal, then $r(A)= \| A \| $, therefore,
$sp(A)\ni a$ with $|a|=\| A \| $, where $a\in \bf K$. Due to Theorem
3.40 $a=\rho (A)$ for some state $\rho $ on ${\cal M}$, also then
$|\rho (A)|= \| A \| $.
\par {\bf 42. Notations.} For a subset $Y$ in a dual space
${\cal M}^{\star }$ (which is linear over $\bf R$) we denote by
${\overline {co}} _{\bf R}(Y) : = \{ z=ax+by: a+b=1, 0\le a\le 1,
0\le b\le 1, x, y\in Y \} $ the closed hull over $\bf R$ for $Y$.
\par {\bf 43. Lemma.} {\it Suppose that $\cal C$ is a $C^*$-algebra
over ${\bf K}=\bf H$ or ${\bf K}=\bf O$, also let $\cal M$ be
selfadjoint over $\bf K$ subspace in $\cal C$ with $I\in \cal M$,
${\cal L}({\cal M})$ be a set of all states on $\cal M$. If $ \| H
\| = \sup \{ |\rho (H)|: \rho \in {\cal L}({\cal M}) \} $ for each
selfadjoint $H\in \cal M$, then ${\overline {co}} _{\bf R}[({\cal
L}({\cal M})\cup - {\cal L}({\cal M})]$ it is the family of all
Hermitian right and left $\bf K$-linear functionals from the dual
space ${\cal M}^{\star }$.}
\par {\bf Proof.} For each $\bf K$-linear Hermitian functional $\rho $
it is accomplished the equality: \par $\rho ((bA)^*)=\rho
(A^*{\tilde b})=\rho (A^*){\tilde b}=(\rho (A))^{\tilde .}{\tilde
b}=(b\rho (A))^{\tilde .}$ \\ for each $A\in \cal M$ and $b\in \bf
K$. In view of the additivity of $\rho $ its Hermiticity is
preserved under taking $\bf K$-linear combinations of elements from
$\cal M$, also Hermitian functionals form the $\bf R$-linear space,
moreover, each Hermitian right and left $\bf K$-linear functional is
uniquely reconstructed from its part $\rho _0$ in the decomposition
(see Lemma 3.33 and Note 3.36). The set of all Hermitian functionals
in the unit ball with the centre at zero in ${\cal M}^{\star }$ is
convex over $\bf R$, it is weakly $*$ closed and contains in itself
$({\cal L}({\cal M})\cup - {\cal L}({\cal M})$, therefore, it
contains also  ${\overline {co}} _{\bf R}[({\cal L}({\cal M})\cup -
{\cal L}({\cal M})]$. \par Suppose that there exists a Hermitian
right and left $\bf K$-linear functional $\rho ^0\notin {\overline
{co}} _{\bf R}[({\cal L}({\cal M})\cup - {\cal L}({\cal M})]$ с $\|
\rho ^0\| \le 1$. Apply the Hahn-Banach theorem to the real parts
$\rho _0$ of Hermitian functionals $\rho $, then weakly continuous
right and left $\bf K$-linear functionals on ${\cal M}^{\star }$
arise from elements belonging to $\cal M$, since it is sufficient to
consider for such functionals ${\cal M}^{\star }_0$ and ${\cal
M}_0$. There exist $A\in \cal M$ and $b\in \bf R$, for which $Re
\rho ^0(A)>b$, also $Re \rho (A)<b$ for $\rho \in {\overline {co}}
_{\bf R}[({\cal L}({\cal M})\cup - {\cal L}({\cal M})]$. For each
Hermitian $\bf K$-linear functional $\rho $ and the real part $H$ of
the element $A\in \cal M$ it is accomplished: $\rho (H)=[\rho
(A)+\rho (A^*)]/2=Re \rho (A)$, therefore, $\rho ^0(A)>b$, $\rho
(H)\le b$ for $\rho \in {\overline {co}} _{\bf R}[({\cal L}({\cal
M})\cup - {\cal L}({\cal M})]$. Thus, $|\rho (H)|\le b$ for $\rho
\in {\cal L}({\cal M})$ and $b<\rho ^0(H)\le \| H \| =\sup \{ |\rho
(H)|: \rho \in {\cal L}({\cal M}) \} \le b$, that leads to the
contradiction.
\par {\bf 44. Lemma.} {\it Let be given a $C^*$-algebra $\cal C$ over
${\bf K}=\bf H$ or ${\bf K}=\bf O$, let also $\cal M$ be a $\bf
K$-vector subspace in $\cal C$, then $\cal M$ is selfadjoint.}
\par {\bf Proof.} In view of the fact that $\cal M$ is the vector
space over $\bf K$ it follows that it has the decomposition into the
direct sum ${\cal M}=\sum_{p=0}^m{\cal M}_pi_p$, where ${\cal M}_p$
and ${\cal M}_q$ are pairwise isomorphic spaces over $\bf R$ for
different $p$ and $q$, $ \{ i_0,...,i_m \} $ is the set of standard
generators of the algebra $\bf K$. At the same time $A\in \cal M$
can be decomposed into the sum $A=\sum_{p=0}^mA_pi_p$, where each
element $A_p$ is Hermitian. Due to Theorem 3.22(ii) the subalgebra
over $\bf K$, generated by the element $A_p$ is isomorphic with
$C(sp(A_p),{\bf K})$, moreover, $sp(A_p)\subset \bf R$. But each
real valued function $f_s$ on $sp(A_p)$ commutes with $i_q$ for each
$q$, consequently, $A_p$ commutes with all generators $i_q$. In view
of $A^*=A_0-\sum_{p=1}^mi_pA_p= A_0-\sum_{p=1}^mA_pi_p$, then
$A^*=(m-1)^{-1} \{ -A +\sum_{p=1}^m i_p(Ai_p^*) \} $, where $m=3$
for $\bf H$ and $m=7$ for $\bf O$.
\par {\bf 45. Theorem.} {\it If $\cal C$ is a $C^*$-algebra over ${\bf K}=\bf H$
or ${\bf K}=\bf O$, $\cal M$ is a $\bf K$-vector subspace in $\cal
C$, $I\in \cal M$, then each Hermitian right and left $\bf K$-linear
functional $\rho $ on $\cal M$ can be expressed in the form $\rho =
\rho ^+ - \rho ^-$, where $\rho ^+$ and $\rho ^-$ are nonnegative
right and left $\bf K$-linear functionals on $\cal M$ and $ \| \rho
\| = \| \rho ^+ \| + \| \rho ^- \| $. If ${\cal M}=\cal C$, then
these conditions define $\rho ^+$ and $\rho ^-$ uniquely.}
\par {\bf Proof.} Without any restriction of the generality suppose that
$ \| \rho \| =1$. Due to Lemma 3.44 the set $\cal M$ is selfadjoint,
that is, $A^*\in \cal M$ for each $A\in \cal M$. Due to Theorem 3.41
$ \| A \| = \sup \{ |\eta (A)|: \eta \in {\cal L}({\cal M}) \} $ for
each $A=A^*\in \cal M$, since $\| \eta \| =1$ for each $\eta \in
{\cal L}({\cal M})$. In accordance with Lemma 3.43 $\rho \in
{\overline {co}} _{\bf R}[({\cal L}({\cal M})\cup - {\cal L}({\cal
M})]$. On the other hand, the subset $F:=\{ a\xi -b \eta : \xi ,
\eta \in {\cal L}({\cal M}), a, b \in [0,\infty ), a+b=1 \} $ is
convex and it is contained in ${\overline {co}} _{\bf R}[({\cal
L}({\cal M})\cup - {\cal L}({\cal M})]$, moreover, $F\supset [({\cal
L}({\cal M})\cup - {\cal L}({\cal M})]$. Moreover, $F$ is weakly $*$
compact, since it is the continuous image for $[{\cal L}({\cal
M})]^2\times [0,1]$ relative to the mapping $(\eta , \xi ,a)\mapsto
a\eta -(1-a)\xi \in {\cal M}^{\star }$. Therefore, $F={\overline
{co}} _{\bf R}[({\cal L}({\cal M})\cup - {\cal L}({\cal M})]$. Thus,
$\rho =a\eta -b\xi $ for some $\eta , \xi \in {\cal L}({\cal M})$,
$a, b \in [0,1]$, $a+b=1$. Then take $\rho ^+=a\eta $ and $\rho
^-=b\xi $, consequently, $ \| \rho ^+ \| + \| \rho ^- \| =a+b=1= \|
\rho \| $.
\par Let now ${\cal M}=\cal C$. Prove the uniquiness
of the decomposition $\rho =\eta -\xi =\eta '- \xi '$ with
nonnegative right and left $\bf K$-linear functionals. For a given
$\epsilon >0$ choose a selfadjoint element $H\in \cal C$, for which
$\rho (H)> \| \rho \| -\epsilon ^2/2$ and let $B=(I-H)/2$. Then
$0\le B\le I$, $\eta (I)+\xi (I)= \| \eta \| + \| \xi \| =\| \rho \|
<\rho (H)+\epsilon ^2/2=\eta (H)-\xi (H)+\epsilon ^2/2$, $\eta
(I-H)+\nu (I-H)<\epsilon ^2/2$, $\eta (B)+\xi (I-B)<\epsilon ^2/4$.
In view of the fact that functionals $\eta $ and $\xi $ are
nonnegative, also $B, I-B\in {\cal C}^+$, then $0\le \eta
(B)<\epsilon ^2/4$ and $0\le \xi (I-B)<\epsilon ^2/4$. For $A\in
\cal C$ the Cauchy-Schwarz-Bunyakovskii inequality (see \S 3.37)
gives $|\eta (BA)|^2=|\eta (B^{1/2}(B^{1/2}A))|^2\le \eta (B)\eta
((A^*B^{1/2})(B^{1/2}A))\le \epsilon ^2 \| A \| ^2/4$, as well as
$|\xi ((I-B)A)|^2=\le \xi (I-B)\xi
((A^*(I-B)^{1/2})((I-B)^{1/2}A))\le \epsilon ^2 \| A \| ^2/4$. From
these inequalities and analogous arguments for $\eta '$ and $\xi '$
we get: $|\eta (BA)|\le \epsilon \| A \| /2$, $|\eta '(BA)|\le
\epsilon \| A \| /2$, $|\xi ((I-B)A)|\le \epsilon \| A \| /2$, $|\xi
'((I-B)A)|\le \epsilon \| A \| /2$. From $\eta -\eta '=\xi - \xi '$
it follows, that $\eta (A)- \eta '(A)=\eta (BA)-\eta '(BA)+\xi
((I-B)A) -\xi '((I-B)A)$ and therefore $|\eta (A)-\eta
'(A)|<2\epsilon \| A \| $. Due to the arbitrariness of $\epsilon
>0$ from this it follows, that $\eta =\eta '$, also then and $\xi =\xi '$.
\par {\bf 46. Corollary.} {\it If $\cal C$ is a $C^*$-algebra over ${\bf K}=\bf H$ or
${\bf K}=\bf O$, $\cal M$ is a $\bf K$-vector subspace in $\cal M$,
$I\in \cal M$, then each bounded quasilinear functional is a $\bf
K$-linear combination of at most $n$ states on $\cal M$, where $n=8$
for $\bf H$, $n=16$ for $\bf O$.}
\par {\bf Proof.} Each bounded quasilinear
functional $\rho $ on $\cal M$ has the form $\rho =\sum_{p=0}^m\rho
_pi_p$, where $\rho _p(A)\in \bf R$ and $\rho _p (A^*)=\rho _p(A)$
for each $A\in \cal M$ and $p=0,1,...,m$. but due to Theorem 3.45
each functional $\rho _p$ is a real linear combination of two states
on $\cal M$.
\par {\bf 47. Note.} In view of the fact that the space of states ${\cal
L}({\cal M})$ is convex over $\bf R$ and it is weakly $*$ compact,
then by the Krein-Milman theorem ${\cal }({\cal M})$ coincides with
the closure of the convex over $\bf R$ hull ${\overline {co}}_{\bf
R}({\cal P}({\cal M}))$ of the set ${\cal P}({\cal M})$ of all its
extreme points. Elements from ${\cal P}({\cal M})$ are called pure
states for $\cal M$, also a weak $*$ closure $cl ({\cal P}({\cal
M}))$ is called the space of pure states for $\cal M$. In the
general case ${\cal P}({\cal M})$ is not closed in ${\cal M}^{\star
}$ and the space of pure states has states, which are not pure.
\par {\bf 48. Theorem.} {\it  Suppose that $\cal C$ is a $C^*$-algebra over
${\bf K}=\bf H$ or ${\bf K}=\bf O$ with the unit $I\in \cal C$,
$\cal M$ is a $\bf K$-vector subspace in $\cal C$, $I\in \cal M$,
$A\in \cal M$. \par $(i)$. If $\rho (A)=0$ for each pure state $\rho
$ on $\cal M$, then $A=0$. \par $(ii)$. If $\rho (A)\in \bf R$ for
each pure state $\rho $ for $\cal M$, then $A=0$.
\par $(iii)$. If $\rho (A)\ge 0$ for each pure state $\rho $
for $\cal M$, then $A\in {\cal M}^+$. \par $(iv)$. If $A$ is a
normal element, then there exists a pure state $\rho ^0$ such that
$|\rho ^0(A)|= \| A \| $.} \par {\bf Proof.} Statements $(i-iii)$
follow from Theorem 3.41, since each state is a weak $*$ limit of
convex combinations of pure states. Due to Lemma 3.44 the set $\cal
M$ is selfadjoint. Let now $A$ be a normal element, then by Theorem
3.41 there exists a constant $b\in \bf K$ and a state $\eta $ for
$\cal M$ such that $\eta (A)=b$, $|b|=\| A \| $. Let now $\hat T$ be
a weak $*$ continuous right or left $\bf K$-linear functional on
${\cal M}^{\star }$ such that $\rho (T)=:{\hat T}(\rho )$, $a\in \bf
K$, $|a|=1$ and at the same time $\eta (aA)=|b|=\| A \| $. It is
known that, if $X$ is a nonvoid compact subset in a locally convex
over $\bf R$ space $Y$, also $\rho _0$ is a continuous $\bf
R$-linear functional on $Y$, then there exists an extreme point
$x_0$ in $X$ such that $\rho _0(x)\le \rho _0(x_0)$ for each $x\in
X$. Therefore, there exists $\rho ^0\in {\cal P}({\cal M})$ such
that \par $\| A \| \ge |\rho ^0(A)|\ge Re (aA)^{\hat .}(\rho ^0)\ge
\sup \{ Re (aA)^{\hat .} (\rho ): \rho \in {\cal L}({\cal M}) \}$
\par $ \ge Re(aA)^{\hat .}(\eta )= Re [\eta (aA)]=\| A \| $.
\par {\bf 49. Theorem.} {\it If $\cal C$ is a $C^*$-algebra
over ${\bf K}=\bf H$ or ${\bf K}=\bf O$, $I\in \cal C$, $\cal M$ is
a $\bf K$-linear subspace in $\cal C$, ${\cal L}_0\subset {\cal
L}({\cal M})$, then the following four conditions are equivalent:
\par $(i)$. If $A\in \cal M$ and $\rho (A)\ge 0$ for each $\rho
\in {\cal L}_0$, then $A\in {\cal M}^+$.
\par $(ii)$. $\| H \| = \sup \{ |\rho (H)|: \rho \in {\cal L}_0 \} $
for each selfadjoint element $H$ in $\cal M$.
\par $(iii)$. ${\overline {co}}_{\bf R}({\cal L}_0)={\cal L}({\cal
M})$. \par $(iv)$. ${\cal P}({\cal M})$ is contained in a weak $*$
closure $ cl({\cal L}_0)$ in ${\cal M}^{\star }$.}
\par {\bf Proof.} Due to Lemma 3.44 the set $\cal M$
is selfadjoint. If $H=H^*\in \cal M$, then we define $b\le \| H \| $
by the formula $b=\sup \{ |\rho (H)|: \rho \in {\cal L}_0 \} $ and
mention that $\rho (bI+cH)=b+c\rho (H)\ge 0$ for $c=1$ and $c=-1$
and for each $\rho \in {\cal L}_0$. If this is satisfied, than
$bI+cH\in {\cal M}^+$, $-bI\le H\le bI$, consequently, $ \| H \| \le
b$, therefore, $b=\| H \| $, that is, from $(i)$ it follows $(ii)$.
Due to Lemma 3.43 and Theorem 3.45 from $(ii)$ it follows $(iii)$,
also from Theorem 3.48 it follows the implication
$(iv)\Longrightarrow (i)$.
\par {\bf 50. Corollary.} {\it If $H$ is a selfadjoint operator
acting on a Hilbert space $X$ over ${\bf K}=\bf H$ or ${\bf K}=\bf
O$, then $ \| H \| =\sup \{ |<Hx,x>|: x\in X, \| x \| =1 \} $. If
$\cal M$ is a $\bf K$-vector subspace in $L_q(X)$ containing $I$ and
${\cal L}_0$ is a set of all vector states on $\cal M$, then ${\cal
P}({\cal M})\subset cl ({\cal L}_0)$ and ${\cal L}({\cal
M})={\overline {co}}_{\bf R}({\cal L}_0)$.}
\par {\bf Proof.} Due to Lemma 2.44 the set $\cal M$ is selfadjoint.
If $A\in \cal M$ and $\rho (A)\ge 0$ for each $\rho \in {\cal L}_0$,
then $<Ax,x>\ge 0$ for each $x\in X$ and thus, $A\in {\cal M}^+ =
{\cal M}\cap L_q(X)^+$ due to Theorem 3.30. Then in accordance with
Theorem 3.49 ${\overline {co}}_{\bf R}({\cal L}_0)={\cal L}({\cal
M})$, ${\cal P}({\cal M})\subset cl ({\cal L}_0)$ and \par $\| H \|
=\sup \{ |\rho (H)|: \rho \in {\cal L}_0 \} =\sup \{ |<Hx,x>| : x\in
X, \| x \| =1 \} $ \\  for each selfadjoint element $H$ from $\cal
M$.
\par {\bf 51. Definitions and notes.} Let a $C^*$-algebra
$\cal C$ be given, also $X$ be a Hilbert space over ${\bf K}=\bf H$
or ${\bf K}=\bf O$. Then a $*$-homomorphism $\phi : {\cal C}\to
L_q(X)$ is called a representation of a $C^*$-algebra on a Hilbert
space $X$. If in addition a $*$-homomorphism $\phi $ is bijective,
that is, it is a $*$-isomorphism, then it is called an exact
representation.
\par Due to our convention about $*$-homomorphisms $\phi (I)=I$
preserves the unit of the algebra, therefore, from Theorem 3.22 it
follows, that $\| \phi (A) \| \le \| A \| $ for each $A\in \cal C$,
consequently, $\phi $ is continuous, $ \| \phi (A) \| = \| A \| $,
if $\phi $ is the exact representation. In general consider the set
$ \{ A\in {\cal C}: \phi (A) =0 \} $, which is the closed two sided
ideal in $\cal C$, that is called the kernel of the representation
$\phi $. If there exists a vector $x\in X$ such that the $\bf
K$-vector space $ \phi ({\cal C})x:= \{ B(x): B=\phi (A), A\in {\cal
C} \} $ is dense in $X$, then the representation $\phi $ is called
cyclic, also $x$ is called the cyclic vector (or the generating
vector) for $\phi $.
\par From the definition of the $*$-homomorphism $\phi $ it follows, that
$\phi \in L_r({\cal C},L_q(X))\cap L_l({\cal C},L_q(X))$. For a
subset $V$ in $X$ we denote by $[V]$ the least closed $\bf K$-vector
subspace in $X$ containing $V$.
\par As an example consider a Banach space $L^p
((S,{\cal F},{\hat \mu }),{\bf K})$ of all \\ $\mu $-measurable
functions $f: S\to {\bf K}$ with the finite norm $ \| f \| _p<\infty
$, where $S$ is a set, $\cal F$ is a $\sigma $-algebra of its
subsets, $\hat \mu $ is a $\sigma $-finite measure (which may be
noncommutative or nonassociative, see \S 2.23) on $(S,{\cal F})$
with values in $\bf K$, $1\le p\le \infty $,
$$\| f \|_p^p :=\sum_{m,l} \int_S|f_m(x)|^p |\mu _{m,l}|(dx)$$
for $1\le p<\infty $, $|\mu _{m,l}|$ denotes the variation of real
valued measure $\mu _{m,l}$, $$\| f \|_{\infty }:=\max_{m,l}
ess-\sup_{|\mu _{m,l}|}|f_m(x)|.$$  As a $C^*$-algebra we take
$L^{\infty } ((S,{\cal F},{\hat \mu }),{\bf K})$ with the pointwise
multiplication of functions and with the involution induced by the
operation of the conjugation in the algebra $\bf K$. As the Hilbert
space $X$ we also take $L^2 ((S,{\cal F},{\hat \mu }),{\bf K})$,
supposing that
$$(f,q):=\sum_{m,l}{\tilde f}_m(x)q_m(x)|\mu _{m,l}|(dx).$$
Then each function $f\in L^{\infty } ((S,{\cal F},{\hat \mu }),{\bf
K})$ generates the operator $M_f\in L_q(X)$ of the multiplication on
the function $f$, $M_f(g):=f(x)g(x)$ for each $g\in X$, $x\in S$.
Thus, the representation  $f\mapsto M_f$ is exact.
\par {\bf 52. Proposition.} {\it If $\rho $ is a right and left
$\bf K$-linear state of a $C^*$-algebra $\cal C$ over ${\bf K}=\bf
H$ or ${\bf K}=\bf O$, then the set ${\cal L}_{\rho } := \{ A\in
{\cal C}: \rho (A^*A)=0 \} $ is the closed ideal in $\cal C$, also
$\rho (B^*A)=0$ for each $A\in {\cal L}_{\rho }$ and $B\in \cal C$.
At the same time there exists a scalar product $<A+{\cal L}_{\rho
},B+ {\cal L}_{\rho }>= \rho (B^*A)$ on the quotient space ${\cal
C}/{\cal L}_{\rho }$.}
\par {\bf Proof.} From the proof of Proposition 3.37
it follows, that it can be taken an internal product $<A,B>_0 :=
\rho (A^*B)$ for each $A, B\in \cal C$, moreover, ${\cal L}_{\rho }
= \{ A\in {\cal C}: <A,A>_0=0  \} $, since the functional $\rho $ is
positive, and hence it is Hermitian. Then ${\cal L}_{\rho }$ is the
$\bf K$-linear subspace in $\cal C$ and the equation $<A+{\cal
L}_{\rho },B+{\cal L}_{\rho }> = <A,B>_0=\rho (A^*B)$ defines the
scalar product on ${\cal C}/{\cal L}_{\rho }$. If $A\in {\cal
L}_{\rho }$, $B\in \cal C$, then due to the analog of the
Cauchy-Schwarz-Bunyakovskii inequality which was proved in \S 3.37
over $\bf K$ (instead of $\bf C$) $|\rho (B^*A)|^2\le \rho
(B^*B)\rho (A^*A)=0$, consequently, $\rho (B^*A)=0$. The element
$B^*B$ is selfadjoint, therefore, in accordance with Theorems 2.22
and 2.31 it commutes with each generator $i_p$ of the algebra $\bf
K$. Each element $A\in \cal C$ has the decomposition
$A=\sum_pi_pA_p$, where each $A_p$ is selfadjoint, then
$((B^*B)A)^*A=(BA)^*(BA)$ due to the alternativity of the algebra
$\bf K$, consequently, $\rho ((BA)^*(BA))=\rho (((B^*B)A)^*A)=0$ for
$A\in {\cal L}_{\rho }$ and $B\in \cal C$, that is, $BA\in {\cal
L}_{\rho }$. Thus, ${\cal L}_{\rho }$ is the closed left ideal in
$\cal C$, since the functional $\rho $ is continuous.
\par We call ${\cal L}_{\rho }$ the left kernel of the functional $\rho $.
\par {\bf 53. Theorem.} {\it  If $\rho $ is a left and right $\bf K$-linear
state of a $C^*$-algebra $\cal C$ over ${\bf K}=\bf H$ or ${\bf
K}=\bf O$, then there exists a cyclic left and right $\bf K$-linear
representation $\pi _{\rho }$ of the algebra $\cal C$ on a Hilbert
space $X$ over $\bf K$, also there exists the unit cyclic vector
$x_{\rho }$ for $\pi _{\rho }$, such that $\rho =w_{x_{\rho }}\circ
\pi _{\rho }$, that is, $\rho (A)= <\pi _{\rho }(A)x_{\rho },x_{\rho
}>$ for each $A\in \cal C$.}
\par {\bf  Proof.} Consider the kernel ${\cal L}_{\rho }=\rho ^{-1}(0)$ of the
functional $\rho $, then ${\cal C}/{\cal L}_{\rho }$ is the
preHilbert space over $\bf K$ relative to the internal product $<A +
{\cal L}_{\rho }, B + {\cal L}_{\rho }> = \rho (B^*A)$, $A, B\in
\cal C$ due to Proposition 3.52. Its completion is the Hilbert space
$X=X_{\rho }$. In view of the fact that $X$ is the vector space over
$\bf K$, it follows that for each vector $x\in X$ there exists a
vector $x_0\in X_0$, such that the sets ${\bf K}x$ and $x{\bf K}$
and ${\bf K}x_0$ and $x_0\bf K$ coincide, moreover, $ax_0=x_0a$ for
each $a\in \bf K$, where $X=X_0\oplus i_1X_1\oplus ... \oplus
i_mX_m$, $X_0,...,X_m$ are Hilbert spaces over $\bf R$, $X_p$ is
isomorphic with $X_0$ for each $p$, $m=3$ for $\bf H$, $m=7$ for
$\bf O$. Therefore, without loss of generality it can be taken
$x_{\rho }\in X_0$.
\par If $A, B_1, B_2\in \cal C$ and $B_1+{\cal L}_{\rho }=B_2+{\cal
L}_{\rho },$ then $B_1-B_2\in {\cal L}_{\rho }$ and $AB_1-AB_2\in
{\cal L}_{\rho }$, since ${\cal L}_{\rho }$ is the left ideal in
$\cal C$. Then the equation $\pi (A) (B+{\cal L}_{\rho }) = AB+{\cal
L}_{\rho }$ defines the unique operator $\pi (A)$ acting on the
preHilbert space ${\cal C}/{\cal L}_{\rho }$, moreover, $\pi
(a_1A_1+a_2A_2) = a_1\pi (A_1)+a_2\pi (A_2)$, that is, $\pi (A)$ is
left $\bf K$-linear by $A$, $\pi (A) (B_1b_1+B_2b_2) = \pi
(A)(B_1b_1) + \pi (A)(B_2b_2)$ for each $a_1, a_2, b_1, b_2 \in \bf
K$, $A_1, A_2, B_1, B_2\in \cal C$, that is, $\pi (A) B$ by the
variable $B$ is right $\bf K$-linear in the case ${\bf K} = \bf H$
and it is right $\bf K$-alternative in the case of ${\bf K}=\bf O$.
It can be used the decomposition ${\cal C}={\cal C}_0\oplus i_1{\cal
C}_1\oplus ... \oplus i_m{\cal C}_m$, where each element $G\in {\cal
C}_p$ is selfadjoint for each $p=0,1,...,m$, such that $G$ commutes
with each generator $i_q$ of the algebra $\bf K$ (see also the proof
of Proposition 3.52). Due to the alternativity of the algebra $\bf
O$ it is accomplished the equality $\rho (C^*C)=0$ for
$C:=A(Bb)-(AB)b$ for each $A, B\in \cal C$ and $b\in \bf O$,
consequently, $A(Bb)-(AB)b\in {\cal L}_{\rho }$, that is, $\pi (A)B$
induces the right $\bf O$-linear operator by $B$ on $X$ over $\bf
O$, as well as over $\bf H$. Consider now three arbitrary elements
$A, B, F\in \cal C$ and put $C=(AB)F-A(BF)$, then due to the
alternativity of the algebra $\bf O$ and the aforementioned
decomposition of the algebra $\cal C$ it is accomplished the
identity $\rho (C^*C)=0$, that is, $C\in {\cal L}_{\rho }$. In the
case of the algebra $\cal C$ over the skew field $\bf H$ from the
associativity of the skew field $\bf H$ of quaternions it follows,
that $C=0$. Thus, for ${\bf K}=\bf H$ and ${\bf K}=\bf O$ the
element $A$ induces the multiplicative operator $\pi (A)$ on $X$,
that is, $\pi (AB)=\pi (A)\pi (B)$, since $(AB)F-A(BF)\in {\cal
L}_{\rho }$. It can be applied this also to $C= (AB)^*F-B^*(A^*F)$
for the algebra $\cal C$ over $\bf O$, then  $<\pi _{\rho
}(A)(B+{\cal L}_{\rho }),C+{\cal L}_{\rho }>=<AB+{\cal L}_{\rho
},C+{\cal L}_{\rho }>=\rho ((AB)^*C)=\rho (B^*(A^*C))=<B+{\cal
L}_{\rho },A^C+{\cal L}_{\rho }>=<B+{\cal L}_{\rho },\pi _{\rho
}(A^*)(C+{\cal L}_{\rho })>$ due to Proposition 3.37. For the
algebra $\cal C$ over $\bf H$ we consider the zero element $C=0$ and
the identity $<\pi _{\rho }(A)(B+{\cal L}_{\rho }),C+{\cal L}_{\rho
}>=<B+{\cal L}_{\rho },\pi _{\rho }(A^*)(C+{\cal L}_{\rho })>$ is
accomplished in this case also. Thus, $\pi (A^*)=(\pi (A))^*$ on
$X$.
\par Then $ \| A \| ^2I - A^*A = \| A^*A \| I -A^*A\in {\cal C}^+$,
consequently, $B^* (( \| A \| ^2I-A^*A)B)\in {\cal C}^+$, therefore,
$ \| A \| ^2 \| B+{\cal L}_{\rho } \| ^2 - \| \pi (A) (B+{\cal
L}_{\rho } \| ^2 = \| A \| ^2\rho (B^*B)- \rho (B^*((A^*A)B))= \rho
(B^*((\| A \| ^2I-A^*A)B))\ge 0$ for each $A, B\in \cal C$,
consequently, $\pi (A)$ is bounded, with $ \| \pi (A) \| \le \| A \|
$ and $\pi (A)$ can be extended by the continuity up to the bounded
right $\bf K$-linear operator $\pi (A)$ acting on $X$. In view of
$\pi (I)=I$ on ${\cal C}/{\cal L}_{\rho }$, it follows, that $\pi
_{\rho }(I)=I$ on $X_{\rho }$. Thus, the representation $\pi _{\rho
}$ is left and right $\bf K$-linear on $X$, that is, $\pi
(a_1A_1+a_2A_2)x = a_1(\pi (A_1)x)+a_2(\pi (A_2)x)$, $\pi
(A_1a_1+A_2a_2)x=((\pi (A_1)a_1)x)+((\pi (A_2)a_2)x)$ for each $A_1,
A_2 \in \cal C$, $a_1, a_2\in \bf K$, $x\in X$. Consider $x_{\rho }
:= I+{\cal L}_{\rho }$, evidently, that $x_{\rho }\in X_0$. Then
$\pi _{\rho }(A)x_{\rho }=\pi _{\rho }(A)(I+{\cal L}_{\rho
})=A+{\cal L}_{\rho }$, where $A\in \cal C$, consequently, $\pi
_{\rho }({\cal C})x_{\rho }$ is everywhere dense in $X$, that is,
the vector $x_{\rho }$ is cyclic. At the same time $\pi _{\rho
}(A_1a_1+A_2a_2)(I+{\cal L}_{\rho })=A_1a_1+{\cal L}_{\rho
}+A_2a_2+{\cal L}_{\rho }=\pi (A_1)(I+{\cal L}_{\rho })a_1+\pi
(A_2)(I+{\cal L}_{\rho })a_2$ for each $a_1, a_2\in \bf K$, $A_1,
A_2\in \cal C$. Moreover, $<\pi _{\rho }(A)x_{\rho },x_{\rho
}>=<A+{\cal L}_{\rho },I+{\cal L}_{\rho }>=\rho (A)$, where $A\in
\cal C$, in particular, $ \| x _{\rho } \| ^2=\rho (1)=1$. Thus, the
representation $\pi _{\rho }$ is left and right $\bf K$-linear. In
view of $x_{\rho }\in X_0$, it folows, that the composition
$w_{x_{\rho }}\circ \pi _{\rho }$ is left and right $\bf K$-linear.
\par The given above construction is the noncommutative nonassociative analog
for the $C^*$-algebra over ${\bf K}=\bf H$ and ${\bf K}=\bf O$ of
the construction used by Gelfand, Naimark and Segal in the
particular case of a $C^*$-algebra over the commutative associative
field of complex numbers.
\par {\bf 54. Proposition.} {\it Suppose that  $\rho $
is a right and left $\bf K$-linear state of a $C^*$-algebra $\cal C$
over ${\bf K}=\bf H$ or ${\bf K}=\bf O$, also $\pi $is a (right and
left $\bf K$-linear) cyclic representation of a $C^*$-algebra $\cal
C$ on a Hilbert space $X$ over $\bf K$, such that $\rho = w_x\circ
\pi $ for some cyclic vector $x$ for $\pi $. If $X_{\rho }$, $\pi
_{\rho }$ and $x_{\rho }$ are a Hilbert space over $\bf K$, a cyclic
representation and the unit cyclic vector obtained from $\rho $ by
the way of the construction of the proof of Theorem 3.53, then there
exists a right and left $\bf K$-linear isomorphism $U$ from $X_{\rho
}$ onto $X$ such that $x=U(\pi _{\rho }(A)U^*)$ for each $A\in \cal
C$.}  \par {\bf Proof.} For each element $A\in \cal C$ there are
accomplished the equalities $$\| \pi (A) x \| ^2= <\pi (A)x,\pi
(A)x>= <\pi (A^*A)x,x> =\rho (A^*A)$$  $$= <\pi _{\rho
}(A^*A)x_{\rho },x_{\rho }>= \| \pi _{\rho }(A)x_{\rho } \| ^2.$$
 If $A, B\in \cal
C$ and $\pi _{\rho }(A)x_{\rho } =\pi _{\rho }(B)x_{\rho }$, then
$\pi (A)x=\pi (B)x$. Then the equation $U_0 \pi _{\rho }(A) x_{\rho
} = \pi (A)x$ for each $A\in \cal C$ defines a norm preserving
operator $U_0$ from $\pi _{\rho } ({\cal C}) x_{\rho }$ onto $\pi
({\cal C})x$, which is left and right $\bf K$-linear, since without
restriction of the generality it can be chosen $x\in X_0$ and
$x_{\rho }\in (X_{\rho })_0$, for which $bx=xb$ and $bx_{\rho
}=x_{\rho }b$ for each $b\in \bf K$. In view of $[\pi _{\rho }({\cal
C})x_{\rho }]=X_{\rho }$ and $[\pi ({\cal C})x]=X$, where $[S]$
denotes the closure of the $\bf K$-linear span of a subset $S$ in
the corresponding $\bf K$-vector topological space, then $U_0$ has
an extension by the continuity up to the isomorphism $U$ from
$X_{\rho }$ on $X$, also $Ux_{\rho }=U_0(\pi _{\rho }(I)x_{\rho
})=\pi (I)x=x$. For $A, B\in \cal C$ it is accomplished the
equality: \par $U(\pi _{\rho }(A)\pi _{\rho }(B))x_{\rho }=U(\pi
_{\rho }(AB)x_{\rho })= \pi (AB)x=\pi (A)(\pi (B)x)=\pi (A)(U(\pi
_{\rho }(B)x_{\rho }))$. \\  In view of the fact that the set of
vectors of the form $\pi _{\rho }(B)x_{\rho }$ with arbitrary $B\in
\cal C$ forms the everywhere dense subset in $X_{\rho }$, then $U\pi
_{\rho }(A)=\pi (A)U$, consequently, $\pi (A)=U(\pi _{\rho
}(A)U^*)$.
\par {\bf 55. Note and Definition.} Representations $\phi $ and
$\psi $ (right and left $\bf K$-linear) for a $C^*$-algebra $\cal C$
on Hilbert spaces $X$ and $Y$ are called unitary equivalent, if
there exists an isomorphism $U$ from $X$ onto $Y$ such that $\psi
(A)=U(\phi (A)U^*)$ for each $A\in \cal C$.  \par  If $\rho $ is a
right and left $\bf K$-linear state of a $C^*$-algebra $\cal C$,
$\pi $ is a (right and left $\bf K$-linear) cyclic representation
for $\cal C$, $\rho =w_x\circ \pi $ for some unit cyclic vector $x$
for $\pi $, then from Proposition 3.54 it follows, that $\pi $ is
equivalent with $\pi _{\rho }$ obtained from  $\rho $ by the way of
the construction from the proof of Theorem 3.53. Moreover, an
isomorphism $U$ can be chosen such that $Ux_{\rho }=x$.
\par {\bf 56. Corollary.} {\it If $x$ is a unit vector in a Hilbert
space $X$ over ${\bf K}=\bf H$ or ${\bf K}=\bf O$, $\cal C$ is a
$C^*$-subalgebra in $L_q(X)$, also $\rho $ is a vector state
$w_x|_{\cal C}$ with $x\in X_0$, then the representation $\pi _{\rho
}$ obtained from $\rho $ by the way of the contruction of \S 3.53 is
equivalent to the representation $A\mapsto A|_{[{\cal C}x]}$ for
$\cal C$ on a Hilbert space $[{\cal C}x]$. The isomorphism $U:
X_{\rho }\to [{\cal C}x]$, which accomplishes this equivalence can
be chosen such that $Ux_{\rho }=x$.}
\par {\bf  Proof.} It follows from Proposition 3.54, since $x$ is the unit cyclic vector
for the representation $\pi : A\mapsto A|_{[{\cal C}x]}$, where
$\rho =w_x\circ \pi $.
\par {\bf 57. Proposition.} {\it If $A$ is a nonzero element
of a $C^*$-algebra $\cal C$ over ${\bf K}=\bf H$ or ${\bf K}=\bf O$,
then there exists a pure state $\rho $ for $\cal C$, such that $\pi
_{\rho }(A)\ne 0$, where $\pi _{\rho }$ is a representation obtained
from $\rho $ by the way of the construction from \S 3.53.}
\par {\bf Proof.} Due to Theorem 3.48 there exists a pure
state $\rho $ for $\cal C$, such that $\rho (A)\ne 0$, that is
equivalent to: $<\pi _{\rho }(A)x_{\rho },x_{\rho }>\ne 0$, from
which it follows, that $ \pi _{\rho }(A)\ne 0$.
\par {\bf 58. Notations.} Let $\cal C$ be a $C^*$-algebra over${\bf K}=\bf H$ or
${\bf K}=\bf O$, $\{ X_b: b \in {\bf B} \} $ is a family of Hilbert
spaces, also $\phi _b$ is a representation of a $C^*$-algebra $\cal
C$ on a Hilbert space $X_b$, $b\in \bf B$, such that the direct sum
$\bigoplus_{b\in \bf B}\phi _b(A)$ is the bounded linear operator
acting on the Hilbert space $\bigoplus_{b\in \bf B}X_b$. The mapping
$\phi : A\mapsto \bigoplus_b\phi _b(A)$ is the representation for
$\cal C$ on $\bigoplus_bX_b$. We call $\phi $ the direct sum of the
family $\{ \phi _b: b\in {\bf B} \} $ of representations for $\cal
C$, also we denote $\phi =\bigoplus_b\phi _b$.
\par {\bf 59. Theorem.} {\it Each $C^*$-algebra over ${\bf K}=\bf H$
or ${\bf K}=\bf O$ has an exact representation.}
\par {\bf Proof.} Consider any family ${\cal L}_0$ of
right and left $\bf K$-linear states of a $C^*$-algebra $\cal C$
containing all pure right and left $\bf K$-linear states. Put $\phi
:=\bigoplus \{ \pi _{\rho }: \rho \in {\cal L}_0 \} $, where $\pi
_{\rho }$ is a state obtained from $\rho $ by the construction of
the  proof of Theorem 3.53. If $A\in \cal C$, also $\phi (A)=0$,
then $\pi _{\rho }(A)=0$ for each $\rho \in {\cal L}_0$, since $\phi
(A)=\bigoplus_{\rho }\pi _{\rho }(A)$, in particular, $\pi _{\rho
}(A)=0$ for each pure state $\rho $ of the $C^*$-algebra $\cal C$,
moreover, $A=0$ due to Proposition 3.57, consequently, $\phi $ is
the exact representation for $\cal C$.
\par {\bf 60. Note.} If $\phi $ is an exact representation
of a $C^*$-algebra $\cal C$ over ${\bf K}=\bf H$ or ${\bf K}=\bf O$
on a Hilbert space $X$ over $\bf K$, then $\phi $ is isometrical,
also $\phi ({\cal C})$ is the $C^*$-subalgebra in $L_q(X)$ due to
Theorem 3.22.
\par Let $\cal L$ be the family of all right and left $\bf
K$-linear states of the $C^*$-algebra $\cal C$. For ${\cal L}_0=\cal
L$ it can be obtained the exact representation $F := \bigoplus_{\rho
\in \cal C} \pi _{\rho }$, which we call the universal
representation for $\cal C$. For $\eta \in \cal L$ there exists a
vector $x_{\eta }$ in $X_{\eta }$ such that $\eta =w_{x_{\eta
}}\circ \pi _{\eta }$ and $\| x_{\eta } \| =1$, that is, $\eta
=w_y\circ F$, where $y=\bigoplus_{\rho \in \cal C}y_{\rho }$,
$y_{\eta }=x_{\eta }$, $y_{\rho }=0$ for $\rho \ne \eta $. In view
of the fact that the mapping $\eta \mapsto \eta \circ F^{-1}$
transfers ${\cal L}={\cal L}({\cal C})$ into ${\cal L}(F({\cal
C}))$, where $F({\cal C})\subset L_q(X_F)$, $X_F:= \bigoplus_{\rho
\in \cal L}X_{\rho }$, then each right and left $\bf K$-linear state
for $F({\cal C})$ is a vector state.
\par {\bf 61. Corollary.} {\it  If $\cal C$ is a $C^*$-algebra over
${\bf K} = \bf H$ or ${\bf K} = \bf O$, $A\in {\cal C}^+$, $B\in
\cal C$, then $B^*(AB)\in {\cal C}^+$.}
\par {\bf Proof.} Due to Theorem 3.59 and Note 3.60, also
due to Theorem 3.30 there exists $A^{1/2}\in {\cal C}^+$, such that
$B^*(AB)=(A^{1/2}B)^*(A^{1/2}B)$, since $<B^*(A(B(x))),x> =
<A(B(x)),B(x)>=<Ay,y>\ge 0$, where $y=B(x)\in X$, also $\cal C$ is
considered as a subalgebra in $L_q(X)$ for a Hilbert space $X$ over
$\bf K$, $x\in X$.
\par {\bf 62. Proposition.} {\it Let $A$ and $B$ be selfadjoint elements
in a $C^*$-algebra $\cal C$. \par $(i)$. If $-B\le A\le B$, then $
\| A \| \le \| B \| $. \par $(ii)$. If $0\le A\le B$, then
$A^{1/2}\le B^{1/2}$. \par $(iii)$. If $0\le A\le B$ and $A$ is
invertible, then $B$ is invertible and $B^{-1}\le A^{-1}$.}
\par {\bf Proof.} In the algebra $\bf K$ each polynomial has roots (see Theorem 3.17
\cite{luoystoc}). $(i)$ follows from inequality $ - \| B \| I\le
-B\le A\le B\le \| B \| I$.
\par $(ii, iii)$. If $0\le A\le B$ and $A$ is invertible, then $A\ge bI$ for some
$b>0$, since $sp(A)\subset [a,\infty )$, where $a>0$ if and only if
$0\notin sp(A)$, that is, $A$ is invertible. At the same time $0\le
B^{-1/2}(AB^{-1/2})\le B^{-1/2}(BB^{-1/2})=I$, also $\|
B^{-1/2}(AB^{-1/2}) \| \le 1$ due to $(i)$. Thus, $$(iv)\quad \|
A^{1/2}B^{-1/2} \| = \| (A^{1/2}B^{-1/2})^*(A^{1/2}B^{-1/2}) \|
^{1/2} = \| B^{-1/2}(AB^{-1/2}) \| ^{1/2} \le 1,$$
 consequently,
$A^{1/2}(B^{-1}A^{1/2})\le I$, also $A^{1/2}\le
B^{1/4}IB^{1/4}=B^{1/2}$. That is, statements $(ii, iii)$ are proved
for invertible element $A$. In the general case for $A, B\in \cal C$
with $0\le A\le B$ it is accomplished the inequality $A+bI\le B+bI$
and $A+bI$ is invertible for each $b>0$. Then
$$(v)\quad (A+bI)^{1/2}\le (B+bI)^{1/2}.$$
Take $G\in {\cal C}^+$ and $f_b\in C(sp (G),{\bf K})$ by the formula
$f_b(t):=(t+b)^{1/2}$, then $f_b(G)\in {\cal C}^+$ and
$(f_b(G))^2=G+bI$. Thus, $(G+bI)^{1/2}=f_b(G)$, since $\lim_{b\to
0}f_b(t)=t^{1/2}$ converges uniformly on $sp(G)$, then $\lim_{b\to
0} \| (G+bI)^{1/2}-G^{1/2} \| =0$. Thus, due to Theorem 3.26 the
limit while $b$ converges to zero gives due to $(v)$ that
$A^{1/2}\le B^{1/2}$.
\par {\bf 63. Proposition.} {\it If ${\cal G}$ is a closed left ideal in
a $C^*$-algebra $\cal C$, then each element $S\in \cal G$ can be
expressed in the form $S=AB$ with $A\in \cal C$ and $B\in {\cal
G}\cap {\cal C}^+$.}
\par {\bf Proof.} If $F\in {\cal G}\cap {\cal
C}^+$, then $F^{1/2}\in {\cal G}\cap {\cal C}^+$, since $\cal G$
contains all $p(K)$ for each $K\in \cal G$ and each polynomial $p$
with the zero constant term. There exists a sequence of polynomials
$ \{ p_n(t): n\in {\bf N} \} $ converging uniformly to $t^{1/2}$ on
$sp (F)$ due to the Stone-Weierstrass theorem (see also \S 2.7 in
\cite{luoystoc}). In view of the fact that $\cal G$ is the closed
ideal, then $$F^{1/2} =\lim_np_n(F)\in \cal G.$$ For $S\in \cal G$,
let $H=(S^*S)^{1/2}$ and $B=H^{1/2}$, then $S^*S\in {\cal G}\cap
{\cal C}^+$, thus, $H, B\in {\cal G}\cap {\cal C}^+$. Take the
elements $A_n=S(I/n+H)^{-1/2}$, then $S=A_n(I/n+H)^{1/2}$, since $H$
commutes with each generator $i_p$ of the algebra $\bf K$, also the
algebra $\bf K$ is alternative. Therefore, \par $ \| A_m-A_n \| = \|
S[(I/m+H)^{-1/2}-(I/n+H)^{-1/2}] \| =$\par $ \|
[(I/m+H)^{-1/2}-(I/n+H)^{-1/2}]S^*] [S[(I/m+H)^{-1/2}-
(I/n+H)^{-1/2}] \| ^{1/2} =$  $\| [(I/m+H)^{-1/2}-(I/n+H)^{-1/2}]
(H^2[(I/m+H)^{-1/2}-(I/n+H)^{-1/2}] \| ^{1/2}$ \par $= \|
[f_{m,n}(H)]^2\| ^{1/2}= \| f_{m,n}(H) \| $ \\  due to Theorem 2.22
for the algebra generated by $S^*S$, also due to Corollary 3.46,
where $f_{m,n}$ is the continuous function on $(0,\infty )$
described by the formula $f_{m,n}(t)=t[(m^{-1}+t)^{-1/2}
-(n^{-1}+t)^{-1/2}]$. Thus, $$ \| A_m-A_n \| =\sup \{ |f_{m,n} (t)|:
t\in sp (H) \} \mbox{, moreover, }lim_{\min (m,n)\to \infty
}f_{m,n}(t)=0$$ converges uniformly on $sp (H)$, that is, $ \{ A_n:
n\in {\bf N} \} $ is the Cauchy sequence in $\cal C$. For
$A=\lim_nA_n$ it gives $S=AH^{1/2}=AB$.
\par {\bf 64. Corollary.} {\it Each closed two sided ideal
$\cal G$ in a $C^*$-algebra $\cal C$ over ${\bf K}=\bf H$ or ${\bf
K}=\bf O$ is selfadjoint. A closed two sided ideal $\cal S$ in $\cal
G$ is a two sided ideal in $\cal C$.}
\par {\bf Proof.} Due to the closedness of the left ideal in $\cal
C$ each element $S\in \cal G$ has the form $S=AB$, where $A\in \cal
C$, also $B=(S^*S)^{1/4}\in {\cal G}\cap {\cal C}^+$ due to
Proposition 3.63. In view of the fact that $\cal G$ is also the
right ideal in $\cal C$, then $\cal G$ is selfadjoint.
\par Let $S\in \cal S$, then $S\in \cal G$, also since $\cal S$ is the
left ideal in $\cal G$, then $S^*S\in {\cal S}\cap {\cal C}^+$.
Then, as in the proof of Proposition 3.63, ${\cal S}\cap {\cal C}^+$
contains the square root of each its term. Therefore,
$B^{1/2}=(S^*S)^{1/8}\in \cal S$. If $F\in \cal C$, then
$F(AB^{1/2})$ and $B^{1/2}F\in \cal G$, since $B^{1/2}\in \cal G$,
also $\cal G$ is a two sided ideal in $\cal C$. In view of the fact
that $\cal S$ is the left ideal in $\cal C$, also $B^{1/2}\in \cal
S$, then $FS=F(AB)=F(A(B^{1/2}B^{1/2}))\in \cal S$, that is, $\cal
S$ is the left ideal in $\cal C$. Thus, $AB^{1/2}\in \cal S$, also
since $\cal S$ is the right ideal in $\cal G$, then
$SF=A(BF)=A((B^{1/2}B^{1/2})F)\in \cal S$. This proves that $\cal S$
is the two sided ideal in $\cal C$.
\par {\bf 65. Lemma.} {\it Suppose that $\cal G$ is a closed
$C^*$-subalgebra (not necessarily with the unit) in a $C^*$-algebra
$\cal C$, and let $\Upsilon := \{ B\in {\cal G}\cap {\cal C}^+: \| B
\| < 1 \} $. \par $(i)$. If $B_1, B_2\in \Upsilon $, then there
exists an element $B\in \Upsilon $ such that $B_1\le B$ and $B_2\le
B$.
\par $(ii)$. If $F\in \cal G$ and $b>0$, then there exists an element $J\in
\Upsilon $ such that from $B\in \Upsilon $ and $B\ge J$ it follows,
that $ \| F-FJ \| <b$.}
\par {\bf Proof.} $(i)$. For given $B_1, B_2\in \Upsilon $
there exists $b>0$ such that $ \| (1+b)B_j\| \le 1$ for $j=1$ and
$j=2$. For each $n\in {\bf N} := \{ 1,2,3,... \} $ and each $0\le
t\le 1$ it is accomplished the inequality $t\le t^{1/n}$, thus,
$(1+b)B_j\le [(1+b)B_j]^{1/n}$. When $n$ has the form $2^q$ for some
$q\in \bf N$, then a repeated application of Proposition 3.62 leads
to the inequality $ [(1+b)(B_1+B_2)]^{1/n}\ge [(1+b)B_j]^{1/n}\ge
(1+b)B_j$. Then $B_1\le B$ and $B_2\le B$, where
$B=(1+b)^{-1}[(1+b)(B_1+B_2)]^{1/n}$ for $n=2^q$. Thus,
$(1+b)(B_1+B_2)\in \cal G$, $0\le (1+b)(B_1+B_2)\le 2I$, also the
function $t^{1/n}$ is the uniform limit on the segment $[0,2]$ of
polynomials without constant terms. In view of the fact that $\cal
G$ is closed, then $B\in \cal G$, moreover, $0\le B\le
(1+b)^{-1}2^{1/n}I$. When $n$ is sufficiently large, then
$(1+b)^{-1}2^{1/n}<1$, hence $B\in \Upsilon $.
\par $(ii)$. For a given $F\in \cal G$ and $b>0$, let $F=aQ$, where
$a> \| F \| $, such that $Q\in \cal G$ and $ \| Q \| <1$. Take
$J=(Q^*Q)^{1/n}$ in $\Upsilon $, where $n$ is sufficienly large,
that to satisfy the inequality $t(1-t)^{1/n}<b^2a^{-2}$ for each
$0\le t\le 1$. If $B\in \Upsilon $ and $B\ge J$, then $ \| F-FB \|
=a \| Q(I-B) \| = a \| Q((I-B)^2Q^*) \| ^{1/2}$, since $B\ge 0$,
also the algebra $\bf K$ is alternative, $\cal C$ has the embedding
into $L_q(X)$ for some Hilbert space $X$ over $\bf K$ due to Theorem
3.59. Then $0\le I-B\le I-J\le I$, thus, $0\le (I-B)^2\le I-B\le
I-J=I-(Q^*Q)^{1/n}$, such that $0\le Q((I-B)^2Q^*)\le
Q([I-(Q^*Q)^{1/n}]Q^*)$, consequently, \par $ \| Q((I-B)^2Q^*) \|
\le \| Q([I-(Q^*Q)^{1/n}]Q^*) \| =r(Q([I-(Q^*Q)^{1/n}]Q^*))= $
\par $=r(Q^*Q[I-(Q^*Q)^{1/n}])\le $  $\sup \{ t(1-t^{1/n}): 0\le t\le
1 \} < b^2a^{-2}$ \\  due to Propositions 3.13, 3.62 and Theorem
3.20, since $sp (Q^*Q)\subset [0,1]$. Therefore, $ \| F-FB \| =a \|
Q((I-B)^2Q^*) \| ^{1/2}<b$.
\par {\bf 66. Definition.} For a closed subset $\cal S$
in a $C^*$-algebra $\cal C$ over ${\bf K}=\bf H$ or ${\bf K}=\bf O$
we say that a net $ \{ F_v: v\in V \} $ of selfadjoint elements
$F_v\in \cal S$ is the increasing right approximation of the unit
for $\cal S$, if $$ lim_{v\in V} \| A-AF_v \| =0$$ for each $A\in
\cal S$, moreover, $0\le F_v\le F_u\le I$ for each $v\le u\in V$,
where $V$ is a directed set. Analogously left increasing and two
sided approximations of the unit are defined.
\par {\bf 67. Proposition.} {\it A closed left ideal in a $C^*$-algebra
$\cal C$ over ${\bf K}=\bf H$ or ${\bf K}=\bf O$ has an increasing
right approximation of the unit, also a closed $C^*$-subalgebra has
an increasing two sided approximation of the unit.}
\par {\bf Proof.} Let $\cal G$ be a closed
$C^*$-subalgebra in $\cal C$. The set $\Upsilon $ from Lemma 3.65 is
directed in the usual way by partial ordering on selfadjoint
elements from $\cal C$. Take $F_B=B$ for $B\in \Upsilon $, then we
take an increasing net in the ball of the unit radius with the
centre at zero in $\cal G$. Due to Lemma 3.65 $\lim_{B\in \Upsilon }
\| A-AF_B \| =0$ for each $A\in \cal G$, also since $\cal G$ is
selfadjoint, then $$\lim_{B\in \Upsilon } \| A-F_bA \| = \lim_{B\in
\Upsilon } \| A^*-A^*F_B \| =0.$$ Thus, $ \{ F_B: B\in \Upsilon \} $
is an increasing two sided approximation of the unit.
\par Let now $\cal S$ be a closed left ideal in $\cal C$.
For the closed right ideal ${\cal S}^* := \{ A^*: A\in {\cal S} \} $
we get  ${\cal S}\cap {\cal S}^*$, which is the closed
$C^*$-subalgebra in $\cal C$, that contains all selfadjoint elements
from $\cal S$. Due to the proof given above ${\cal S}\cap {\cal
S}^*$ has an increasing two sided approximation of the unit $ \{
F_v: v \in V \} $. Due to Proposition 3.63 each element $A\in \cal
S$ has the from $A=BG$ with $B\in \cal S$ and $G\in {\cal S}\cap
{\cal C}^+ \subset {\cal S}\cap {\cal S}^*.$ In view of $ \| A -AF_v
\| = \| B(G-GF_v) \| \le \| B \| \| G-GF_v \| $, then $ \{ F_v: v\in
V \} $ is the increasing right approximation of the unit for $\cal
S$.
\par {\bf 68. Theorem.} {\it  If a quasicommutative algebra $\cal A$ is contained
in $L_q(X)$ for a Hilbert space $X$ over ${\bf K}=\bf H$ or ${\bf
K}=\bf O$, also $\cal A$ contains the unit $I$, it is closed
relative to the weak operator topology, then $\cal A$ is
isometrically $*$-isomorphic with the algebra $C(S,{\bf K})$ for a
totally disconnected compact Hausdorff space $S$.}
\par {\bf Proof.} Due to Theorem 2.22 $\cal A$
is isometrically $*$-isomorphic with the algebra $C(S,{\bf K})$,
where $S$ is a compact Hausdorff topological space. This isomorphism
gives the ordering of the set of selfadjoint operators from $\cal A$
by the way of the pointwise ordering on the set of real valued
continuous functions on $S$. Each increasing net $ \{ f_a: a\in
\Upsilon \} $ from $C(S,{\bf K})$, which is bounded by the constant
$w$, there corresponds an increasing net $ \{ \mbox{ }_aA: a\in
\Upsilon \} $ of selfadjoint operators from $\cal A$, bounded from
above by the operator $wI$. Due to Lemma 2.22.4 $ \{ \mbox{ }_aA:
a\in \Upsilon \} $ has a supremum $A$ in $\cal A$. There exists a
function $f\in C(S,{\bf K})$ corresponding to the operator $A$,
moreover, $f= \sup_{a\in \Upsilon } f_a$, $f(t)\le w$ for each $t\in
S$. For the family $\Omega $ of all finite subsets from the directed
set $\Upsilon $ and each $\omega \in \Omega $ the function
$\max_{a\in \omega }f_a$ belongs to $C(S,{\bf K})$, moreover, $w\ge
\max_{a\in \omega }f_a$.
\par Thus, $S$ is characterized by the conditions, that each
bounded from above (below) net $ \{ f_a: a\in \Upsilon \} \subset
C(S,{\bf R})\subset C(S,{\bf K})$ has $\sup_{a\in \Upsilon } f_a\in
C(S,{\bf R})$ ($\inf_{a\in \Upsilon } f_a\in C(S,{\bf R})$
respectively). Let $U$ be an open subset in $S$, $cl_S(U)$ be its
closure in $S$. Consider the family of all functions ${\cal
G}\subset C(S,{\bf R})$ such that $0\le f\le 1$ on $S$ and $f(t)=0$
for each $t\in S\setminus cl_S(U)$ and each $f\in \cal G$. Let
$f_0:=\sup \{ f: f\in {\cal G} \} $. Then $f_0(t)\le 1$ for each
point $t\in S$. If $t\in U$, then there exists a function $f\in \cal
G$ such that $f(t)=1$, consequently, $f_0(t)=1$ for each point $t\in
U$, also this means that for each point $t\in cl_S(U)$ due to the
continuity of the function $f_0$. If otherwise $t\in S\setminus
cl_S(U)$, then there exists a function $g\in C(S,{\bf R})$ such that
$0\le g\le 1$ on $S$, $g(t)=0$ for each $t\in S\setminus cl_S(U)$,
$g(t)=1$ for each $t\in cl_S(U)$. That is, $g=\sup \{ f: f\in {\cal
G} \} $ and $f_0\le g$. Thus, $f_0(t)=1$ on $cl_S(U)$ and $f_0(t)=0$
on $S\setminus cl_S(U)$. Due to the continuity of the function $f_0$
the set $cl_S(U)$ is open, consequently, $S$ is totally
disconnected.
\par {\bf 69. Definitions.} Let ${\cal F}\subset L_q(X)$
for a Hilbert space $X$ over ${\bf K}=\bf H$ or ${\bf K}=\bf O$. If
two operators $A$ and $B$ belong to $L_q(X)$, then we denote by $ \{
A, B \} $ the minimal subalgebra over $\bf K$ in $L_q(X)$,
containing $A$ and $B$. The subalgebra $${\cal F}^{\star }:= \{ A\in
L_q(X), \{ A, B \} \mbox{ is quasicommutative for each }B\in {\cal
F} \} $$
 we call supercommutant over $\bf K$ for a family (or a
subalgebra) $\cal F$.
\par Let $\cal N$ be a $C^*$-algebra of operators from $L_q(X)$ on
a Hilbert space $X$ over ${\bf K}=\bf H$ or ${\bf K}=\bf O$,
moreover, $I\in \cal N$, $\cal N$ is closed in the weak operator
topology. If the centre $Z({\cal N})$ consists of ${\bf R}I$, then
we say that $\cal N$ is the factor.
\par For example, $L_q(X)$ is the factor.
\par {\bf 70. Theorem.} {\it If $\cal C$ is a selfadjoint $C^*$-subalgebra
in the algebra $L_q(X)$ for a Hilbert space $X$ over ${\bf K}=\bf H$
or ${\bf K}=\bf O$, $I\in \cal C$, then the closure of a
$C^*$-algebra $\cal C$ relative to the strong and weak operator
topologies coincide with $({\cal C}^{\star })^{\star }$.}
\par {\bf Proof.} At first we mention that ${\cal C}^{\star }$
is closed relative to the weak operator topology, consequently,
$({\cal C}^{\star })^{\star }$ also is closed relative to the weak
operator topology. Due to Theorem 2.22.2 the strong and weak
closures of a $\bf R$-convex subset $Y$ in $L_q(X)$ for a Hilbert
space $X$ over ${\bf K}=\bf H$ or ${\bf K}=\bf O$ coincide. In view
of the fact that $\cal C$ is the $\bf R$-convex set, it follows that
the closures of $\cal C$ relative to the weak and strong operator
topologies coincide. We prove that $({\cal C}^{\star })^{\star }$ is
the strong operator closure of the $C^*$-algebra $\cal C$. In view
of ${\cal C}\subset ({\cal C}^{\star })^{\star }$, also $({\cal
C}^{\star })^{\star }$ is closed relative to the strong and weak
operator topologies, then the closure of the $C^*$-algebra $\cal C$
relative to the strong operator topology is contained in $({\cal
C}^{\star })^{\star }$.
\par Let $T\in ({\cal C}^{\star })^{\star }$ and vectors
$x_1,...,x_n\in \cal C$ be given. We find an operator $T_0\in \cal
C$ such that $ \| ((T-T_0)(x_j) \| < 1$ for each $j=1,...,n$. Let
$X^{\oplus n}:=X\oplus ... \oplus X$ denotes the $n$ times direct
sum of $n$ copies of the Hilbert space $X$ over $\bf K$. For an
operator $F\in L_q(X)$ we denote $F^{\oplus n}:=F\oplus ... \oplus
F$. Then there exists a selfadjoint $C^*$-algebra ${\cal C}^{\oplus
n} := \{ F^{\oplus n}: F\in {\cal C} \} $. Then a $\bf K$-vector
space $[{\cal C}(x^{\oplus n})] := cl_{X^{\oplus n}}({\cal
C}(x^{\oplus n}))$ is invariant relative to ${\cal C}^{\oplus n}$,
since ${\cal C}(x^{\oplus n})$ is a $\bf K$-vector subspace in
$X^{\oplus n}$ due to Corollary 46. In view of Proposition 2.22.20
there exists a graded operator of projection ${\hat E}^{\oplus n}\in
{\cal C}^{\oplus n}$ such that it is the greatest, than others
graded operators of projections in ${\cal C}^{\oplus n}$ and ${\hat
E}^{\oplus n}A^{\oplus n}=A^{\oplus n}{\hat E}^{\oplus n}=A^{\oplus
n}$ for each $A^{\oplus n}\in {\cal C}^{\oplus n}$ such that
$(A^{\oplus n})^{-1}(0)$ and $A^{\oplus n}(X^{\oplus n})$ are $\bf
K$-vector subspaces in $X^{\oplus n}$. Thus, ${\hat E}^{\oplus n}\in
({\cal C}^{\oplus n})^{\star }$, since ${\hat E}_0^{\oplus n}$
commutes with each $A_0^{\oplus n}$ for $A_0\in {\cal C}_0$.
\par We mention that $T^{\oplus n}(x^{\oplus n})=(T(x))^{\oplus n}$
for each $x\in X$. Then $(({\cal C}^{\oplus n})^{\star })_0$
consists of all operators of the form
$F=\bigoplus_{i,j=1}^nF_{i,j}$, $F_{i,j}\in {\cal C}^{\star }_0$,
$F_{i,j}: X^i\to X^j_0$, where $X^i$ is isomorphic with $X$ as a
Hilbert space over $\bf K$ for each $i=1,...,n$, $X=X_0\oplus
X_1i_1\oplus ... \oplus X_mi_m$. Therefore, $((({\cal C}^{\oplus
n})^{\star })^{\star })_0$ consists of all operators
$F=\bigoplus_{i,j=1}^nF_{i,j}$, $F_{i,j}\in ({\cal C}^{\star
})^{\star }_0$, $F_{i,j}: X^i\to X^j_0$, $F_{i,j}=0$ for each $i\ne
j=1,...,n$, $F_{i,i}=F_{1,1}$ for each $i=1,...,n$. Then if $T\in
({\cal C}^{\star })^{\star }$, then $T^{\oplus n}\in (({\cal
C}^{\oplus n})^{\star })^{\star }$. From this it follows, that
${\hat E}^{\oplus n}T^{\oplus n}=T^{\oplus n}{\hat E}^{\oplus n}$,
if $(T^{\oplus n})^{-1}(0)$ and $T^{\oplus n}(X^{\oplus n})$ are
$\bf K$-vector spaces. That is, the range of values of the operator
${\hat E}^{\oplus n}$ is invariant relative to the action of the
operator ${\hat T}^{\oplus n}$. In view of the fact that ${\cal
C}^{\oplus n}(x^{\oplus n})$ is dense in $cl_{X^{\oplus n}}({\cal
C}^{\oplus n}(x^{\oplus n}))$ it follows, that there exists an
operator $T_0\in \cal C$ such that $ \| ((T-T_0)(x_j) \| < 1$ for
each $j=1,...,n$, since ${\cal T}^{\oplus n}(x^{\oplus n})$ belongs
to the region of values of the operator ${\hat E}^{\oplus n}$, also
$x^{\oplus n}=(Ix_1,...,Ix_n)\in {\hat E}^{\oplus n}({\hat
X}^{\oplus n})$.
\par {\bf 71. Proposition.} {\it Each continuous $\bf K$-valued function
on $\bf K$ for ${\bf K}=\bf H$ or ${\bf K}=\bf O$ while an extension
on a bounded subset of normal operators from $L_q(X)$ is continuous
relative to the strong operator topology on a Hilbert space $X$ over
$\bf K$.}
\par {\bf Proof.} In the general case the set of bounded operators
is contained in the ball \\ $B(L_q(X),0,r)$ of the radius $r>0$ in
$L_q(X)$ with the centre at zero. Let $\mbox{ }_0T$ be a normal
operator in this ball, $b>0$, $x_1,...,x_n\in X$. If $ \|
f(T)-f(\mbox{ }_0T) (x/ \| x \| ) \| <b/\| x \| $, then $ \|
f(T)-f(\mbox{ }_0T) (x) \| <b$. That is, it can be supposed that $
\| x \| =1$. Due to the Stone-Weierstrass theorem (see also \S 2.7
\cite{luoystoc}) there exists a polynomial $g$ by the variable $z$,
such that $ \| f-g \|_{C(B({\bf K},0,r),{\bf K})} <b/3$. The
multiplication on  bounded subsets of operators and the operation of
taking adjoints are continuous on a bounded set of normal operators.
Therefore, there exists $y_1,...,y_m\in X$, $v>0$, such that $ \|
(g(T)-g(\mbox{ }_0T))x \| < b/3$ for a marked vector $x\in X$, if $
\| (T-\mbox{ }_0T)y_j \| <b$, $T$ is normal, $ \| T \| \le r$.
Further these arguings can be prolonged for $x=x_1,...,x_m$ by the
way of the mathematical induction. In this case \par $ \|
(f(T)-f(\mbox{ }_0T))x \| \le \| (f(T)-g(T))x \| $ $+ \|
(g(T)-g(\mbox{ }_0T) x \| $ \par  $+ \| (g(\mbox{ }_0T)-f(\mbox{
}_0T) \| +b/3$ $\le 2 \| f-g \| _{C(B({\bf K},0,r),{\bf K})}
+b/3<b$.\\ For the conclusion, that $ \| f(T)-g(T) \| $ and $ \|
f(\mbox{ }_0T)-g(\mbox{ }_0T) \| $ are majorized by the way of $ \|
f-g \|_{C(B({\bf K},0,r),{\bf K})}$ in the preceding inequalities it
can be done the transition to the functional representation of a
quasicommutative algebra generated by $T$ and $T^*$ in accordance
with Theorem 2.22. This representation is accomplished on $sp (T)$,
moreover, to $T$ there corresponds the variable $z$, also to $T^*$
there corresponds $z^*$. Since this representation is the isometry,
then $ \| f(T)-g(T) \| = \| f-g \|_{C(sp(T),{\bf K})}$ and
analogously $ \| f(\mbox{ }_0T) -g (\mbox{ }_0T) \| = \| f-g
\|_{C(sp(T),{\bf K})}$.
\par {\bf 72. Note.} For a selfadjoint operator $T\in L_q(X)$ for
a Hilbert space $X$ over ${\bf K}=\bf H$ or ${\bf K}=\bf O$ and each
$M\in \bf K$ with $|M|=1$ and $Re (M)=0$ there exists the operator
$$U_M(T):= (T-MI)(T+MI)^{-1}.$$ In the particular case of ${\bf K}=\bf C$ and $M=i$
this construction is known as the Cayley transformation.
\par {\bf  73. Proposition.} {\it  The transformation $T\mapsto U_M(T)$ from
\S 72 is continuous relative to the strong operator topology on the
subset of all selfadjoint operators in $L_q(X)$, moreover, the
operator $U_M(T)$ is unitary for each selfadjoint operator $T$ from
$L_q(X)$ and each $M\in \bf K$, $|M|=1$, $Re (M)=0$.}
\par {\bf Proof.} Due to Lemma 2.27 $sp(T)\subset \bf R$,
also due to Theorem 2.22 the algebra, generated by the operator $T$
is isometrically isomorphic with $C(sp(T),{\bf K})$. Then the
operator $U_M(T)$ is correctly defined and belongs to $L_q(X)$. In
view of the alternativity of the algebra $\bf K$ and the functional
representation it follows, that
\par $U_M(T)U_M(T)^*=((T-MI)(T+MI)^{-1})((T-MI)^{-1}(T+MI))$ \par
$=(T-MI)((T^2+I)^{-1}(T+MI))=(T^2+I)(T^2+I)^{-1}=I$ \\
and analogously $U_M(T)^*U_M(T)=I$. Moreover, there are accomplished
relations \par $ (T+MI)((U_M(T)-U_M(A))(A+MI))=2M(T-A)$ \\ for each
selfadjoint operators $A, T$ from $L_q(X)$ (see also the proof of
Lemma 44). Then \par $ \| (U_M(T)-U_M(A))x \| = 2 \|
(T+MI)^{-1}((T-A)(A+MI)^{-1}(x)) \| $ \par $\le 2 \|
(T-A)((A+MI)^{-1}(x)) \| $  $\le 2 \| T-A \| \| x \| $, \\ since
$<x,y>=\sum_{p,q}<x_p,y_q>i_p^*i_q$, also $ \| (T+MI) ^{-1} \| \le
1$ due to the functional representation of the algebra, generated by
the operator $T$.
\par {\bf 74. Theorem.} {\it If $f$ is a continuous $\bf R$-valued function
on $\bf R$ and there exists the limit $\lim_{|x|\to \infty }
f(x)=0$, then $f$ is continuous relative to the strong operator
topology on the set of selfadjoint operators from $L_q(X)$ for a
Hilbert space $X$ over ${\bf K}=\bf H$ or ${\bf K}=\bf O$.}
\par {\bf Proof.} Let $g(z)=f(-(z-1)^{-1}((z+1)M))$ for
$z\ne 1$ and $|z|=1$, where $M\in \bf K$, $|M|=1$, $Re (M)=0$. Put
$g(1)=0$, then the function $g$ is continuous on the set $
\partial B({\bf K},0,1)=S({\bf K},0,1):= \{ z: z\in {\bf K}; |z|=1
\} $, moreover, $g(S({\bf K},0,1)) \subset \bf R$. In view of
\par $(U^*_M(T)+I)(U_M(T)-I)=
U_M(T)-U_M(T^*)=-(U^*_M(T)-I)(U_M(T)+I)$, then \par
$T^*=(M(U^*_M(T)+I))(U^*_M(T)-I)^{-1}=-(U_M(T)-I)^{-1}((U_M(T)+I)M)=T$,\\
since $M^*=-M$. Then $g(U_M(T))=f(T)$ for each selfadjoint operator
$T$. In view of the fact that the function $g$ is continuous on
$S({\bf K},0,1)$ it follows, that it gives a function on a bounded
set of unitary operators, which is continuous relative to the strong
operator topology in accordance with Proposition 71. Since $f$ is
the composition of the transformations $U_M(T)$ and $g$, also due to
Proposition 73 the function $f$ is continuous relative to the strong
operator topology.
\par {\bf 75. Theorem.} {\it If $\cal C$ is a selfadjoint $C^*$-algebra
of operators contained in $L_q(X)$ for a Hilbert space $X$ over
${\bf K}=\bf H$ or ${\bf K}=\bf O$, then for each $T\in B(cl_s({\cal
C}),0,1)$ from the ball of the unit radius with the centre at zero
in $cl_s ({\cal C})$ (the closure of the $C^*$-algebra $\cal C$ in
$L_q(X)$ relative to the strong operator topology) it is
accomplished the inclusion $T\in cl_s (B({\cal C},0,1))$. If an
operator $T$ is selfadjoint in $B(cl_s({\cal C}),0,1)$, then $T$
belongs to the strong operator closure of the set of selfadjoint
operators  from $B({\cal C},0,1)$.}
\par {\bf Proof.} If a selfadjoint operator $T$ belongs
to $cl_s({\cal C})$, also $ \{ F_a : a \in \Upsilon \} $ is a net of
operators converging in $\cal C$ to $T$ relative to the weak
operator topology, then the net $ \{ F_a+F_a^*: a \in \Upsilon \} $
consists of selfadjoint operators and also converges to $T$ relative
to the weak operator topology. In view of the fact that the set
$\cal C$ is convex over $\bf R$, then due to Theorem 2.22.2 $T$
belongs to $cl_s({\cal C})$.
\par Let an operator $T$ be selfadjoint and belongs to $B(cl_s({\cal
C}),0,1)$, also the net of selfadjoint operators $ \{ F_a : a \in
\Upsilon \} $ from $\cal C$ converges to $T$ relative to the strong
operator topology. Take the function $f(t)=t$ on $[-1,1]\subset \bf
R$ and $f(t)=1/t$ for $|t|>1$, where $t\in \bf R$, then
$\lim_{|t|\to \infty }f(t)=0$. The function $f$ generates the
function $g$ from the proof of Theorem 74. Therefore, the net
$f(F_a)$ converges relative to the strong operator topology to
$f(T)$, where each $F_a$ is selfadjoint. Then $f(T)=T$, that is, $T$
belongs to  the strong operator closure of the set of selfadjoint
elements in $B(cl({\cal C},0,1)$, where the closure $cl (A)$ of a
subset $A$ in $L_q(X)$ is taken relative to the topology of the
operator norm. On the other hand, each selfadjoint element from
$B(cl({\cal C},0,1)$ is the limit relative to the topology of the
operator norm, hence also relative to the strong operator topology
of selfadjoint elements from $B({\cal C},0,1)$.
\par {\bf 76. Corollary.} {\it If ${\cal C}\subset L_q(X)$ is a
selfadjoint $C^*$-subalgebra of operators acting on a Hilbert space
$X$ over ${\bf K}=\bf H$ or ${\bf K}=\bf O$, also $T$ is a
nonnegative operator from $B(cl_s({\cal C}),0,1)$, then $T$ belongs
to the strong operator closure for $B({\cal C}^+,0,1)$.}
\par {\bf Proof.} In view of $T\ge 0$ there exists the equality $T=G^2$ for
a selfadjoint operator $G$ from $B(cl_s({\cal C}),0,1)$. Due to
Theorem 75 $G$ belongs to the strong operator closure of selfadjoint
operators from $B({\cal C},0,1)$. From the continuity of the
multiplication on $B({\cal C},0,1)$ relative to the strong operator
topology it follows, that $G^2$ belongs to the strong operator
closure of the family of operators $F^2$ for each selfadjoint
operator from $B({\cal C},0,1)$, but $F^2\in B({\cal C}^+,0,1)$.
\par {\bf 77. Corollary.} {\it  If ${\cal C}\subset L_q(X)$ is a $C^*$-subalgebra
of operators acting on a Hilbert space $X$ over ${\bf K}=\bf H$ or
${\bf K}=\bf O$, $U$ is a unitary operator from  $cl_s({\cal C})$,
then $U$ belongs to the closure relative to the strong operator
topology of the set of unitary operators from $\cal C$.}
\par {\bf Proof.} Due to Theorem 2.33 $U=\exp (MT)$ for
some selfadjoint operator $T$, $M\in \bf K$, $|M|=1$, $Re (M)=0$. In
accordance with Theorem 75 $T$ is the limit relative to the strong
operator topology of the net $\{ F_a: a\in \Upsilon \} $ of
selfadjoint operators from $B({\cal C},0,\| T \| )$, where $\Upsilon
$ is a directed set. From the continuity of the function ${\bf R}\ni
t\mapsto \exp (Mt)\in \bf K$ it follows, that $ \{ \exp (MF_a): a\in
\Upsilon \} $ converges to $\exp (MF)$ relative to the strong
operator topology.
\par {\bf 78. Definition.} A family $\cal G$ from $L_q(X)$ for a Hilbert space
$X$ over ${\bf K}=\bf H$ or ${\bf K}=\bf O$ acts topologically
irreducibly on $X$, when $ \{ 0 \} $ and $X$ are the unique closed
subspaces in $X$ invariant relative to the action of $\cal G$.
\par {\bf 79. Theorem.} {\it If $\cal G$ is a selfadjoint family
contained in $L_q(X)$ for a Hilbert space $X$ over ${\bf K}=\bf H$
or ${\bf K}=\bf O$, then $\cal G$ acts topologically irreducibly on
$X$ if and only if ${\cal G}^{\star }=\bf K$, that is equivalent to
the condition $({\cal G}^{\star })^{\star }=L_q(X)$.}
\par {\bf Proof.} If ${\cal G}^{\star }=\bf K$, then
$({\cal G}^{\star })^{\star }=L_q(X)$. If $({\cal G}^{\star
})^{\star }=L_q(X)$, then $(({\cal G}^{\star })^{\star })^{\star
}=\bf K$. On the other hand, ${\cal G}\subset ({\cal G}^{\star
})^{\star }$, therefore, $(({\cal G}^{\star })^{\star })^{\star
}\subset {\cal G}^{\star }$, but ${\cal G}^{\star }_0$ commutes with
$(({\cal G}^{\star })^{\star })^{\star }_0$, consequently, ${\cal
G}^{\star }\subset (({\cal G}^{\star })^{\star })^{\star }$. Thus,
${\cal G}^{\star }=(({\cal G}^{\star })^{\star })^{\star }$.
\par For a selfadjoint family $\cal G$ of operators from $L_q(X)$
the $C^*$-algebra ${\cal G}^{\star }$ is closed relative to the weak
operator topology and contains the unit operator $I$. Due to Theorem
2.24 ${\cal G}^{\star }=\bf K$ if and only if each graded operator
of the projection from ${\cal G}^{\star }$ is either zero $0$, or
the unit $I$. In view of the fact that the family $\cal G$ is
selfadjoint it follows, that the graded operator of the projection
$\hat E$ belongs to ${\cal G}^{\star }$ if and only if its region of
values ${\hat E}(X)$ is invariant  relative to ${\cal G}$ (see
Proposition 2.22.7). Thus, $\cal G$ acts topologically irreducibly
on $X$ if and only if ${\cal G}^{\star }=\bf K$.
\par {\bf 80. Lemma.} {\it If $ \{ x_1,...,x_n \} $ are orthonormal vectors
in a Hilbert space $X$ over ${\bf K}=\bf H$ or ${\bf K}=\bf O$ with
the scalar product $<*,*>$ from Definition 2.16, $z_1,...,z_n\in
B(X,0,r)$, $0<r<\infty $, then there exists an operator $F\in
L_q(X)$ such that $ \| F \| \le (2n)^{1/2}r$ and $Fx_j=z_j$ for each
$j=1,...,n$. If $Ax_j=z_j$ for some selfadjoint operator $A$, then
$F$ can be chosen selfadjoint.}
\par {\bf Proof.} Let $\hat E$ be a graded
operator of a projection from $X$ on $span_{\bf K} \{ x_1,...,x_n \}
$, which exists due to Proposition 2.22.7. Put $$T(x) :=
\sum_{j=1}^n <{\hat E}(x),x_j>z_j,\mbox{ then }T{\hat E}=T\mbox{
and}$$
$$ \| T(x) \| \le r (\sum_{j=1}^n |<{\hat
E}(x),x>|^2)^{1/2}(\sum_{j=1}^n1)^{1/2} \le n^{1/2} r \| {\hat E}(x)
\| \le n^{1/2} r \| x \| .$$ Thus, $ \| T \| \le n^{1/2} r$.
\par If $Ax_j=z_j$ for some selfadjoint operator
$A$, then ${\hat E}T={\hat E}(A{\hat E})=({\hat E}A){\hat E}$ due to
Proposition 2.22.9, consequently, ${\hat E}T$ is selfadjoint.
Consider an operator $F=T+T^*(I-{\hat E})={\hat E}T+(I-{\hat
E})T+T^*(I-{\hat E})$ and $Fx_j=Tx_j=z_j$. In view of $T(I-{\hat
E})=0=(I-{\hat E})T^*$, then \par $\| F F^* \| = \| F \| ^2 = \|
FF^*+ F^*((I-{\hat E})F) \| \le 2 \| FF^* \| = 2 \| F \| ^2 \le
2nr^2$, \\
since $T(I-{\hat E})=T-T{\hat E}=0$, also $(I-{\hat
E})^2=(I-{\hat E})$.
\par {\bf 81. Theorem.} {\it If a $C^*$-algebra $\cal C$ acts topologically
irreducibly on a Hilbert space $X$ over ${\bf K}=\bf H$ or ${\bf
K}=\bf O$, $ \{ y_1,...,y_n \}$ is a subset of vectors in $X$, $\{
x_1,...,x_n \} $ is a $\bf K$-linearly independent set of vectors in
$X$, then there exists an operator $F\in \cal C$ such that
$Fx_j=y_j$ for each $j=1,...,n$. If $Bx_j=y_j$ for some selfadjoint
operator $B$, then $F$ can be chosen selfadjoint.}
\par {\bf Proof.} Consider the $\bf K$-vector span
$span_{\bf K} \{ x_1,...,x_n \}$, then in this $\bf K$-vector
subspace it can be chosen an orthonormal basis relative to the
scalar product $<*,*>$. Therefore, without restriction of the
generality we suppose, that vectors $x_1,...,x_n$ are orthonormal.
Take operators $\mbox{ }_0G\in L_q(X)$ such that $\mbox{
}_0Gx_j=y_j$ for each $j=1,...,n$. In view of $cl_s ({\cal C}) =
({\cal C}^{\star })^{\star }=L_q(X)$ in accordance with Theorem 79,
then there exists $\mbox{ }_0S\in {\cal C}$, such that \par  $ \|
\mbox{ }_0Gx_j-\mbox{ }_0Sx_j \| = \| \mbox{ }_0Sx_j-y_j \| \le
[2(2n)]^{-1}$. \\ Due to Lemma 80 there exists an operator $\mbox{
}_1G\in L_q(X)$ such that $\mbox{ }_1Gx_j=y_j-\mbox{ }_0Sx_j$ for
each $j$, moreover, $ \| \mbox{ }_1G \| \le 1/2$. The operator
$\mbox{ }_0S$ can be chosen selfadjoint, if the operator $\mbox{
}_0G$ is selfadjoint, since the family of all selfadjoint operators
is dense relative to the strong operator topology in $cl_s({\cal
C})$. In this case in accordance with Lemma 79 the operator $\mbox{
}_1G$ can be taken selfadjoint. Due to Theorem 75 there exists the
operator $\mbox{ }_1S\in \cal C$ such that $ \| \mbox{ }_1S \| \le
1/2$ and $ \| \mbox{ }_1Sx_j- \mbox{ }_1Gx_j \| \le
[4(2n)^{1/2}]^{-1}$, moreover, $\mbox{ }_1S$ is selfadjoint, if
$\mbox{ }_0G$ is selfadjoint.
\par Suppose that it is constructed an operator $\mbox{ }_kG$ such that
$ \| \mbox{ }_kG \| \le 2^{-k}$, \par $\mbox{ }_kGx_j=y_j-\mbox{
}_0Sx_j- \mbox{ }_1Sx_j-...-\mbox{ }_{k-1}Sx_j$, \\  operator
$\mbox{ }_kG$ is selfadjoint, if $\mbox{ }_0G$ is selfadjoint.
Choose $\mbox{ }_kS$ in $\cal C$, such that \par $\| \mbox{ }_kS\|
\le 2^{-k}$, also $ \| \mbox{ }_kSx_j-\mbox{ }_kGx_j \| \le
[2^{k+1}(2n)^{1/2}]^{-1}$ \\  with the help of Theorem 75, where
$\mbox{ }_kS$ is selfadjoint, if $\mbox{ }_kG$ is selfadjoint. Due
to Lemma 80 there exists $\mbox{ }_{k+1}G$ with $ \| \mbox{ }_{k+1}G
\| \le 2^{-(k+1)}S$ and $\mbox{ }_{k+1}Gx_j=y_j- \mbox{
}_0Sx_j-...-\mbox{ }_kSx_j$, where the operator $\mbox{ }_{k+1}G$ is
selfadjoint, if $\mbox{ }_kS$ is selfadjoint. Then the series
$\sum_{k=0}^{\infty } \mbox{ }_kS$ conbverges relative to the
operator norm to the operator $S \in \cal C$, moreover, the operator
$S$ is selfadjoint, if $\mbox{ }_0G$ is selfadjoint, also
$$y_j-Sx_j=y_j-\sum_{k=0}^{\infty }\mbox{ }_kSx_j= \lim_{k\to \infty
}(y_j-\mbox{ }_0Sx_j-...-\mbox{ }_kSx_j) =\lim_k \mbox{
}_{k+1}Gx_j=0.$$
\par {\bf 82. Corollary.} {\it  If a $C^*$-subalgebra $\cal C$ in $L_q(X)$
acts topologically irreducibly on a Hilbert space $X$ over ${\bf
K}=\bf H$ or ${\bf K}=\bf O$, then it acts algebraically irreducibly
on $X$.}
\par Thus, for a $C^*$-subalgebra in $L_q(X)$ it can be not
distinguished the topological and the algebraic irreducibility of
its action on $X$.
\par {\bf 83. Theorem.} {\it If a $C^*$-subalgebra $\cal C$ in $L_q(X)$ acts
topologically irreducibly on a Hilbert space $X$ over ${\bf K}=\bf
H$ or ${\bf K}=\bf O$, also if a unitary operator $V$ on $X$ is such
that $Vx_1=y_1$,...,$Vx_n=y_n$. Then there exists a selfadjoint
operator $S$ on $X$ such that $Ux_1=y_1$,....,$Ux_n=y_n$, where
$U=\exp (MS)$, $M\in \bf K$, $|M|=1$, $Re (M)=0$.}
\par {\bf Proof.} If $x\in X$, $x\ne 0$, then there exists $a\in \bf K$, $a\ne
0$, such that $ax\in X_0$, where $X=X_0\oplus X_1i_1\oplus ...
\oplus X_mi_m$, $X_0,...,X_m$ are pairwise isomorphic Hilbert spaces
over $\bf R$, $ \{ i_0, i_1,..., i_m \} $ is the set of standard
generators of the algebra $\bf K$. Then for an arbitrary set of $\bf
K$-linearly independent vectors $x_1,...,x_m\in X$ there exist
constants $a_1,...,a_n\in \bf K$, $|a_1|=...=|a_n|=1$, for which
$z_1:=a_1x_1,...,z_n:=a_nx_n\in X_0$. But the Hilbert space can be
also considered over the field of real numbers and the restriction
of the scalar product $<x,y>$ from $X$ on $X_0$ gives the scalar
product $(x,y)$ on $X_0$. Then the Schmidt procedure of the
orthogonalization can be applied on $X_0$ to vectors $z_1,...,z_n$.
This gives the orthonormalized system of vectors $q_1,...,q_n\in X$,
for which there exists the operator $S\in L_q(X)$ such that
$S(x_j)=q_j$ for each $j=1,...,n$. In view of the fact that $n$ can
be in the general case taken arbitrary, then by the transfinite
induction in $X$ there exists the orthonormal basis.
\par Without restriction of the generality $\{ x_1,...,x_n \} $
is an orthonormal set of vectors in $X$. The same can be supposed
about $ \{ y_1,...,y_n \} $, since the operator $V$ is unitary,
$<V(x),V(z)>=<x,z>$ for each $x, z \in X$. Denote by $ \{
x_1,...,x_m \} $ and $ \{ y_1,...,y_m \} $ the extension of the
families $ \{ x_1,...,x_n \} $ and $ \{ y_1,...,y_n \} $
respectively up to the orthonormed bases of the space $span_{\bf K}
\{ x_1,...,x_n; y_1,...,y_n \} $. Then it can be prescribed the
operator $U$ on \par $span_{\bf K} \{ x_1,...,x_n; y_1,...,y_n \} $
by the formula $Ux_j=y_j$ for each $j$. \\
It can be taken another
orthonormed basis such that $Ux_j=c_jx_j$ for each $j$ with $
|c_j|=1$ due to Theorems 2.28 and 2.33, moreover, $U = \exp (M T)$
for some selfadjoint operator $T\in L_q(X)$, $M\in \bf K$, $|M|=1$,
$Re (M)=0$. Due to Theorem 81 this operator $T$ can be taken in
$\cal C$ such that $Tx_j = a_jx_j$ for each $j$, where $a_j\in \bf
R$, $\exp (M a_j)=c_j$. In view of the fact that $U$ is the limit of
polynomials by $T$ relative to the norm topology, then $U\in \cal
C$, since $\cal C$ in accordance with Definition 2.16 is the Banach
algebra, that is, complete relative to the topology of the operator
norm.

\par The author is sincerely grateful to Professors Hans de Groote and
Fred van Oystaeyen for discussions of the work and the hospitality
at the Mathematical Departments of the Frankfurt-am-Main and Antwerp
Universities.

\end{document}